\pdfoutput=1
\documentclass[12pt, twoside]{article}
\let\oldphi\varphi \let\varphi\phi \let\phi\oldphi
\let\oldepsilon\varepsilon \let\varepsilon\epsilon \let\epsilon\oldepsilon

\usepackage{amsmath}
\usepackage{imakeidx}
\usepackage{amsfonts}
\usepackage[utf8]{inputenc}
\usepackage{graphicx}
\usepackage[T1]{fontenc}
\usepackage[english]{babel}
\usepackage{graphicx}
\usepackage{fancyhdr}
\usepackage{float}
\usepackage{afterpage}
\usepackage{placeins}
\usepackage{esvect}
\usepackage{multicol}
\usepackage{titlesec}
\usepackage{setspace}
\usepackage[font=small,labelfont=bf]{caption}
\usepackage{url}
\usepackage[left=2cm,right=2cm,top=2.5cm,bottom=2.5cm]{geometry}
\usepackage{musicography}
\usepackage{tikz}
\usepackage{tikz-cd}
\usetikzlibrary{tqft}
\usepackage{dsfont}
\usepackage{amssymb}
\usepackage{enumitem}
\usepackage{makeidx}

\usepackage{amsthm}
\DeclareMathAlphabet{\mymathbb}{U}{BOONDOX-ds}{m}{n}

\usepackage{xstring}
\usepackage{subfigure}
\usepackage{multirow}

\graphicspath{
	{./pics/}
}

\makeindex

\newcommand{\dref}[1]{\textit{definition \ref{#1}}}
\newcommand{\ssecref}[1]{\textit{subsection \ref{#1}}}

\newcommand{\remref}[1]{\textit{remark \ref{#1}}}

\newcommand{\RR}{\ensuremath{\mathbb{R}}}
\newcommand{\CC}{\ensuremath{\mathbb{C}}}
\newcommand{\NN}{\ensuremath{\mathbb{N}}}
\newcommand{\ZZ}{\ensuremath{\mathbb{Z}}}
\newcommand{\ZC}{\ensuremath{\mathcal{Z}}}
\newcommand{\LL}{\ensuremath{\mathcal{L}}}
\newcommand{\FF}{\ensuremath{\mathcal{F}}}
\newcommand{\HH}{\ensuremath{\mathcal{H}}}
\newcommand{\OO}{\ensuremath{\mathcal{O}}}

\newcommand{\gf}{\ensuremath{\mathfrak{g}}}
\newcommand{\kk}{\ensuremath{\mathds{k}}}
\newcommand{\Ci}{\ensuremath{\mathcal{C}^{\infty}}}

\newcommand{\tens}[1]{\ensuremath{\otimes_{#1}}}
\newcommand{\tensor}{\ensuremath{\otimes}} 
\newcommand{\ttensor}{\ensuremath{\widetilde{\tensor}}} 
\newcommand{\pair}[2]{\ensuremath{\langle #1, #2 \rangle}}
\newcommand{\triple}[3]{\ensuremath{\langle #1, #2, #3 \rangle}}
\newcommand{\cop}{\coprod\nolimits}
\newcommand{\ins}{ {\ \cdot \ } }

\newcommand{\ra}{\rightarrow}

\newcommand{\lra}{\longrightarrow}
\newcommand{\Lra}{\Longrightarrow}

\newcommand{\lmap}{\longmapsto}
\newcommand{\atob}[1]{\overset{#1}{\lra}}
\newcommand{\longhookright}{\lhook\joinrel\longrightarrow}
\newcommand{\nat}{\mathrel{\vbox{\offinterlineskip\mathsurround=0pt\ialign{\hfil##\hfil\cr\normalfont\scalebox{2.5}{.}\cr$\lra$\cr}}}}

\newcommand{\cat}[1]{\ensuremath{\mathtt{#1}}}       
\newcommand{\tqft}{\ensuremath{\cat{nTQFT}_\kk}}
\newcommand{\twoft}{\ensuremath{$2$\cat{TQFT}_\kk}}
\newcommand{\comfrob}{\mathtt{comFrob}}
\newcommand{\symfrob}{\mathtt{symFrob}}
\newcommand{\comk}{\ensuremath{\comfrob_\mathds{k}}}
\newcommand{\cob}{\mathtt{cob}}
\newcommand{\cobtwo}{\ensuremath{\mathtt{cob}_2}}
\newcommand{\cobn}[1]{\ensuremath{\cob_{#1}}}
\newcommand{\vect}{\cat{Vect}}
\newcommand{\vectr}{\cat{Vect}_\RR}
\newcommand{\vectk}{\ensuremath{\vect_\mathds{k}}}
\newcommand{\alg}{\cat{Alg}_\kk}
\newcommand{\ract}[1]{\cat{Ract}_{#1}}
\newcommand{\lact}[1]{\cat{Lact}_{#1}}
\newcommand{\rmod}[1]{\cat{rMod}_{#1}}
\newcommand{\lmod}[1]{\cat{lMod}_{#1}}
\newcommand{\finvect}{\cat{FinVect}_\kk}
\newcommand{\mon}{\cat{Mon}}
\newcommand{\smon}{\cat{SymMon}}

\DeclareMathOperator\End{End}
\DeclareMathOperator\Hom{Hom}
\DeclareMathOperator\Aut{Aut}

\DeclareMathOperator\id{id}
\DeclareMathOperator\tr{tr}

\DeclareMathOperator\cod{cod}
\DeclareMathOperator\dom{dom}
\DeclareMathOperator\eval{eval}
\DeclareMathOperator\modulo{mod}
\DeclareMathOperator*{\amp}{\&}
\DeclareMathOperator{\pr}{pr}

\newtheoremstyle{note}
{1.2em}                
{1.2em}                
{}                     
{}                     
{\scshape\bfseries}    
{ }                    
{.5em}                 
{}                     
\theoremstyle{note}
\newtheorem{theorem}{Theorem}[section]
\newtheorem{definition}[theorem]{Definition}
\newtheorem{notation}[theorem]{\textbf{NOTATION}}
\newtheorem{lemma}[theorem]{Lemma}
\newtheorem{corollary}[theorem]{Corollary}
\newtheorem*{proofof}{Proof}
\newtheorem{prop}[theorem]{Proposition}

\newtheorem{example}[theorem]{Example}
\newtheorem{construction}[theorem]{Construction}
\newtheorem{principle}[theorem]{Principle}
\newtheorem{remark}[theorem]{Remark}

\begin{document}

\setlength{\parindent}{0pt}

\begin{titlepage}
  \begin{center}
   \vspace*{1cm}
   \textbf{\large Semester Project}

   \vspace{1.6cm}
   \textbf{\Huge Introduction to 2-dimensional Topological Quantum Field Theory}

   \vspace{5.6cm}
   \textbf{Submitted by:}

   \vspace{0.3cm}
   {\large Leon Menger (né Geiger)}

   \vspace{0.8cm}
   \textbf{Supervisor: Prof. Will J. Merry}

   \vspace{0.8cm}
   Submitted as mandatory part of the programm\\
   \textit{Master of Science ETH}

   \vspace{0.8cm}
   \includegraphics[width=0.7\textwidth]{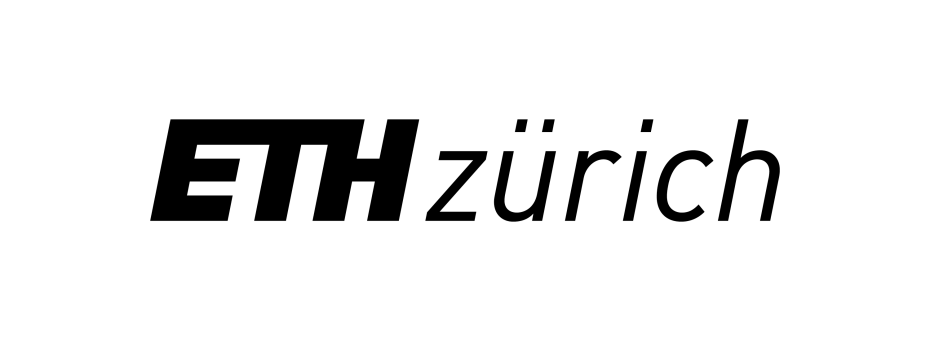}

   \vspace{0.8cm}
   Eidgenössische Technische Hochschule Zürich\\
   Switzerland, 30.06.2020
  \end{center}
\end{titlepage}

\newpage

\tableofcontents\thispagestyle{fancy}
\setcounter{page}{2}
\clearpage
\markboth{\thepage}{Introduction}
\thispagestyle{fancy}

\section{Introduction}
\label{sec:Introduction}

This semester project is a manifest to my current standing in my master studies in physics. Over the course of my undergrad studies my preference for mathematics and especially mathematical physics became more and more apparent. During my studies at ETH I could further delve into this tendency feeling that it was were physics and mathematics are closest that I want to learn more about. One of the topics that sparked my interest during my self-studies was category theory which seemed to provide a sheer endless amount of top-down applications, a new language to generalise and abstract existing structures. For me the most fascinating part about category theory is the natural way in which previously loosely related or seemingly unrelated fields smoothly join hands on this common mathematical background describing their innate structure.\\

Will Merry, whom I got to know as the lecturer of Differential Geometry $I$ during the first Corona-semester, included bits and pieces of category theory in the bonus sections of \cite{Diffgeo_Merry} not without making clear he'd like for this language to be more well-known. After asking him for a semester project relating to category theory he gave me a brief overview of the basic concepts of TQFTs. The topic immediately grabbed my attention since I grew more and more interested in the mathematical background of physical field theories over the course of my physical studies. Using this as a perfect opportunity to both lean more towards mathematical physics and get a better grip on category theory we decided on the format of a script aiming to be as as self-contained as possible and adequate for the workload. One could argue that this got a bit out of hand once I decided to treat $2$-dimensional TQFTs more generally than is standard in literature.\\

TQFT really is a stunning field of mathematical physics. It provides an abstract description of actual physical field theories in a purely mathematical setting and allows to derive certain properties of such theories without ever needing to investigate a particular one. Next to the extensive use of category theory whose fundaments will be provided in \ref{sec:Categorial_Preliminaries} and the subsequent applications of those concepts to cobordisms \ref{sec:Cobordisms} and Frobenius Algebras \ref{sec:Comfrob} the main attraction of this project is found in the section on TQFTs \ref{sec:TQFT}. It provides many interesting detours into and discussions of the physical analogies of the axioms, properties and even two popular examples of TQFTs. Particularly enlightening for these considerations were related discussions of Atiyah \cite{Atiyah_1} and Witten \ref{Witten89a}, the interested reader is wholeheartedly referred to those great texts for further insights.\\

I hope this script can convey even a part of the fun I had while working on it and possibly get a few fellow physicists more interested in the beauty and elegance of category theory. Before the script begins I want to express my sincerest gratitude to Will Merry for his enlightening introduction to the topic and endorsing my transfer to mathematical physics, to Maarten for listening to my excited explanations and the enlightening discussions about physical analogies, to Bryan for the honest reality checks on the dojo matt and to Laura for bearing with my absorbed self.

\newpage

\newpage
\section{Categorial Preliminaries}
\label{sec:Categorial_Preliminaries}

The goal of this chapter is to give a compact and insightful yet not extensive introduction to \emph{Category Theory}. While traces of the terminology and concepts of category theory can be found in most modern mathematics, it will be used in a very concrete manner in this project. Thus the presented preliminaries do not only serve to construct the mathematical setting but directly contribute some of the ideas of \emph{Topological Quantum Field Theory}.\\

We will begin by giving a basic introduction to \emph{categories}, \emph{functors} and \emph{natural transformations} discussing their definitions and highlighting some of their main properties, branching off where needed for upcoming work. Next we introduce the important notion of \emph{monoidal} categories which allow for concatenation of their elements and thus qualify for a description of many interesting categories. After discussing \emph{monoidal functors} and \emph{monodial natural transformations}, we turn towards \emph{braidings}, which can be thought of as "swap" or "twist" operations with respect to the product of a monoidal category, and adapt functors to braidings by defining \emph{braided monoidal functors}. We conclude by giving the full definition of a \emph{symmetric monoidal category} which will turn out to be the fitting setting for the categories involved in \emph{Topological Quantum Field Theory}.\\

This chapter mainly draws from \emph{Categories for the Working Mathematician} by \emph{Saunders Mac Lane} \cite{Categories} which poses an invaluable self-contained introduction to the topic. It will also make use of the bonus material on categories from \emph{Will J. Merry's} Introduction to Differential Geometry \cite{Diffgeo_Merry} due to its compact yet insightful style. When treating monoidal categories, braided monoidal categories and symmetric monoidal categories, it will further make use of the technical review on \emph{Tensor Categories} by \emph{P. Etingof, S. Gelaki, D. Nikshych, and V. Ostrik} \cite{Tensor_Categories}.

\newpage
\subsection{Introduction to Category Theory}
\label{subsec:introduction_categories}

Category Theory might be best described as an effective theory that, in contrary to others, aims at generalizing and unifying different concepts rather than making them more concrete. Serving as a type of "Metatheory" for mathematics, it succeeds at modelling, in a very abstract fashion, collections of objects, maps between them and even mappings between such systems. While certainly not the most concrete piece of mathematics its stunning simplicity allows for a description of a myriad of mathematical structure. As mentioned before, this project will ultimately work with a very concrete application of category theory that showcases both its remarkable flexibility and its inspiring efficiency.\\

While we will soon give the commonly used definition of a category using Set-theoretic notions, we will, for the sake of completeness and deeper insight into the structure of the field, start by giving the more abstract definition seen in \cite{Categories}.

\begin{definition}[Metagraphs]
\label{def:metagraph}
  A \textbf{metagraph}\index{Metagraph} consists of the follwing data: A collection of objects\index{Object} $a,b,c, ...$, a collection of arrows\index{Morphism!arrow} $f,g,h,...$ and two operations
  \begin{itemize}
    \item \textit{Domain}\index{Morphism!domain}, which assigns to each arrow $f$ an object $a = \dom(f)$
    \item \textit{Codomain}\index{Morphism!codomain} (or \textit{Target}), which assigns to each arrow $f$ an object $b = \cod(f)$
  \end{itemize}
\end{definition}

This rather abstract definition of a metagraph can be readily visualized using the suggestive namescheme. We thus indicate an arrow $f$ by an actual arrow and write
$$ f: \dom(f) \lra \cod(f) \quad or \quad a \overset{f}{\lra} b  \quad where \quad a=\dom(f), \ b=\cod(f)$$
While this justifies the name, it also allows for easy representation of finite metagraphs. Expanding the notion of a metagraph with additional structure brings us to the definition of a

\begin{definition}[Metacategory]
\label{def:metacategory}
A \textbf{metacategory}\index{Category!meta-} is a metagraph equipped with the following two additional operations:
  \begin{itemize}
    \item The \textit{Identity} operation, which assigns to each object $a$ an arrow $\id_a = 1_a: a\lra a$ called the \textbf{identity arrow}
    \item The \textit{Composition} operation, which assigns to a pair of arrows $(f,g)$ with $\cod(f) = \dom(g)$ their \textbf{composite arrow} denoted by $g\circ f: \dom(f) \lra \cod(g)$
  \end{itemize}
Those operations are subject to the following two axioms:
  \begin{itemize}
    \item \textbf{Associativity} given three arrows $f,g,h$ such that
    $$ a \atob{f} b \atob{g} c \atob{h} d $$
    implies
    $$ h \circ (g \circ f) = (h \circ g) \circ f $$
    \item The \textbf{Identity Law}: For $a \atob{f} b$ and $b \atob{g} c$ we have
    $$ \id_b \circ f = f \quad \quad g \circ \id_b = g $$
  \end{itemize}
\end{definition}

Now this extension to a metacategory already hints at some common mathematical operations like the composition of smooth functions between manifolds or continuous functions between topological spaces. Moving on to concrete categories, we further expand upon the notion of categories being defined by their arrows. Note that while this definition already uses set theory, it is not yet the promised "Set-version" of a category.

\begin{definition}[Category - Abstract Version]
\label{def:category_abstract}
  Let $A$ be a set of arrows and $\OO$ be a set of objects. Let there be two functions $\cod$ and $\dom$ as before that form a special  metagraph (sometimes called \emph{directed graph})
  \begin{center}
  \begin{tikzcd}
    A\arrow[r, shift left, "dom"]\arrow[r, shift right, "cod"'] &B
  \end{tikzcd}
  \end{center}

  Now in this metagraph the set of composable pairs of arrows in $A$ can be denoted as
  $$ A \times_\OO A := \{ (f,g)| \ f,g \in A \quad s.t. \quad \cod(f) = \dom(g) \}$$
  Now comparing the aforementioned definition of a metacategory, we aim to implement similar properties on the above metagraph. Thus a \textbf{category}\index{Category!abstract} is a graph as constructed above together with two operations
  \begin{align*}
  \OO &\overset{\id}{\lra} A, \quad \quad \quad \ \ c \longmapsto id_c\\
  A\times_\OO A &\lra A, \quad \quad (f,g) \longmapsto g \circ f
  \end{align*}
  again called \textbf{identity} and \textbf{composition} such that
  $$ \dom(\id_a) = a = \cod(\id_a), \quad \dom(g\circ f) = \dom(f), \quad \cod(g\circ f) = \cod(g) $$
  and such that the axioms of \emph{Associativity} and the \emph{Identity Law} hold.
\end{definition}

\begin{notation}
  Let $C$ be a category. In the following we will not explicitely mention the sets $A$ and $\OO$ and instead write $c \in C$ to denote that $c$ is an object of the category $C$ and $f \ in \ C$ to denote that $f$ is an arrow of $C$. Instead of $A$ we write the set of arrows between two objects $a,b \in C$ as
  $$ \Hom(a,b) = \{ f | f \ in \ C, \ \dom(f) = a, \ \cod(f) = b \} $$
  As we form here on work with the more concrete notion of categories in a set theory context, arrows will often be called "morphisms"\index{Morphism} hinting at the generalization of the notion of a homomorphisms.
\end{notation}

Now as promised, we move to the definition directly coming from Set theory \cite{Diffgeo_Merry}. Note that this amounts to defining composition acting on the $\Hom$-sets rather than the set of arrows:

\begin{definition}[Category - Set Version]
  \label{def:category_set}
  A \textbf{category}\index{Category} $\cat{C}$ consists of the following data: A set of objects denoted by $\cat{obj}(\cat{C})$, for each ordered pair of objects $(a,b)$ a set $\Hom(a,b)$ of morphisms\index{Morphism!Hom set} from $a$ to $b$ (defined as above) and a composition rule\index{Morphism!composition} which associates to an ordered triple of objects $(a,b,c)$ a map
  $$ \Hom(a,b) \times \Hom(b,c) \lra \Hom(a,c) $$
  such that for $a\atob{f}b$ and $b \atob{g} c$ the map is denoted as
  $$ (f,g) \longmapsto g\circ f $$
  This data is subject to the following three axioms:
  \begin{itemize}
    \item[i)]   The sets $\Hom$ are pairwise disjoint, thus every $f\in\Hom(a,b)$ has unique domain $a$ and codomain $b$.
    \item[ii)]  Composition is associative if well-defined, thus
    $$ a\atob{f}b \atob{g} c \atob{h} d \quad \Longrightarrow \quad h \circ (g \circ f) = (h \circ g) \circ f $$
    \item[iii)] For each $a \in \cat{obj}(\cat{C})$ there exists a unique morphism $\id_a \in \Hom(A,A)$ called the \textbf{identity morphism}\index{Morphism!identity}. For every $a \atob{f} b$ it satisfies
    $$\id_b \circ f = f = f \circ \id_a$$
  \end{itemize}
\end{definition}

\begin{remark}
\label{rem:small_categories}
  Note that in some books (see \cite{Categories}), one also considers the situation of $\cat{obj}(\cat{C})$ being a class rather than a set which allows for an even more general definition of a category. One then considers \textbf{small categories}\index{Category!small} to be those where both the collection of objects and the collection of arrows turn out to be sets and not classes. Categories are considered \textbf{large categories}\index{Category!large} otherwise. While for the scope of this project, small categories will suffice, some sidenotes will make use of the notion of large categories. A short discussion and directions towards further material using classes can be found in \cite[page 23]{Categories}.
\end{remark}

To further illuminate the definition of a category as defined in \ref{def:category_set}, we give some examples that also show how widely applicable the notion of a category is.

\begin{example}~
  \begin{itemize}
    \item Taking topological spaces as objects and the set of continuous functions $C(X,Y)$ as $\Hom(X,Y)$ for two topological spaces $X,Y$ forms the category \cat{Top}\index{The category of!topological spaces} using the usual composition of functions.

    \item We define the category of groups \cat{Grp}\index{The category of!groups} using groups as objects. Given two groups $G,H$ we define $\Hom(G,H)$ to be the set of all group homomorphisms from $G$ to $H$ and use the usual composition to obtain a full category.

    \item Smooth manifolds form a category \cat{Man}\index{The category of!smooth manifolds}: To see this, take objects smooth manifolds and, using the usual composition of smooth maps, define for any two smooth manifolds $M,N$ the set $\Hom(M,N)$ to be $\Ci(M,N)$.

    \item The most important example for later will be the category of vector spaces over a field $\mathds{k}$ denoted by $\vectk$\index{The category of!vector spaces}. Its objects are $\mathds{k}$-vector spaces and its morphisms are the $\mathds{k}$-linear maps. In many ways, $\vectk$ is a special category: One of the most interesting properties is that since for $V,W \in \vectk$ the set $\Hom(V,W)$ naturally is a vector space, it is also an element of the category. As will be apparent after \ssecref{subsec:Monoidal_categories}, $\vectk$ is also an example of a \textbf{monoidal} category. Due to its importance in TQFT, we will continue to investigate its properties.
  \end{itemize}
\end{example}

Looking at the above examples the question arises, how the arrows of a category can capture concepts like \emph{bijections}, \emph{diffeomorphisms} or \emph{invertible linear transformations}. Each of these embodied the idea of an \emph{isomorphism} in their respective category, allowing for a two-sided inverse:

\begin{definition}[Isomorphisms]
\label{def:object_iso}
  Let $\cat{C}$ be a category and $f \colon a \lra b$ an arrow in $\cat{C}$. We call $f$ an \textbf{isomorphism of objects}\index{Object!isomorphism} or just \textbf{isomorphism} if there exists an arrow $g \colon b \lra a$ such that
  $$ f \circ g = \id_a \quad \quad \text{and} \quad \quad g \circ f = \id_b $$
  Given two objects $a,b \in \cat{C}$ we call them \textbf{isomorphic} if there exists an isomorphism between them.
\end{definition}

For further use, we define special types of categories which will turn out to be deeply tied to groups:

\begin{definition}
\label{def:monoids}
(Monoids and Groups)
\begin{itemize}
  \item A \textbf{discrete category}\index{Category!discrete} is a category which only has identity arrows. Thus it is uniquely characterized by its set of objects, hence discrete categories are sets.

  \item A category with only one object $a$ is called a \textbf{monoid}\index{Category!monoid}. Thus it is uniquely characterized by its arrows, the identity arrow and its composition law. Since every arrow has domain and codomain $a$, any two arrows have a composite. Thus a monoid can be described as a set $M$ with a binary associative operation $M \times M \ra M$ and a unit. Thus a monoid is a a \textbf{semigroup} with unit. Note that for any category $\cat{C}$ and any object $c\in \cat{C}$ the set of \textbf{endomorphisms} $\End(c) = \Hom(c,c)$ forms a monoid under composition.

  \item A monoid where every arrow has a two-sided inverse under the given composition law expands to and is thus called a \textbf{group}\index{Category!group}.
\end{itemize}
\end{definition}

When talking about categories as collections of objects and the arrows between them, the next natural step would be to investigate arrows between categories. This leads us to the notion of a \emph{functor}, a morphism between categories. As such it will need to send objects in the domain category to objects in the codomain category and likewise with functions. While it is possible to define a functor purely using arrows, we follow the definition of \cite{Categories} which uses objects and arrows alike.

\begin{definition}[Functors]
\label{def:functors}
Let $\cat{A}, \cat{B}$ be two categories. A \textbf{covariant functor}\index{Functor} of categories $\cat{F}$ is a morphism of categories $ \cat{A} \atob{\cat{F}} \cat{B}$ consisting of two functions, both denoted by $\cat{F}$ in an abuse of notation:
\begin{itemize}
  \item[1.] An \textbf{arrow function} assigning to every morphism $f:a \ra a^\prime$ in $\cat{A}$ a morphism $\cat{Ff}: \cat{F}a \ra \cat{F}a^\prime$ in $\cat{B}$ such that for $a \atob{f}b \atob{g} c$
  $$ \cat{F}(g\circ f) = \cat{F}(g) \circ \cat{F}(f) $$

  \item[2.] An \textbf{object function} assigning to every $a\in \cat{A}$ an object $\cat{F}a \in B$ such that
  $$ \cat{F}(\id_a) = \id_{\cat{F}(\cat{A})} \quad \forall a \in \cat{A} $$
\end{itemize}
\end{definition}

Note that we call such a functor "covariant"\index{Functor!covariant} because it preserves the direction of the morphisms as can be seen from its action on the composition of morphisms. A "contravariant"\index{Functor!contravariant} functor reverses the order. In the following material, a functor refers to a covariant functor whereas the contravariant version will be explicitely denoted as such.

\begin{example}
  Again turning towards the category $\vectk$, we find that the Tensor product turns out to be a covariant functor: The product $\otimes : (\vectk, \vectk) \ra \vectk$ assigns $(V,W) \mapsto V \otimes W$ and functions alike. Note that this product is associative up to isomorphism which will later enable us to define the category $\vectk$ as a monoidal category.
\end{example}

An interesting class of functors are \textbf{forgetful}\index{Functor!forgetful} functors which, as the name implies, forget part of the structure of the domain category. An easy example is the forgetful functor $\cat{Top} \lra \cat{Sets}$\index{The category of!small sets} which assigns to a topological space $X$ the respective set it is defined on and to a continuous function the same function viewed as a function between sets. A similar forgetful functor can be constructed as $\cat{Man} \lra \cat{Top}$. It "forgets" the manifold structure of the underlying topological space and sends the smooth maps of $\cat{Man}$ to the same maps viewed as continous maps of topological spaces. We will come back to forgetful functors later on when constructing TQFTs \ref{sec:TQFT}.\\

Like morphisms can be composed to provide new morphisms, we can inspect the composition of functors\index{Functor!composition}. Given the categories $\cat{A},\cat{B}$ and $\cat{C}$ and two functors
$$\cat{A} \atob{\cat{F}} \cat{B} \atob{\cat{G}} \cat{C}$$
one can inspect the composition of two functor using the induced object and arrow function. On objects $a\in\cat{A}$ and morphisms in $\cat{A}$ they are defined as the maps
$$ a \lmap \cat{G}(\cat{F}(a)) \quad \quad \quad \quad f \lmap \cat{G}(\cat{F}(f)) $$
Thus they define a \textbf{composite functor} denoted by $\cat{G}\circ \cat{F}: \cat{A} \lra \cat{C}$ called the composition of $\cat{G}$ with $\cat{F}$. The composition functor for contravariant functors is built accordingly. To no surprise one can for every category $\cat{A}$ construct the according identity functor $\cat{I}_\cat{A}: \cat{A} \lra \cat{A}$ that acts as an identity for functor composition. Comparing the data at hand with the definition of a category brings us to the notion of a \textit{"category of categories"}. And indeed, using small categories as objects (compare \remref{rem:small_categories}) one can form the category $\cat{Cat}$\index{The category of!categories} with functors as arrows between its objects. It uses the composition constructed above to compose its arrows.\\

Next, we turn towards special types of functors that play special roles when comparing two categories. First and foremost, we want to have a tool at hand to express that "two categories are the same", namely a well-defined form of isomorphisms.

\begin{definition}[Isomorphism of categories]
\label{def:isomorphism_categories}
  A functor between two categories $\cat{A} \atob{\cat{F}} \cat{B}$ is called an \textbf{isomorphism of categories}\index{Category!isomorphism}\index{Functor!isomorphism} if and only if there exists a functor $\cat{B} \atob{\cat{G}} \cat{A}$ such that the possible composition functors satisfy
  $$ \cat{G} \circ \cat{F} = \cat{I}_\cat{A}  \quad \quad and \quad \quad \cat{F} \circ \cat{G} = \cat{I}_\cat{B}$$
  We call $\cat{G}$ the \textbf{two-sided inverse} of $\cat{F}$ and often denote it as $\cat{G} = \cat{F}^{-1}$. Note that the composition of two isomorphisms is again an isomorphism.
\end{definition}

We will present two very important notions \emph{weaker} than that of an isomorphism of categories which will also find their place in later chapters. They can be read as more general but weaker forms of \emph{surjectivity} and \emph{injectivity}:

\begin{definition}
\label{def:full_faithful}
(Full and faithful functors)
\begin{itemize}
  \item A functor $\cat{A} \atob{\cat{F}} \cat{B}$ is a \textbf{full functor}\index{Functor!full} if for every pair of objects $a,b \in \cat{A}$ and every morphism $g: \cat{F}a \lra \cat{F} a^\prime$ of the codomain category $\cat{B}$ there exists a morphism $f:a \lra a^\prime$ in $\cat{A}$ such that $\cat{F} f = g$.

  \item A functor $\cat{A} \atob{\cat{F}} \cat{B}$ is a \textbf{faithful functor}\index{Functor!faithful} if for every pair of objects $a,b \in \cat{A}$ and any pair of morphisms $f_1, f_2: a \lra b$ the equality $\cat{F}f_1 = \cat{F}f_2$ implies $f_1 = f_2$.
\end{itemize}
\end{definition}

\begin{remark}
  Full functors portray a weaker form of surjectivity because not all objects in $\cat{B}$ need to be of the form $\cat{F}a$ for some $a \in \cat{A}$, thus there may exist morphisms not covered by a full functor. The composition of full functors is again a full functor.\\
  Meanwhile faithful functors, sometimes called \emph{"embeddings"}, represent a weaker form of injectivity: Since two morphisms in the domain category may be mapped to the same morphism in the codomain category, a faithful functor need not be injective on morphisms. Again, the composition of faithful functors is again a faithful functor.
\end{remark}

With a theory at hand describing the very essence of mappings between mathematical objects and their collections, may they host sets, topological spaces, groups, manifolds or something completely different, it is only natural to ask for an elevated concept to that of a subset, subspace, subgroup or submanifold. The notion of a subcategory is in many ways just what one one would expect:

\begin{definition}[Subcategories]
\label{def:subcategory}
Given a category $\cat{B}$, a \textbf{subcategory}\index{Category!subcategory} $\cat{A} \subseteq \cat{B}$ is a collection of objects and morphisms of $\cat{B}$ such that for each object $a$ the respective identity morphism, for each morphism $f$ both $\dom(f), \cod(f) \in \cat{A}$ and for each pair of composable arrows their composite arrow also lies in $\cat{A}$.
\end{definition}

Since these properties make $\cat{A}$ a category, we can use the inclusion functor $\cat{A} \longhookright \cat{B}$ that sends both objects and morphisms of $\cat{A}$ to themselves in $\cat{B}$ to further classify subcategories using \dref{def:full_faithful}.

\begin{definition}
\label{def:full_subcategory}
  Let $\cat{A}$ be a subcategory of $\cat{B}$. We call $\cat{A}$ a \textbf{full subcategory}\index{Category!subcategory!full} of $\cat{B}$, if the inclusion functor $\cat{A} \longhookright \cat{B}$ is a full functor. Thus a full subcategory is uniquely defined by its set of objects, since it simply contains every arrows between these objects that also lies in the base category.
\end{definition}

When closely inspecting full subcategories one might ask, if they can "capture" the full structure of the ambient category. Indeed one can investigate such structures which in many cases simplifies the investigation of the full category.

\begin{definition}[Skeleton of a Category]
\label{def:skeleton}
  Given a category $\cat{C}$ a \textbf{skeleton}\index{Skeleton of a category} $\cat{S}$ of $\cat{C}$ is a full subcategory such that every object $c \in \cat{C}$ is isomorphic \ref{def:object_iso} to exactly one object in $\cat{S}$. Thus the inclusion functor $\cat{S} \hookrightarrow \cat{C}$ is an equivalence of categories. Further we call a category \textbf{skeletal} when it is its own skeleton.
\end{definition}

A similar concept can be defined for morphisms of a category. Namely given the full set of morphisms of a category, which subset can form any arrow using the usual composition of arrows.

\begin{definition}[Generating Set of a Category]
\label{def:generating_set}
  For a category $\cat{C}$, a \textbf{generating set} or \textbf{set of generators}\index{Generators of a category} is a set $\cat{A}$ of arrows such that any arrow of $\cat{C}$ can be obtained by composing arrows of $\cat{A}$.
\end{definition}

Both the \emph{skeleton} and the \emph{generating} set of a category will become utterly helpful when dealing with conrete categories: Identifying the skeleton and then finding a set of generators for a tremendously simplified treatment of categories since only a subset of objects and arrows has to be considered. While it is in general not trivial to find a skeleton or prove that a certain set of arrows is indeed a set of generators, we will encounter cases where both can be done with manageable effort.\\

The last fundamental concept of category theory we will introduce in this section is that of a \emph{natural transformation}. A natural transformation transforms between functors thus posing a tool to compare them. One might think about natural transformations being "functors between functors". And indeed, in the category $\cat{Cat}$ of small categories, natural transformations act like functors do in our usual categories. Note that one can define higher-order categories and regard $\cat{Cat}$ as a $2$-category where natural transformations pose the respective $2$-morphisms. Additional material on this interpretation of $\cat{Cat}$ can be found in \emph{chapter II.5} of \cite{Categories}.

\begin{definition}[Natural Transformations]
\label{def:natural_transformation}
  Given two categories $\cat{A}, \cat{B}$ and two functors $\cat{F},\cat{G}:\cat{A} \lra \cat{B}$ between them, a \textbf{natural transformation}\index{Natural transformation} $\tau: \cat{F} \nat \cat{G}$ is a function that assigns to each object $a \in \cat{A}$ an arrow $\tau_a = \tau a: \cat{F}a \lra \cat{G} a $ and, for each arrow $f:a \lra a^\prime$ in $\cat{A}$ makes the following diagram commute:
  \begin{center}
  \begin{tikzcd}[sep=large]
    a \arrow[d,"f"] & \cat{F} a \arrow[r,"\tau_a"] \arrow[d,"\cat{F(f)}"']       & \cat{G}a \arrow[d,"\cat{G}(f)"] \\
    a^\prime        & \cat{F} a^\prime \arrow[r,"\tau_{a^\prime}"']              & \cat{G}a^\prime
  \end{tikzcd}
  \end{center}

\end{definition}

Now a very important type of natural transformation is that extending the notion of isomorphisms of categories (see \dref{def:isomorphism_categories}).

\begin{definition}[Natural isomorphisms]
\label{def:natural_isomorphism}
  Given two categories $\cat{A}, \cat{B}$, two functors $\cat{F},\cat{G}:\cat{A} \lra \cat{B}$ between them and a natural transformation $\tau: \cat{F} \nat \cat{G}$, we call $\tau$ a \textbf{natural isomorphism}\index{Natural transformation!natural isomorphism}, if for every $a \in \cat{A}$ the respective map $\tau_a$ is invertible in $\cat{B}$.
  We denote this by $\tau: \cat{A} \cong \cat{B}$ and observe that the inverses $(\tau_a)^{-1}$ define a natural isomorphism $\tau^{-1}:\cat{B} \nat \cat{A}$ or rather $\tau^{-1}:\cat{B} \cong \cat{A}$.
\end{definition}

In many ways, natural transformations make precise statements about certain "natural" or "canonical" isomorphisms arising both in physics and pure mathematics. That is when saying that two objects are naturally isomorphic, there exists a natural isomorphism between the respective categories. Taking $\vectk$ as our working example, we show that this is indeed the case for the double dual functor, following the proof of \cite{Diffgeo_Merry}:

\begin{lemma}
\label{lemma:double_dual}
  Let $V$ be a finite dimensional vector space over the field $\mathds{k}$ and denote the category of finite dimensional vector spaces over $\kk$ by $\finvect$. Then the double dual $V^{**}$ is naturally isomorphic to $V$.
\begin{proof}
  Let $\cat{F}: \vect \lra \finvect$ be the functor that acts on objects $V\in \finvect$ as
  $$ \cat{F}(V) := V^{**} = \Hom(\Hom(V,\kk),\kk) $$
  For an arrow of $\finvect$, that is a linear map $l: V \lra W$ for $V,W \in \finvect$, the arrow functor of $\cat{F}$ acts via
  \begin{align*}
    \cat{F}(l): \cat{F}(V) &\lra \cat{F}(W),\\
    l &\lra l^{**} (:V^{**}\lra W^{**})
  \end{align*}
  with $l^{**}$ defined for $\alpha \in V^{**}$ and $\beta \in W^{*}$ as
  $$ l^{**}(\alpha)(\beta) = \alpha(\beta \circ l) $$
  Thus we have established a valid functor representing the double dual acting on a vector space. To conclude the proof, we need to build a natural isomorphism between this functor and the evident identity functor. One evident candidate is the function that assigns to each $V \in \finvect$ the arrow
  $$ \eval_V: \id(V) \lra \cat{F}(V) = V^{**}, \quad \quad \eval_V (v)(\zeta) := \zeta(v) \quad \text{for} \ \ v \in V, \zeta \in V^* $$
  Now all that's left to show is that for any $l:V\lra W$ the following diagram commutes:
  \begin{center}
  \begin{tikzcd}[sep=large]
    V \arrow[d,"l"] & V \arrow[r,"\eval_V"] \arrow[d,"l"']       & V^{**} \arrow[d,"l^{**}"] \\
    W        & W \arrow[r,"\eval_W"']              & W^{**}
  \end{tikzcd}
  \end{center}
  Thus let $v\in V$ and $\zeta \in W^*$. We indeed obtain:
  \begin{align*}
    l^{**}(\eval_V(v))(\zeta) &= \eval_V(v)\left( \zeta \circ l \right) = (\zeta \circ l) (v)\\
    &= \zeta(lv) = \eval_V(lv)(\zeta)
  \end{align*}
  which readily completes the proof.
\end{proof}
\end{lemma}

While we already introduced the notion of an \emph{isomorphism of categories} \ref{def:isomorphism_categories}, we are now in a position to comfortably introduce the \emph{equivalence of categories} without having to supress natural transformations.

\begin{definition}[Equivalence of categories]
\label{def:equivalence_categories}
  Let $\cat{A}, \cat{B}$ be categories and $\cat{A} \atob{\cat{F}} \cat{B}$ a functor between them. $\cat{F}$ is an \textbf{equivalence of categories}\index{Functor!equivalence of categories}\index{Category!equivalence} (thus $\cat{A}$ and $\cat{B}$ are equivalent), if there exists a functor $\cat{G}: \cat{B} \lra \cat{A}$ such that the two evident compositions
  $$ \cat{F} \circ \cat{G}: \cat{A} \lra \cat{A} \quad \quad and \quad \quad \cat{G} \circ \cat{F}: \cat{B} \lra \cat{B}$$
  are are naturally isomorphic to the respective identity functor $\cat{I}_\cat{A}$ and $\cat{I}_\cat{B}$. Note that this also makes $\cat{G}$ an equivalence of categories.
\end{definition}

Looking back at full and faithful functors, it seems only natural to make precise the idea of a bijective functor being an equivalence of categories. We thus aim to prove

\begin{theorem}
\label{theo:equivalence}
  Let $\cat{A} \atob{\cat{F}} \cat{B}$ be a functor between categories. Then $\cat{F}$ is an equivalence of categories if and only if $\cat{F}$ is full, faithful and every object in $\cat{B}$ is isomorphic to $\cat{F}a$ for some $a \in \cat{A}$.
\begin{proof}
  Assume $\cat{F}$ is an equivalence of categories. Let $\cat{G}$ be an "inverse" functor as in \ref{def:equivalence_categories} and $\tau: \cat{G} \circ \cat{F} \cong \cat{I}_\cat{A}$ and $\tau^{-1}: \cat{F} \circ \cat{G} \cong \cat{I}_\cat{B}$ the respective natural isomorphisms. Now since for every element $b\in \cat{B}$, $b \cong \cat{F}(\cat{G}(b))$ holds true, we see that $b \cong \cat{F}(a)$ for $a = \cat{G}(b) \in \cat{A}$.
  Further using $\tau$, we take any arrow $f:a \lra a^\prime$ in $\cat{A}$ to obtain the commutative diagram

  \begin{center}
  \begin{tikzcd}[sep=large]
    a \arrow[d,"f"] & \cat{G}\circ \cat{F}(a) \arrow[r,"\tau_a"] \arrow[d,"\cat{G}\cat{F}(f)"']       & a \arrow[d,"f"] \\
    a^\prime        & \cat{G}\circ \cat{F}(a^\prime)  \arrow[r,"\tau_{a^\prime}"']              & a^\prime
  \end{tikzcd}
  \end{center}

  Thus we obtain $f = \tau_{a^\prime} \circ \cat{GF}(f) \circ \tau_a^{-1}$. Thus, if for some arrow $g: a \lra a^\prime$ the equation $\cat{F}(f) = \cat{F}(g)$ holds true, the arrows must agree which proves that $\cat{F}$ is faithful. Completely symmetrical, we can show that $\cat{G}$ is faithful. Now to show that $\cat{F}$ is full, we consider any arrow in $\cat{B}$ of the form $h: \cat{F}a \lra \cat{F}a^\prime$ for some $a, a^\prime \in \cat{A}$. Now setting $f = \tau_{a^\prime} \circ \cat{G}(h) \circ \tau_a^{-1}$ leaves us with $\cat{GF}f = \cat{G}(h)$ and since we know that $\cat{G}$ is faithful, $\cat{F}f = h$ which makes \cat{F} a full functor.\\

  The converse statement is a bit easier: Since the object function of $\cat{F}$ is by construction surjective and a full faithful functor is (up to isomorphism) injective on objects, it is bijective on objects. Now since $\cat{F}$ is full and faithful, its arrow function is bijective on all arrows of the form $g: \cat{F}a \lra \cat{F}a^\prime$ for some $a,a^\prime \in A$. Since all arrows are necessarily of that form, $\cat{F}$ has a two-sided inverse and thus, by rule of \ref{def:isomorphism_categories} and \ref{def:equivalence_categories}, is an equivalence of categories.
\end{proof}
\end{theorem}

To end this subchapter, note that one of the founding fathers of category theory, \emph{Saunders Mac Lane}, stated that \emph{"[...] 'category' has been defined in order to be able to define 'functor' and 'functor' has been defined in order to be able to define 'natural transformation'."} \cite{Quote_MacLane}. This will turn out to be true for most of our applications of category theory. After finding certain assignments that naturally arise in physics to be functors, thus applying the framework of category theory to a physical setting, the main point of interest will indeed be the comparison of such functors and the underlying categories using natural transformations. That being said, the most important type of natural transformations will be natural isomorphisms which allow for a precise formulation of "canonical" identification.

\newpage

\subsection{Monoidal Categories and Functors}
\label{subsec:Monoidal_categories}

It is undeniable that vector spaces and their morphisms play an important role in theoretical physics. As we will see in the subsequent chapters, the thus defined category $\vectk$ is a natural candidate ingredient for the notion of a \emph{Topological Quantum Field Theory}. To make $\vectk$ tangible in the context of category theory, this subchapter will generalize many of the properties it inherits from its objects. We will see that many exterior products like the tensor product or the direct product are guises of only a few underlying categorial concepts. After this subchapter we will be able to conretely identify $\vectk$ as a \emph{monoidal category} and understand morphisms within and natural transformations between such categories.\\

We start by making precise the often encountered notion of commonly enctountered expressions of the form $U\times V$ for two sets $U, V \in \cat{Set}$ or $V \times W$ for two vector spaces $V,W \in \vect$. This direct product is widely used in mathematical and physical notation and while its use is clear by informal arguments, it stems from a more precise categorial concept:

\begin{definition}[Product categories]
\label{def:product_categories}
  Let $\cat{A}, \cat{B}$ be two categories. The \textbf{product category}\index{Category!product category} denoted by $\cat{A}\times \cat{B}$ is the category whose objects are pairs of the form $\pair{a}{b}$ for $a \in \cat{A}$ and $b\in \cat{B}$. An arrow $\pair{a}{b} \lra \pair{a^\prime}{b^\prime}$ of $\cat{A} \times \cat{B}$ is a pair $\pair{f}{g}$ of arrows $f:a\lra a^\prime$ and $g:b \lra b^\prime$, the composite of two such arrows is defined by
  $$ \pair{f^\prime}{g^\prime} \circ \pair{f}{g} = \pair{f^\prime \circ f}{g^\prime \circ g} $$
  The two evident functors defined on objects and arrows by
  \begin{align*}
    \cat{A} \overset{\cat{P}}{\longleftarrow} \cat{A} &\times \cat{B} \atob{\cat{Q}} \cat{B}\\
    \cat{P}\pair{f}{g} = f,& \quad \quad \cat{Q}\pair{f}{g} = g
  \end{align*}
  are called \textbf{projection functors}\index{Functor!projection}.
\end{definition}

\begin{remark}
  Note that for any category $\cat{C}$ equipped with two functors
  $$ \cat{A} \overset{\cat{L}}{\longleftarrow} \cat{C} \atob{\cat{R}} \cat{B} $$
  there exists a unique functor $\cat{F}: \cat{C} \lra \cat{A} \times \cat{B}$ defined on arrows $h$ in $\cat{C}$ as $\cat{F}h = \pair{\cat{L}h}{\cat{R}h}$ thus satisfying $\cat{P}\cat{F} = \cat{L}, \cat{Q}\cat{F} = \cat{R}$. This can be visualized in the following commuting diagram:
  \begin{center}
  \begin{tikzcd}[sep=large]
    & \arrow[dl, "\cat{L}"'] \cat{C} \arrow[d, dashed, "\cat{F}"]  \arrow[dr, "\cat{R}"]&\\
    \cat{A}  &  \arrow[l, "\cat{P}"] \cat{A} \times \cat{B} \arrow[r, "\cat{Q}"'] &  \cat{B}
  \end{tikzcd}
  \end{center}
  Properties like this can be generalized to the notion of \textbf{general properties}\index{Morphism!universal} or \textbf{general arrows}. The pair of projection functors is universal among pairs of functors to $\cat{A}$ and $\cat{B}$.
\end{remark}

In the above remark we already saw a functor mapping a category to a product category. Socalled \textbf{bifunctors}\index{Functor!bi-} that map a product category to a usual category will be indispensable for the discussion of \emph{monoidal categories}. One can interpret a bifunctor as a functor whose object and arrow function depend on two instances respectively. Taking the usual cartesian product $X \times Y$ of two sets $X,Y \in \cat{Set}$ as an example, we immediately see that it denotes nothing else but the object function of the evident bifunctor $\cat{Set}\times \cat{Set} \lra \cat{Set}$.\\

Having a rather precise notion of product categories, we aim to show that the product $\times$ can be understood as a bifunctor itself. Note that while we discuss it as a functor on the product category
$$ \times: \cat{Cat} \times \cat{Cat} \lra \cat{Cat} $$
there are very similar functors $\cat{Set} \times \cat{Set} \lra \cat{Set}$, $\cat{Top} \times \cat{Top} \lra \cat{Top}$, $\cat{Man} \times \cat{Man} \lra \cat{Man}$ etc.which, in an abuse of notation, are often collectively refered to as \emph{the product functor} $\times$. To identify the product functor, take two functors $\cat{F}: \cat{A} \lra \cat{B}$ and $\cat{G}: \cat{C} \lra \cat{D}$ and note that the respective "product functor" $\cat{F} \times \cat{G}: \cat{A} \times \cat{C} \lra \cat{B} \times \cat{D}$ can be uniquely defined on arrows and objects as
$$ (\cat{F} \times \cat{G})\pair{f}{g} = \pair{\cat{F}f}{\cat{G}g}, \quad \quad (\cat{F} \times \cat{G})\pair{a}{c} = \pair{\cat{F}a}{\cat{G}c} $$
Thus we can identify the product $\times$ as a bifunctor carrying the previously defined product of categories as its object function and the above defined product functors as its arrow function. Thus it constitutes a bifunctor $\cat{Set} \times \cat{Set} \lra \cat{Set}$. Thus having a clear understanding of the product of categories, we can use $\times$ unambigously which enables the definition of many additional structures on categories.\\

As the title of this subchapter suggests, our main goal is to approach \emph{monoidal categories} which can be thought of as categories that carry an internal product product between their objects and arrows that is associative up to isomorphism. In the construction of such categories, we will make use of bifunctors and product categories.

\begin{remark}
  Thinking back to \emph{Monoids} \ref{def:monoids}, we aim to generalize the idea of a category equipped with an associative binary operation and a unit to categories with more than one object. This concept is thus parallel to that of a \emph{monoid} in set theory. While this roughly explains the name of \emph{monoidal categories}, another beautiful explanation comes yet again from higher category theory (again, see \textit{chapter II.5} of \ref{Categories}). Here a monoidal 1-category can be considered as a 2-category with only one object satisfying the analogous properties of a monoid for 2-categories.
\end{remark}

Without further ado, we turn towards the definition of monoidal categories, keeping in mind the structure we want to generalize.

\begin{definition}[Monoidal categories]
\label{def:monoidal_category}
A \textbf{monoidal category}\index{Monoidal category} $\cat{M} = \pair{\cat{M}}{\otimes, e, \alpha, \lambda, \rho}$ consists of a category $\cat{M}$, a bifunctor $\otimes: \cat{M} \times \cat{M} \lra \cat{M}$, an object $e\in \cat{M}$ acting as a unit of the product and three natural isomorphisms $\alpha, \lambda, \rho$ on $\cat{M}$ subject to the following conditions:
\begin{itemize}
  \item For any $a\in\cat{A}$ the natural isomorphisms $\lambda,\rho$ satisfy
  \begin{align*}
    \lambda_a : a \otimes e \cong a, \quad \quad \rho_a: e \otimes a \cong a, \quad \quad \lambda_e = \rho_e : e\otimes e \lra e
  \end{align*}

  \item For any $a,b,c,d \in\cat{M}$ the following commuting diagram, called the \textbf{pentagon identity} or \textbf{pentagon equation}\index{Monoidal category!pentagon equation}, arises from $\alpha_{a,b,c} : a\otimes (b \otimes c) \cong (a\otimes b) \otimes c$:
  \begin{center}
  \begin{tikzcd}[sep=huge]
     & (a\otimes b) \otimes (c \otimes d) \arrow[dr, "\alpha_{a \otimes b, c, d}"] & \\
    a\otimes (b \otimes (c \otimes d)) \arrow[d, "\id_a \otimes \alpha_{b, c, d}"'] \arrow[ur, , "\alpha_{a, b, c \otimes  d}"] &   & ((a\otimes b) \otimes c) \otimes d \\
     a\otimes ((b \otimes c) \otimes d) \arrow[rr, , "\alpha_{a, b  \otimes c, d}"] &  & (a\otimes (b \otimes c)) \otimes d \arrow[u, , "\alpha_{a, b, c} \otimes \id_d"']
  \end{tikzcd}
  \end{center}

  \item For any $a,b \in \cat{A}$ the following diagram, called the \textbf{triangle identity}\index{Monoidal category!triangle identity}, commutes:
  \begin{center}
  \begin{tikzcd}[sep=large]
      & a \otimes b & \\
    a \otimes (e \otimes b) \arrow[ur,"\id_a \otimes \rho_b"] \arrow[rr, "\alpha_{a,e,b}"] &             & (a \otimes e) \otimes b \arrow[ul, "\lambda_a \otimes \id_b"']
  \end{tikzcd}
  \end{center}
\end{itemize}
\end{definition}

Note that the bifunctor is suggestively denoted by $\otimes$ which hints at one of its popular names, the \emph{"tensor product"}\index{Monoidal category!tensor product}. While this name stems from some famous applications of the concept of a monoidal category, other authors (see \cite{Categories}) denote the defining bifunctor using $"\square"$ or $"\times"$ while refraining from calling it the tensor product to avoid any terms that hint towards special applications. While the name "tensor product" won't be used unless fitting the context, the bifunctor will continue to be denoted by $"\otimes"$.\\

To illuminate the widespread occurence of monoidal categories, partly in very familiar settings, we provide some examples that will go without detailed proof:

\begin{example}~
\begin{itemize}
  \item Taking the category of abelian groups $\cat{Ab}$\index{The category of!abelian groups}, one can use the usual tensor product of of its objects $A,B$ defined by the mapping
  $$ A \times B \lra A \otimes B, \quad \quad (a,b) \lmap a \otimes b $$
  together with the unique natural isomorphism $\alpha_{ABC}: A \otimes (B \otimes C) \lra A \otimes (B \otimes C)$ and the two natural isomorphisms $\lambda_A: A \otimes \ZZ \cong A$ and $\rho_A: \ZZ \otimes A \cong A$ to form the monoidal category of abelian groups $\pair{\cat{Ab}}{\otimes, \ZZ, \alpha, \lambda, \rho}$\index{The monoidal category of!abelian groups}. Note that $\ZZ$ acts as the unit of the tensor product.

  \item Since we already discussed the action of $\times$ on the category of small sets $\cat{Set}$ in the context of product categories, it is only natural to expand the well-known isomorphisms that come with $\times$. First we denote that $\times$ is a bifunctor on $\cat{Set}$. As a unit for this bifunctor, we can use a singleton set $\{*\}$. Together with the usual identifications $\alpha: U \times (V \times W) \cong (U \times V) \times W$ and $\rho: \{*\} \times V \cong V \cong V \times \{*\}: \lambda$ for any $U,V,W \in \cat{Set}$, we obtain the monoidal category $\pair{\cat{Set}}{\times, \{*\}, \alpha, \lambda, \rho}$\index{The monoidal category of!small sets}.

  \item Turning to one of the most important categories for our purposes, the category of vector spaces over a field $\mathds{k}$ denoted by $\vectk$, we can once again use the usual tensor product of vector spaces. For this tensor product, the field $\mathds{k}$ acts as a unit which allows for the following natural isomorphisms derived from the known isomorphisms for the tensor product for $\mathds{k}$-vector spaces $U,V,W$:
  $$ \alpha_{UVW}: U \otimes (V \otimes W) \lra U \otimes (V \otimes W), \quad \quad \rho_V: \mathds{k} \otimes V \cong V \cong V \otimes \mathds{k}: \lambda_V $$
  We thus obtain the monoidal category $\pair{\vectk}{\otimes, \mathds{k}, \alpha, \lambda, \rho}$\index{The monoidal category of!vector spaces}.
\end{itemize}
\end{example}

While monoidal categories provide a very powerful tool to identify categories equipped with an "almost associative" bifunctor, we can restrict to categories whose tensor product is "completely associative":

\begin{definition}[Strict monoidal categories]
\label{def:strict_monoidal}
  A monoidal category $\pair{\cat{C}}{\otimes, e, \alpha, \lambda, \rho}$ is a \textbf{strict monoidal category}\index{Monoidal category!strict} if the three natural isomorphisms $\alpha,\lambda, \rho$ are naturally isomorphic to the respective identity isomorphisms. We usually omit them from our notation and denote a strict monoidal category by $\pair{\cat{C}}{\otimes, e}$.
\end{definition}

We will soon see that we can in fact treat monoidal categories as their strict counterparts. To unambigously formulate such an equivalence we need to formulate the notion of functors and natural transformations in the monoidal context. Thus our goal is to restrict \ref{def:functors} and \ref{def:natural_transformation} to respect the monoidal structure.

\begin{definition}[(Strict) monoidal functors]
\label{def:monoidal_functors}
  Given (strict) monoidal categories $\pair{\cat{C}}{\otimes, e, \alpha, \lambda, \rho}$ and $\pair{\cat{D}}{\widetilde{\otimes}, \widetilde{e}, \widetilde{\alpha}, \widetilde{\lambda}, \widetilde{\rho}}$ we call a functor $F: \cat{C} \lra \cat{D}$ a \textbf{(strict) monoidal functor}\index{Functor!strict monoidal}\index{Functor!monoidal} and denote it by
  $$\cat{F} : \pair{\cat{C}}{\otimes, e, \alpha, \lambda, \rho} \lra \pair{\cat{D}}{\widetilde{\otimes}, \widetilde{e}, \widetilde{\alpha}, \widetilde{\lambda}, \widetilde{\rho}}$$
  if it satisfies the following properties for any objects $a,b,c \in\cat{C}$ and arrows $f,g$ in $\cat{C}$:
  \begin{align*}
    \cat{F}(a \otimes b) \ &\cong\  \cat{F}a \widetilde{\otimes} \cat{F}b, \quad &\quad \cat{F}(f \otimes g) \ &\cong\  \cat{F}f \widetilde{\otimes} \cat{F}g, \quad &\quad \cat{F}e \ &=\  \widetilde{e}\\
    \cat{F}\alpha_{abc} \ &=\  \widetilde{\alpha}_{\cat{F}a \cat{F}b \cat{F}c }, \quad &\quad \cat{F}\lambda_a \ &=\  \widetilde{\lambda}_{\cat{F}a}, \quad &\quad  \cat{F}\rho_a \ &=\  \widetilde{\rho}_{\cat{F}a}
  \end{align*}
  We call a monoidal functor an \textbf{equivalence of monoidal categories}\index{Monoidal category!equivalence} if it is an equivalence of catgeories in the usual sense. For such an equivalence of monoidal categories we call the respective categories \textbf{monoidally equivalent}. Also note that the composition of monoidal functors is again monoidal.
\end{definition}

Unpacking the above definition reveals that a monoidal functor does indeed respect the monoidal structure in that it commutes with the bifunctor $\otimes$ and maps its unit and the respective natural isomorphisms onto one another. Having morphisms between monoidal catgeories at hand, we can build the category of small monoidal categories, $\cat{Moncat}$\index{The category of!small monoidal categories}. One can indeed show that this category hosts the full subcategory \ref{def:full_subcategory} of strict monoidal categories denoted by $\cat{Moncat}^\cat{S}$\index{The category of!small strict monoidal categories}. Further note that for strict monoidal categories, we can completely omit the second row of conditions in the definition of a monoidal functor and replace all occurences of "$\cong$" by "=".\\

We complete the transition to the monoidal setting by defining \emph{natural transformations of monoidal functors}. The definition adapts the one given in \cite[chapter 2.4.]{Tensor_Categories} to the chosen approach to monoidal categories which stems from \cite{Categories}:

\begin{definition}[Monoidal natural transformation]
\label{def:monoidal_transformation}
  Let $\pair{\cat{C}}{\otimes, e, \alpha, \lambda, \rho}$ and $\pair{\cat{D}}{\widetilde{\otimes}, \widetilde{e}, \widetilde{\alpha}, \widetilde{\lambda}, \widetilde{\rho}}$ be two monoidal categories and $\cat{F}, \cat{G}: \cat{C} \lra \cat{D}$
  two monoidal functors between them. We call a natural transformation $\eta: \cat{F} \nat \cat{G}$ a \textbf{morphism of monoidal functors} or a \textbf{monoidal natural transformation}\index{Natural transformation!monoidal}, if $\eta_e = \widetilde{e}$ and the following diagram commutes for all $a,b \in \cat{C}$:
  \begin{center}
  \begin{tikzcd}[sep = large]
    \cat{F}(a \otimes b) \arrow[r, leftrightarrow, "\cong"] \arrow[d, "\eta_{a \otimes b}"'] & \cat{F}a \ \widetilde{\otimes}\  \cat{F}b \arrow[d, "\eta_a \widetilde{\otimes} \eta_b"]\\
    \cat{G}(a \otimes b) \arrow[r, leftrightarrow, "\cong"] & \cat{G}a \ \widetilde{\otimes}\  \cat{G}b
  \end{tikzcd}
  \end{center}
  For a natural transformation of strict monoidal functors, we can reduce these conditions to $\eta_e = \widetilde{e}$ and $\eta_{a \otimes b} = \eta_a \widetilde{\otimes} \eta_b$.
\end{definition}

\begin{remark}
  To no surprise, for an equivalence of monoidal categories $\cat{F}: \cat{C} \lra \cat{D}$ there exists an inverse equivalence of monoidal categories $\cat{F}^{-1}: \cat{D} \lra \cat{C}$ such that $\cat{F} \circ \cat{F}^{-1}$ and $\cat{F}^{-1} \circ \cat{F}$ are naturally monoidally isomorphic to the monoidal identity functors $\cat{I}_\cat{C}$ and $\cat{I}_\cat{D}$ respectively.
\end{remark}

Using the notion of natural monoidal transformations we can now formulate a very fundamental theorem when working with monoidal categories.

\begin{theorem}[Strictification Theorem]\index{Monoidal category!Strictification Theorem}
\label{theo:strictification_theorem}
  Every monoidal category is monoidally equivalent to a strict monoidal category.
\begin{proof}
  This proof adapts the one of \emph{Theorem 2.8.5.} in \cite{Tensor_Categories} to the given definition of a monoidal category \ref{def:monoidal_category}. Thus let $\pair{\cat{C}}{\otimes, e, \alpha, \lambda, \rho}$ be a monoidal category. Our goal is to construct a strict monoidal category using $\cat{C}$ and then showing that there exists a monoidal equivalence of categories between them.\\

  Thus let $\widetilde{\cat{C}}$ be the category defined as follows: Its objects are pairs $\pair{\cat{F}}{\eta}$ where $\cat{F}\colon \cat{C} \lra \cat{C}$ is a functor and for any $a,b\in \cat{C}$
  $$ \eta_{ab} \colon \cat{F}(a \tensor b) \cong \cat{F}(a) \tensor b$$
  is a natural isomorphism such that the following diagram commutes for all $a,b,c \in \cat{C}$:
  \begin{center}
  \begin{tikzcd}[sep = large]
    & (\cat{F}(a) \tensor b ) \tensor c & \\
    \cat{F}(a \tensor b) \tensor c \arrow[ur, "\eta_{ab} \tensor \id_c"] &   & \cat{F}(a) \tensor (b \tensor c) \arrow[ul, "\alpha_{\cat{F}(a), b, c}"'] \\
    \cat{F}((a\tensor b) \tensor c) \arrow[u, "\eta_{a \tensor b c}"] &  & \cat{F}(a\tensor (b \tensor c))\arrow[ll,"\cat{F}\left(\alpha_{abc}\right)"'] \arrow[u, "\eta_{a b \tensor c}"']
  \end{tikzcd}
  \end{center}
  We define the arrows of $\widetilde{\cat{C}}$, denoted as $\tau\colon (\cat{F}_1, \eta^1) \lra (\cat{F}_2, \eta^2)$, to be natural transformations $\tau: \cat{F}_1 \nat \cat{F}_2$ which make the following diagram commute for any $a,b \in \cat{C}$:
  \begin{center}
  \begin{tikzcd}[sep = large]
     \cat{F}_1(a \tensor b) \arrow[r,"\eta^1_{ab}"] \arrow[d,"\tau_{a\tensor b}"'] & \cat{F}_1(a) \tensor b \arrow[d,"\tau_a \tensor \id_b"] \\
     \cat{F}_2(a \tensor b) \arrow[r,"\eta^1_{ab}"] & \cat{F}_2(a) \tensor b
  \end{tikzcd}
  \end{center}
  Now all that is left to make $\widetilde{\cat{C}}$ a category is the composition of arrows. But since we defined its arrows to be natural transformations, we can simply make use of their usual composition. The next step is to show that $\widetilde{\cat{C}}$ is indeed a monoidal category:\\

  Let $(\cat{F}_1, \eta^1) $ and $(\cat{F}_2, \eta^2)$ be two objects in $\widetilde{\cat{C}}$. We define a bifunctor $\tensor\colon \widetilde{\cat{C}} \times \widetilde{\cat{C}} \lra \widetilde{\cat{C}}$ with an object function $(\cat{F}_1, \eta^1) \tensor (\cat{F}_2, \eta^2) = (\cat{F}_1 \cat{F}_2, \eta)$ where $\eta$ is defined by the following composition for any $a,b\in \cat{C}$:
  \begin{center}
  \begin{tikzcd}[sep = large]
     \cat{F}_1 \cat{F}_2(a \tensor b) \arrow[r,"\cat{F}_1(\eta^2_{ab})"]& \cat{F}_1 (\cat{F}_2(a) \tensor b) \arrow[r,"\eta^1_{\cat{F}_2(a)b}"] &\cat{F}_1 \cat{F}_2(a) \tensor b
  \end{tikzcd}
  \end{center}
  The tensor product for arrows is, again, the usual composition of natural transformations. This immediately shows that $\widetilde{\cat{C}}$ is a \emph{strict monoidal category} with unit $(\cat{I}_\cat{C}, \id)$.\\

  Thus all that's left is to construct a monoidal equivalence between $\cat{C}$ and $\widetilde{\cat{C}}$. To this end, we define a functor
  \begin{align*}
    \cat{L} \colon \cat{C} & \lra \widetilde{\cat{C}}, & a &\lmap (a \tensor -, \alpha_{a--}), & f &\lmap (f \tensor -)
  \end{align*}
  where $(-)$ denotes an argument. Note that writing the above pentagon diagram for our newly defined functor $\cat{L}$ amounts to writing the pentagon equation \ref{def:monoidal_category} for $\cat{C}$. Further note that since for any $a, b \in \cat{C}$
  $$ e \tensor a \cong a, \quad \quad \text{and} \quad \quad \alpha_{eab}: e \tensor (a \tensor b)  \cong (e \tensor a) \tensor b$$
  any $(\cat{F}, \eta)$ in $\widetilde{\cat{C}}$ is isomorphic to $\cat{L}(\cat{F}(e))$. This property marks the first requirement to satify \ref{theo:equivalence}. Now let $\tau: \cat{L}(a) \lra \cat{L}(b)$ be a morphism in $\widetilde{\cat{C}}$ and define an arrow of $\cat{C}$, $f:a \lra b$ by the composition
  $$ a \atob{\lambda^{-1}_a} a \tensor e \atob{\tau_e} b \tensor e \atob{\lambda_b} b $$
  To show that $\cat{L}$ is full, we need to prove that for all $c \in \cat{C}$ one has $\tau_c = f \tensor \id_c$. Since $\cat{L}(f) = (f \tensor -)$ this would prove $\tau = \cat{L}(f)$ and thus $\tau$ would be full by definition \ref{def:full_faithful}. Indeed, the following diagram commutes for any $c \in \cat{C}$:
  \begin{center}
  \begin{tikzcd}[sep = huge]
    a \tensor c \arrow[d,"\tau_c"'] \arrow[r,"\id_a \tensor \rho^{-1}_c"] & a \tensor (e \tensor c) \arrow[d,"\tau_{e \tensor c}"'] \arrow[r,"\alpha_{aec}"] & (a \tensor e ) \tensor c \arrow[d,"\tau_e \tensor \id_c"] \arrow[r,"\lambda_a \tensor \id_c"] & a \tensor c \arrow[d, "f \tensor \id_c"]\\
    b \tensor c \arrow[r,"\id_b \tensor \rho^{-1}_c"'] & b \tensor (e \tensor c) \arrow[r,"\alpha_{bec}"'] & (b \tensor e) \tensor c \arrow[r,"\lambda_b \tensor \id_c"'] & b \tensor c
  \end{tikzcd}
  \end{center}
  This can be seen as follows: The rows commute by rule of the triangle identity \ref{def:monoidal_category}. The left quare commutes since $\tau$ is a natural transformation \ref{def:natural_transformation}, the right square by definition of $f$. For the remaining square note that $\tau$ is an arrow in $\widetilde{\cat{C}}$ and thus a natural transformation.\\
  To show that $\cat{L}$ is faithful we take two arrows $f,g$ in $\cat{C}$. Now if $\cat{L}(f) = \cat{L}(g)$ we particularly have $f \tensor \id_e = g \tensor \id_e$ and by rule of the natural isomorphism $\lambda$, $f = g$. Thus $\cat{L}$ is faithful. Putting all this together, $\cat{L}$ satisfies the requirements of \ref{theo:equivalence} thus forming an equivalence of categories.\\

  The last part of the proof is to show that $\cat{L}$ is a monoidal functor between categories \ref{def:monoidal_functors}. First we note that $\cat{L}(e)$ is isomorphic to the identity in $\widetilde{\cat{C}}$. Further, for any two arrows $f,g \in \cat{C}$ we have
  $$\cat{L}(f \tens{\cat{C}} g)  = ((f \tens{\cat{C}} g) \tens{\cat{C}} -) \cong (f \tens{\cat{C}} (g \tens{\cat{C}} -)) = \cat{L}f \tens{\widetilde{\cat{C}}} \cat{L}g$$
  For $\cat{L(\alpha_{abc})} \cong \widetilde{\alpha}_{\cat{L}(a)\cat{L}(b)\cat{L}(c)}$ we note that the respective commuting hexagon diagram reduces to the \emph{pentagon identity} in $\cat{C}$. To show $\cat{L}(\lambda_a) \cong \widetilde{\lambda}_{\cat{L}(a)}$ (and thus analogously for $\rho_a$) we can investigate the commuting diagram
  \begin{center}
  \begin{tikzcd}[sep = large]
    & (a \tensor e) \tensor b \arrow[dl, "\lambda_a \tensor \id_b"'] & \\
    a \tensor b \arrow[rr, "\id_a \tensor \lambda^{-1}_b"]& & a \tensor (e \tensor b) \arrow[ul, "\alpha_{aeb}"']
  \end{tikzcd}
  \end{center}
  to see that both hold. Now all that's left is $\cat{L}(a \tensor b) \cong \cat{L}(a) \tens{\widetilde{\cat{C}}} \cat{L}(b) $ which, like for arrows, reduces to composition of functors and thus holds. This concludes the proof.

\end{proof}
\end{theorem}

By rule of \ref{theo:strictification_theorem} we omit the distinction between strict and "relaxed" monoidal categories and just use the term "monoidal" category for the respective strict monoidal category. Thus we denote a monoidal category by $\pair{\cat{C}}{\otimes, e}$ while still using the natural isomorphisms $\alpha, \lambda, \rho$ occasionally. Using this notation, we find that monoidal categories harbour monoids in a natural manner:

\begin{definition}[Monoids of a monoidal categories]
\label{def:monoidal_monoids}
  Let $\triple{\cat{C}}{\otimes}{e}$ be a monoidal category. A \textbf{monoid}\index{Monoidal category!monoid} in $\cat{C}$ is a triple $\triple{c}{\mu}{\eta}$ where $c\in \cat{C}$ is an object, $\mu: c \otimes c \lra c$ a bifunctor often called "multiplication" and $\eta: e \lra c$ an arrow often called a "unit" such that the following two diagrams commute:
  \begin{center}
  \begin{tikzcd}
    c \otimes (c \otimes c) \arrow[r, "\alpha_{ccc}"] \arrow[d, "\id_c \otimes \mu"'] & (c \otimes c) \otimes c \arrow[r, "\mu \otimes \id_c"] & c\otimes c \arrow[d, "\mu"] & e \otimes c \arrow[dr, "\rho"'] \arrow[r, "\eta \otimes \id_c"]& c \otimes c \arrow[d, "\mu"] & c \otimes e \arrow[l, "\id_c \otimes \eta"'] \arrow[dl, "\lambda"]\\
    c \otimes c \arrow[rr, "\mu"] & & c & & c &
  \end{tikzcd}
  \end{center}
\end{definition}

A closely related concept to that of a monoid in a monoidal category is that of a comonoid which we will introduce here to enbable a more categorial point of view in chapter \ref{sec:Comfrob}.

\begin{definition}[Comonoids of a monoidal categories]
\label{def:monoidal_comonoids}
  Let $\triple{\cat{C}}{\otimes}{e}$ be a monoidal category. A \textbf{comonoid}\index{Monoidal category!comonoid} in $\cat{C}$ is a triple $\triple{c}{\delta}{\epsilon}$ where $c\in \cat{C}$ is an object, $\delta: c \lra c \otimes c $ a morphism often called "comultiplication" and $\epsilon: c \lra e$ an arrow often called a "counit" such that the following two diagrams commute:
  \begin{center}
  \begin{tikzcd}
    c \otimes (c \otimes c) \arrow[r, leftarrow, "\simeq"] \arrow[d, leftarrow, "\id_c \otimes \delta"'] & (c \otimes c) \otimes c \arrow[r, leftarrow, "\delta \otimes \id_c"] & c\otimes c  & e \otimes c \arrow[dr, leftrightarrow, "\simeq"'] \arrow[r, leftarrow, "\epsilon \otimes \id_c"]& c \otimes c & c \otimes e \arrow[l, leftarrow, "\id_c \otimes \epsilon"'] \arrow[dl, leftrightarrow, "\simeq"]\\
    c \otimes c & & c \arrow[ll, "\delta"'] \arrow[u, "\delta"'] & & c \arrow[u, "\delta"'] &
  \end{tikzcd}
  \end{center}
\end{definition}

\begin{remark}
  Note that every category $\cat{C}$ has a so-called "opposite category"\index{Category!opposite} $\cat{C}^{op}$ whose objects are the same as those of $\cat{C}$ but a morphism $f \colon a \lra b$ in $\cat{C}$ is the same as an arrow $f \colon b \lra a$ in $\cat{C}^{op}$. Further a composite of arrows $f \circ g$ in $\cat{C}$ is defined to be the composite $g \circ f$ in $\cat{C}^{op}$.
  While we won't further use opposite category, note that every oppposite of a monoidal category naturally inherits a monoidal structure itself. Further comonoids are nothing but the monoids of the opposite category.
\end{remark}

Needless to say there are special arrows that respect the structure of monoids:

\begin{definition}[Morphism of Monoids]
\label{def:monoidal_monoid_arrows}
  Given two monoids $\triple{c}{\mu}{\eta}$ and $\triple{d}{\widetilde{\mu}}{\widetilde{\eta}}$ and an arrow $f: c \lra d$ in the respective monoidal category, we call $f$ a \textbf{morphism of monoids}\index{Monoidal category!monoid!morphism} and denote it by $f: \triple{c}{\mu}{\eta} \lra \triple{d}{\widetilde{\mu}}{\widetilde{\eta}}$ if $f$ satisfies
  $$ f \mu = \widetilde{\mu}(f \otimes f) : c \otimes c \lra d, \quad \quad \quad f \eta = \widetilde{\eta} : e \lra d $$
  In the same manner we can define a morphism of comonoids\index{Monoidal category!comonoid!morphism} $\triple{c}{\delta}{\epsilon}$ and $\triple{d}{\widetilde{\delta}}{\widetilde{\epsilon}}$ as an arrow $f \colon c \lra d$ such that
  $$ (f \tensor f) (\delta) = \widetilde{\delta}(f) : c \lra d \tensor d, \quad \quad \quad \eta = \widetilde{\eta}(f) : c \lra e $$
\end{definition}

Using the morphisms of monoids as arrows, one can define the category of monoids of a monoidal category $\triple{\cat{C}}{\otimes}{e}$ denoted by $\cat{Mon}_\cat{C}$\index{The category of!monoids of a monoidal category} (the same works for comonoids with $\cat{coMon}_\cat{C}$\index{The category of!comonoids of a monoidal category}).
One can define the evident forgetful functor $F:\cat{Mon}_\cat{C} \lra \cat{C}$ that maps $\triple{c}{\mu}{\eta} \lmap c$ and thus "forgets" the monoid structure defined on "c". To shed some light on the notion of monoids in a monoidal category, we provide some examples:

\begin{example}~
\begin{itemize}
  \item Monoids in the monoidal category $\triple{\cat{Ab}}{\otimes}{\ZZ}$ are precisely rings\index{Monoidal category!monoid!in \cat{Ab}}.

  \item A monoid in the monoidal category of small sets $\triple{\cat{Set}}{\times}{\{*\}}$ is a set equipped with an associative multiplication and a respective identity element. Thus a monoid in $\cat{Set}$ agrees with the previously defined notion of a monoid \ref{def:monoids}\index{Monoidal category!monoid!in \cat{Set}}.

  \item Again looking at the monoidal category $\triple{\vectk}{\otimes}{\mathds{k}}$, the monoids are vector spaces equipped with a ring structure that respects the $\mathds{k}$-vector space structure, namely $\kk$-algebras\index{Monoidal category!monoid!in $\vectk$} which we will encounter in detail in a later chapter.
\end{itemize}
\end{example}

Generalizing the ideas of left- and right-modules as actions of a ring on an abelian group, we define special types of bifunctors that utilize monoids:

\begin{definition}[Action of a Monoid]
\label{def:monoid_action}
  Given a monoidal category $\triple{\cat{C}}{\otimes}{e}$ and a monoid $\triple{c}{\mu}{\eta}$ therein, a \textbf{left action of a monoid}\index{Monoidal category!action of a monoid} on an object $a \in \cat{C}$ is an arrow $\sigma: c \otimes a \lra a$ such that the following diagram commutes:
  \begin{center}
  \begin{tikzcd}
  c \otimes (c \otimes a) \arrow[r, "\alpha_{cca}"] \arrow[d, "\id_c \otimes \sigma"'] & (c \otimes c) \otimes a \arrow[r, "\mu \otimes \id_a"] & c\otimes a \arrow[d, "\sigma"] & e \otimes a \arrow[l, "\eta \otimes \id_a"'] \arrow[dl, "\lambda"]\\
  c \otimes a \arrow[rr, "\sigma"] & & a &
  \end{tikzcd}
  \end{center}
  A \textbf{right action of a monoid} is defined in a similar fashion.
\end{definition}

As mentioned prior to the definition of an action of a monoid, there are some interesting familiars emerging from our previously given examples of monoids in monoidal categories:

\begin{example}~
\begin{itemize}
  \item As previously mentioned, the monoidal category $\triple{\cat{Ab}}{\otimes}{\ZZ}$ has rings as its monoids. Taking a Ring $\mathcal{R}$ and an abelian group $\mathcal{A} \in \cat{Ab}$, the usual left (or right) action of $\mathcal{R}$ on $\mathcal{A}$ gives it a left (or right) $\mathcal{R}$-module structure.

  \item For a monoid $\triple{M}{\mu}{\eta}$ in $\triple{\cat{Set}}{\times}{\{*\}}$, the left (or right) action of $M$ on a set $U \in \cat{Set}$ is given by the usual left (or right) $M$-action of a monoid on a set.

  \item For a monoid $\triple{M}{\mu}{\eta}$ and an object $V$ in $\triple{\vectk}{\otimes}{\mathds{k}}$, the left (or right) action of $M$ on $V$ defined the structure of a left (or right) $\mathds{k}$-algebra module.
\end{itemize}
\end{example}

As for every categorial object we investigated so far, we can define specialized morphisms between them:

\begin{definition}[Morphisms of Actions of Monoids]
\label{def:monoid_action_morphisms}
  Given a monoidal category $\triple{\cat{C}}{\otimes}{e}$, three elements $a,b,c \in \cat{C}$, a monoid $\triple{a}{\mu}{\eta}$, a left action of $a$ on $b$ and one on $c$, denoted by $\sigma$ and $\lambda$ respectively, a \textbf{morphism of left actions}\index{Monoidal category!action of a monoid!morphism} $f \colon \sigma \lra \lambda$ is an arrow $f \colon b \lra c$ in $\cat{C}$ such that the following diagram commutes:
  \begin{center}
  \begin{tikzcd}[sep=huge]
    a \tensor b \arrow[d, "\sigma"'] \arrow[r, "\id_a \tensor f"]  & a \tensor c \arrow[d, "\lambda"] \\
    b \arrow[r, "f"'] & c
  \end{tikzcd}
  \end{center}
  Again, a \textbf{morphism of right actions} is defined analogously.
\end{definition}

Note that these morphisms enable a definition of the categories of right and left actions for a given monoid $a$. We will denote them by $\ract{a}$\index{The category of!right actions} and $\lact{a}$\index{The category of!left actions} respectively.\\

To conclude this subchapter, note that up until here, we explored the idea of a category equipped with a bifunctor of varying associativity. Starting with the general notion of product categories, we soon saw that monoidal categories provide an excellent tool to describe "associativity up to isomorphism". After they turned out to be equivalent to their strict counterparts which embodied full associativity, we saw that the concept of monoids from \ssecref{subsec:introduction_categories} naturally occurs within strict monoidal categories. Arguably the most important takeaway from this chapter is that $\vectk$ does indeed qualify as a monoidal category which will enable the entire categorial formulation of chapter \ref{sec:Comfrob}.

\newpage
\subsection{Braided and Symmetric Monoidal Categories}
\label{subsec:braided_categories}

The previous chapter explored \emph{monoidal categories} which made precise the idea of a category equipped with an internal product of varying associativity between its objects thus providing a categorial version of a commutative \emph{monoid}. The next natural step is to formulate the \emph{distributivity} of such a product which will lead us to the notion of \emph{braided monoidal categories}. As a special case of these categories, we will then encounter \emph{symmetric monoidal categories} which mark the final construct needed to understand the remarkable internal structure making $\vectk$ a fitting candidate in constructing the notion of a \emph{TQFT} from a mathematical and physical perspective.\\

This chapter again makes use of definitions appearing in \cite{Categories} and \cite{Tensor_Categories} while adapting to the previously chosen notation. Some of the proactive comments on \emph{TQFT's} will refer to discussions found in \emph{chapter 2} of \cite{Intro_TQFT} and \cite{Atiyah_1}.

\begin{definition}[Braiding]
\label{def:braiding}
  Let $\pair{\cat{C}}{\otimes, e, \alpha, \lambda, \rho}$ be a monoidal category. A \textbf{braiding}\index{Braiding} of $\cat{C}$ is a collection of natural isomorphisms, natural in all $a,b\in \cat{C}$
  $$ \gamma_{ab}\colon a \otimes b \cong b \otimes a $$
  that satisfy the following "coherence" properties\index{Braiding!coherence} with respect to the isomorphisms $\alpha, \lambda, \rho$:
  \begin{center}
  \begin{tikzcd}[sep = large]
    (ab)c \arrow[r,"\gamma_{(ab)c}"] \arrow[d,"\alpha^{-1}_{abc}"'] & c(ab) \arrow[d, "\alpha_{cab}"] & a(bc) \arrow[d,"\alpha_{abc}"'] \arrow[r,"\gamma_{a(bc)}"] & (bc)a \arrow[d,"\alpha^{-1}_{bca}"]  & a  \arrow[r, leftrightarrow, "="]& a\\
    a(bc) \arrow[d,"\id_a \tensor \gamma_{bc}"'] & (ca)b \arrow[d,"\gamma_{ca} \tensor \id_b"] & (ab)c \arrow[d,"\gamma_{ab} \tensor \id_c"'] & b(ca) \arrow[d,"\id_b \tensor \gamma_{ca}"] & &\\
    a(cb) \arrow[r,"\alpha_{acb}"] & (ac)b & (ba)c \arrow[r, "\alpha^{-1}_{bac}"] & b(ac) & a \tensor e \arrow[uu, "\lambda"] \arrow[r, "\gamma_{ae}"]& e \tensor a \arrow[uu, "\rho"']
  \end{tikzcd}
  \end{center}
  In the above diagram, we suppressed many instances of "$\tensor$" to shorten notation. Dropping the indices, we usually denote \emph{a braiding} (namely the collection of isomorphisms) by its respective greek letter.
\end{definition}

\begin{remark}
\label{rem:BraidingNatural}
  Note that the naturality of the braiding\index{Braiding!naturality} in its two factors is a quite strong statement. Namely for any two arrows $f \colon a \lra b$ and $g \colon c \lra d$ in $\cat{C}$ we have $(g \tensor f) \circ \gamma_{ac} \simeq \gamma_{bd} \circ (f \tensor g)$. This property will be extremely helpful when talking about cobordisms \ref{sec:Cobordisms} and Frobenius Algebra \ref{sec:Comfrob} since it allows us to permute certain maps resulting in quite general insights about the respective structures.
\end{remark}

To no surprise, we can specialize a monoidal category when equipping it with a braiding to implement distributivity up to isomorphism:

\begin{definition}[Braided monoidal category]
\label{def:monoidal_braided}
  Given a monoidal category $\triple{\cat{C}}{\tensor}{e}$ and a braiding $\gamma$ of \cat{C}, a \textbf{braided monoidal category}\index{Braiding!braided monoidal category}\index{Monoidal category!braided} is the tuple $\triple{\cat{C}}{\tensor}{e, \gamma}$ which we abbreviate by $\pair{\cat{C}}{\gamma}$.
\end{definition}

It is important to note that one and the same monoidal category can be used to form different braided monoidal categories depending on the respective braiding. To compare and identify the respective structures defined by a braiding, we will again need adapted notions of functors and natural transformations. Before turning to those definitions, we provide two interesting examples of braided categories:

\begin{example}~
\begin{itemize}
  \item The monoidal category of small sets, $\triple{\cat{Set}}{\times}{\{*\}}$, can be equipped with the family of transposition isomorphism $\gamma_{UV}\colon U \times V \cong V \times U$ for any $U,V \in \cat{Set}$. They act like a braiding forming the braided category $\pair{\cat{Set}}{\gamma}$\index{The braided monoidal category of!small sets}\index{The symmetric monoidal category of!small sets}.

  \item We already presented the category of vector spaces over a field $\mathds{k}$ as a monoidal categors. In fact, it can be endowed with a braiding induced by the natural isomorphism of the tensor product of any two objects $V,W$:
  $$ \gamma_{VW}: V \tensor W \atob{\sim} W \tensor V, \quad \quad v \tensor w \lmap w \tensor v$$
  Thus we obtain the braided monoidal category $\pair{\vectk}{\gamma}$\index{The braided monoidal category of!vector spaces}\index{The symmetric monoidal category of!small sets}.
\end{itemize}
\end{example}

To compare the structure of braided monoidal categories, the notion of a functor needs to be adapted to the respective braidings. Very similar to the extra conditions we introduced when defining \emph{monoidal functors} \ref{def:monoidal_functors} this will lead to some coherence equations. This definition and the following definition of a braided monoidal natural transformation are adapted from \cite{Categories}.

\begin{definition}[Braided monoidal functors]
\label{def:monoidal_braided_functor}
  Given two braided monoidal categories $\pair{\cat{A}}{\gamma}$ and $\pair{\cat{B}}{\widetilde{\gamma}}$ we call a monoidal functor $\cat{F}: \cat{A} \lra \cat{B}$ a \textbf{braided monoidal functor}\index{Braiding!braided functor}\index{Functor!braided monoidal} if the following diagram commutes for any $a,b \in \cat{A}$:
  \begin{center}
  \begin{tikzcd}[sep = large]
    Fa \tensor Fb \arrow[d,"\cong"'] \arrow[r,"\widetilde{\gamma}"] & Fb \tensor Fa \arrow[d,"\cong"] \\
    F(a \tensor b) \arrow[r,"F(\gamma)"] & F(b \tensor a)
  \end{tikzcd}
  \end{center}
  Thus a braided monoidal functor "commutes" with the involved braidings. Again the composition of two braided functors is again braided.
\end{definition}

Note that over a base category $\cat{C}$ one can again define a category using braided monoidal functors as arrows and braided monoidal categories as objects, often denoted by $\cat{BrMonCat}_\cat{C}$\index{The category of!braided monoidal categories}.\\

While by now we now have a pretty solid grasp on the categorial structure of $\vectk$ and, as will turn out in Section \ref{sec:Cobordisms}, on that of the category of cobordisms, $\cobn{n}$, there is yet another refinement to the notion of general braidings. Since this will indeed become relevant on $\vectk$ as well as $\cobn{n}$, it is indispensable for our purposes. Considering the idea of a strict monoidal category \ref{def:strict_monoidal} which embodies the idea of a fully associative tensor product, we aim to define braided monoidal categories with a "fully distributive" braiding.

\begin{definition}[Symmetric monoidal categories]
\label{def:monoidal_symmetric}
  Let $\pair{\cat{C}}{\gamma}$ be a braided monoidal category. We call it a \textbf{symmetric monoidal category}\index{Monoidal category!symmetric}, if for every $a,b \in \cat{C}$ the following diagram commutes:
  \begin{center}
  \begin{tikzcd}[sep = large]
    & b \tensor a \arrow[dr, "\gamma_{ba}"] & \\
    a \tensor b \arrow[ur, "\gamma_{ab}"] \arrow[rr, leftrightarrow, "\id_{a \tensor b}"] & & a \tensor b
  \end{tikzcd}
  \end{center}
  This can be formulated as $\gamma_{ab} \circ \gamma_{ba} = \id_{a\tensor b}$\index{Braiding!symmetric}.
\end{definition}

Note that both examples we gave for braided monoidal categories are indeed symmetric monoidal categories. Also the definition of a braided monoidal functor naturally extends to that of a \textbf{symmetric monoidal functor}\index{Functor!symmetric monoidal}. Their composition is also a symmetric monoidal functor.\\

The above definition sheds a lot of light on the name "braiding". A symmetric braiding can be undone by applying the naive inverse braiding morphism. Meanwhile a general braiding can behave like a "twist" operation that is not involutory. While both types of braidings have applications in Quantum Mechanics (see \emph{page 251} of \cite{Categories}), we will mainly work with categories that qualify as symmetric monoidal categories. At this point, we will take a short detour and investigate the notion of dual objects. This provides an interesting result for $\finvect$.

\begin{definition}[Dualisable object]
  Given $\cat{C} \in \cat{MonCat}$ denoted by $\triple{\cat{C}}{\tensor}{e}$ an object $A$ is \textbf{left dualisable} if there exists an object $A^* \in \cat{C}$ and two morphisms:
  \begin{enumerate}
    \item The \textbf{evaluation map}
    \begin{equation}
      \cat{ev}_A \colon A^* \tensor A \lra e
    \end{equation}
    \item The \textbf{coevaluation map}
    \begin{equation}
      \imath_A \colon e \lra  A \tensor A^*
    \end{equation}
  \end{enumerate}
  such that the following two diagrams commute:
  \begin{center}
  \begin{tikzcd}[sep=huge]
    A^* \tensor (A \tensor A^*) \arrow[d,"\simeq"'] & \arrow[d,"\simeq"] \arrow[l, "\id \tensor \imath_A"'] A^* \tensor e  &  (A \tensor A^*) \tensor A \arrow[d,"\simeq"'] & \arrow[d,"\simeq"] \arrow[l, "\imath_A \tensor \id"'] e \tensor A \\
    (A^* \tensor A) \tensor A^* \arrow[r, "\cat{ev}_A \tensor \id"'] & e \tensor A^* &  A \tensor (A^* \tensor A) \arrow[r, "\id \tensor \cat{ev}_A"'] & A \tensor e
  \end{tikzcd}
  \end{center}
  The notion of \textbf{right dualisability} follows analogously.
\end{definition}

\begin{lemma}
\label{lemma:RigidComonoid}
  Let there be a monoidal category $\triple{\cat{C}}{\tensor}{e}$ and a monoid/comonoid $\triple{A}{\mu}{\eta}$ therein. If for $A$ there also exists a comonoid/monoid structure, the object $A$ has a left and a right dual given by the object itself.
\begin{proof}
  The results is immediate when comparing the diagrams for the multiplication in the definition of a monoid \ref{def:monoidal_monoids} and for the comultiplication in that of a comonoid \ref{def:monoidal_comonoids} with the above for left and right dualisability.
\end{proof}
\end{lemma}

\begin{definition}[Invertible Object]
  If $A$ is a left or right dualisable object for which $\cat{ev}_A$ and $\imath_A$ are isomorphisms it is called \textbf{invertible}.
\end{definition}

\begin{remark}
  Note that for a symmetric monoidal category, like $\vectk$ there is no difference between a left and a right dual. This is immediately clear since both the left and the right dual for a vector space $V$ are indeed $V^* = \Hom(V, \kk)$.
\end{remark}

\begin{definition}[Rigid Monoidal Category]
\label{def:RigidMonoidalCat}
  If every $A \in \cat{C}$ for $\cat{C} \in \cat{MonCat}$ has a left and right dual, we call $\cat{C}$ a \textbf{rigid monoidal category}. If it is symmetric every object is invertible and we call it a \textbf{compact closed category}. Note that one sometimes refers to objects being \textbf{rigid} if they have a left and a right dual.
\end{definition}

The following results about the rigidness of a monoidal category immediately follows from \ref{lemma:RigidComonoid}:

\begin{corollary}
  If every $A \in \cat{C}$ for $\cat{C} \in \cat{MonCat}$ admits the structure of a monoid and a comonoid, $\cat{C}$ is rigid.
\end{corollary}

\begin{example}
  In the monoidal category of finite-dimensional vector spaces, $\finvect$, every object is left and right dualisable via the usual dual of vector spaces. Further the evaluation and coevaluation morphisms are obviously isomorphisms to the unit element $\kk$. Note that we had to make a restriction to the full subcategory of finite-dimensional vector spaces, a restriction which is needed for almost all examples of rigid monoidal categories.
\end{example}

With the introduction of \emph{symmetric monoidal categories} done the chapter on categorial introduction (modulo proof of \ref{theo:strictification_theorem}) comes to an end.
While far from complete or extensive even, it presents the most important terms for the upcoming introductions and investigations of cobordisms, $\cobn{n}$, commutative Frobenius algebras, $\comk$, and ultimately of the very notion of \emph{Topological Quantum Field Theory}.

\newpage
\section{The Category of Cobordisms}
\label{sec:Cobordisms}

This chapter will introduce and discuss the category of \emph{cobordisms}, $\cobn{n}$, whose objects are smooth oriented closed $(n-1)$-dimensional real manifolds. It is named after its arrows, equivalence classes of $n$-dimensional manifolds that have the disjoint union of two objects as their boundary. As such, $\cobn{n}$ can be seen as an abstract toy model for spacetime; Its objects represent slices through the equivalence classes of spacetime while the arrows function as the (time-)evolution of these slices. While this hints at a rich geometric structure, note that we are interested in Topological Quantum Field Theory, which will in particular not make use of any local geometric properties.\\

The following presentation and discussion of $\cobn{n}$ will mostly follow \textit{chapter 1} of \textit{Joachim Kock's} book on \textit{Frobenius Algebras and 2D Topological Quantum Field Theories} \cite{FrobAlgebraTQFT} while taking into consideration the interesting comments of \cite{Intro_TQFT}.

\subsection{Cobordisms}
\label{subsec:cobordism_def}

While we will ultimately work with oriented cobordisms, this subchapter starts by introducing unoriented ones. They catch the main ideas and concepts while directly motivating the notion of oriented cobordisms.

\begin{definition}[Unoriented cobordism]
\label{def:cobordism_unoriented}
  Let $\Sigma_1$ and $\Sigma_2$ be two smooth compact $(n-1)$-dimensional manifolds without boundary for some $n \in \NN$. A \textbf{cobordism}\index{Cobordism} between $\Sigma_1$ and $\Sigma_2$ is a smooth compact $n$-dimensional manifold $M$ with boundary the disjoint union $\Sigma_1 \cop \Sigma_2$. The manifolds $\Sigma_1$ and $\Sigma_2$ are then called \textbf{cobordant}.
\end{definition}

As noted in \cite{FrobAlgebraTQFT}, the prefix "co-" does not hint towards the concept of duality but rather stems from its meaning as "together". This is due to the fact that a \textit{bordism} was first defined to be a manifold that is the boundary of another manifold, a \textit{cobordism} expanded this concept to two such manifolds. One can represent cobordisms in lower dimensions using catchy diagrams:

\begin{example}
  Some examples in dimension $n = 2$ which should be read from left to right include
  \begin{center}
  \begin{tikzpicture}[tqft/.cd,cobordism/.style={draw},every lower boundary component/.style={draw},every incoming lower boundary component/.style={draw, dashed},every outgoing lower boundary component/.style={draw, solid}]

    \pic[tqft, incoming boundary components=1, outgoing boundary components=1, offset=0, genus=0, rotate=90, name=A, at={(0,0)}];

    \pic[tqft/pair of pants, rotate=90, name=B, at={(4,0)}];
    \pic[tqft/reverse pair of pants, rotate=90, at=(B-outgoing boundary 1)];

    \pic[tqft/pair of pants, rotate=90,name=D, at={(10,0)}];
  \end{tikzpicture}
  \end{center}
  Note that the rightmost example displays a cobordism between a simply connected and a not simply connected manifold. Other interesting visualizations arise when considering one or both manifolds being the empty manifold $\emptyset$:
  \begin{center}
  \begin{tikzpicture}[tqft/.cd,cobordism/.style={draw},every lower boundary component/.style={draw},every incoming lower boundary component/.style={draw, dashed},every outgoing lower boundary component/.style={draw, solid}]

    \pic[tqft, outgoing boundary components = 2, incoming boundary components = 0, rotate=-90, at={(0,1)}];

    \pic[tqft/cap, rotate = 90, at={(0,0)}, name = A];
    \pic[tqft/cup, rotate=90, at=(A-outgoing boundary 1)];

    \pic[tqft/pair of pants, rotate=90, name=C, at={(5,0)}];
    \pic[tqft/cup, rotate=-90, at=(C-incoming boundary 1)];
  \end{tikzpicture}
  \end{center}
\end{example}

This style of visualization using suggestive diagrams, or rather drawings, will permeate this project. It allows for a graphic and more or less unambigous representation of many of the more complex compositions and equivalences when dealing with $\cobtwo$ in later subchapters \ref{subsec:cobordism_2D} and when working with TQFT \ref{sec:TQFT} in general. To no surprise, the examples will prevalently live in $n=2$ or $n=1$. While the former is due to the fact that $1$-dimensional closed manifolds are circles or disjoint unions of circles, the latter actually only hosts points and their disjoint unions as manifolds. This consideration brings us to two interesting statements about these special cases:

\begin{lemma}
  Two closed $0$-manifolds are cobordant iff they have the same number of points $\modulo 2$. Furthermore, any two closed $1$-manifolds are cobordant.
\begin{proof}
  We begin with $n=1$ thus $0$-manifolds:

  Taking two such manifolds with the same number of points $\modulo 2$, we have two cases. If both have an even number of points, we can join them in pairs on both sides by smooth curves. If both have an odd number of points, we progress as before until one point is left on each side. Then we connect the two remaining points. Now for the other direction note that any $1$-manifold with boundary necessarily has an even number of boundary components. This requires that the two boundary manifolds together have an even number of points which proves the statement.\\

  Next we turn to $n=2$ and thus $1$-manifolds:

  As mentioned before, any $1$-dimensional closed manifold is a disjoint union of circles. Thus we have such manifolds on both sides. Now consider the two prominent cobordisms from the circle to the empty manifold or from the empty manifold to the circle

  \begin{center}
  \begin{tikzpicture}[tqft/.cd, cobordism/.style={draw},every lower boundary component/.style={draw},every incoming lower boundary component/.style={draw, dashed},every outgoing lower boundary component/.style={draw, solid}]
    \pic[tqft/cup, rotate = 90, at={(2,0)}, name = A];
    \pic[tqft/cap, rotate = 90, at={(2,0)}, name = B];
  \end{tikzpicture}
  \end{center}
  By "attaching" the respective cobordism to all circles on either side, one obtains a rather trivial cobordism which proves the statement.
\end{proof}
\end{lemma}

At this point it might already seem overdue to assign an unambigous assignment to the "left" and "right" side of the diagrams. While this will allow us to draw parallels to the direction of time evolution, it also allows to instill the structure of a category to cobordisms of a fixed dimension. Rather than transitions agnostic of direction, an orientation allows cobordisms to be "arrow-like" objects from one boundary manifold to another. To this end, we need to introduce some geometric preliminaries:

\begin{definition}
  Let $M$ be a smooth oriented manifold and $N$ a smooth closed $(n-1)$-dimensional manifold inheriting the orientation of $M$. Given some $p \in N$ and an element $v \in T_p N$ we call $v$ a \textbf{positive normal} if for any positive basis $[v_1, ...,v_{m-1}]$ of $T_p N$ one obtains a positive basis of $T_p M$ via $[v_1, ...,v_{m-1}, v]$.
\end{definition}

\begin{definition}[In- and Out-Boundaries]
  In the same setting as above, let $N$ be a connected component of $\partial M$. We call $N$ an \textbf{out-boundary}\index{Cobordism!out-boundary} if a positive normal of $N$ points outwards relative to $M$ and an \textbf{in-boundary}\index{Cobordism!in-boundary} if it points inwards relative to $M$.
\end{definition}

One convenient property of this definition of in- and out-boundaries is its invariance under orientation-reversal. Thus we can unambigously define

\begin{definition}[(Oriented) Cobordisms]
\label{def:cobordism_oriented}
  Let $n \in \NN$ $\Sigma_{in}$ and $\Sigma_{out}$ be closed $(n-1)$-dimensional manifolds. An \textbf{(oriented) cobordism}\index{Cobordism!oriented} between them is a tuple $(M, \imath_i, \imath_o)$ where $M$ is a smooth compact oriented $n$-dimensional manifold and the maps
  \begin{center}
    \begin{tikzcd}[sep = normal]
      \Sigma_{in} \arrow[r,"\imath_{in}"] & M & \arrow[l, "\imath_{out}"'] \Sigma_{out}
    \end{tikzcd}
  \end{center}
  are smooth orientation preserving diffeomorphisms onto the in- and out-boundary of $M$ respectively. When the diffeomorphisms are understood, we sometimes denote a cobordism by the "connecting" manifold $M$ or as $M \colon \Sigma_0 \Lra \Sigma_1$.
\end{definition}

From here on a \emph{cobordism} always refers to an oriented cobordism unless explicitely stated otherwise. Note that while the above definition preserves the idea of a manifold connecting two codimension $1$ manifolds, it gives us more freedom than the naive extension where $\Sigma_{in}$ and $\Sigma_{out}$ \textbf{are} the in- and out-boundary. Among other things one can construct a cobordism from a manifold to itself (see \cite{FrobAlgebraTQFT}). To illustrate the notion of a cobordism, we consider a prominent toy-model example:

\begin{example}[The Unit Interval]
\label{example:unit_interval}
  Let $I = [0,1]$ be the unit interval equipped with its standard orientation "$+$", which is passed down to the boundary points (and $0$-manifolds) $1$ and $0$. Thus $0$ is the in-boundary and $1$ an out-boundary. This can be generalized by taking any two positively oriented points $p_0, p_1$ and regard them as $0$-manifolds. Mapping $p_0$ to $0$ and $p_1$ to $1$ is an orientation-preserving diffeomorphism onto the boundary of $I$ thus forming the cobordism
  $$ p_0 \lra I \lra p_1 $$
  This can be generalized to any path-manifold $M$ such that there exists an orientation-preserving diffeomorphism $I \atob{\sim} M$. Playing with the other orientation of $I$ denoted by "$-$" gives us four possible diagrams for cobordisms between its boundary points:\index{Cobordism!1-cobordisms}
  \begin{center}
  \begin{tikzcd}[nodes={inner sep=-1pt}]
    + \ \bullet \arrow[r] & \bullet \ + & & - \ \bullet & \arrow[l] \bullet \ - \\
    - \ \bullet \arrow[dd, bend left=90, leftrightarrow] & & & & \bullet \ + \arrow[dd, bend right=90, leftrightarrow]\\
     & \emptyset & & \emptyset &\\
    + \ \bullet & & & & \bullet \ -
  \end{tikzcd}
  \end{center}
\end{example}

A very special and easy to construct type of cobordisms in any dimension is that of a cylinder. Intuitively one would try to generalize the idea of crossing a circle and an interval to an arbitrary closed oriented manifold $\Sigma$ and a path-manifold $M \atob{\sim} I$. Indeed this brings about the \emph{cylinder construction}:

\begin{example}[Cylinder Cobordism]
\label{example:cylinder}
  Let $\Sigma$ be a closed oriented $(n-1)$-manifold and $I$ the unit interval. The product manifold $\Sigma \times I$ has one in- and one out-boundary, namely $\Sigma \times \{0\}$ and $\Sigma \times \{1\}$. Taking the two evident maps
  \begin{center}
  \begin{tikzcd}[sep=normal]
    \imath_{in}: \Sigma \arrow[r, "\sim"] & \Sigma \times \{0\} \ \ \subset \ \ \Sigma \times I \ \ \supset \ \ \Sigma \times \{1\} & \arrow[l, "\sim"'] \Sigma : \imath_{out}
  \end{tikzcd}
  \end{center}
  we obtain a cobordism $(\Sigma \times I, \imath_{in}, \imath_{out})$ from $\Sigma$ to $\Sigma$. Now the same construction can be generalized to a cobordism between two $(n-1)$-manifolds $\Sigma_0$ and $\Sigma_1$ diffeomorphic to $\Sigma$. Further using the generalization of the unit interval to path-manifolds from \ref{example:unit_interval} we can build a cobordism $M \simeq I$ between $\Sigma_0$ and $\Sigma_1$.
\end{example}

The above example already shows that for two closed oriented $(n-1)$-dimensional manifolds there can be a multitude of cobordisms between them. Thus we have an increased interest in collecting certain equivalence classes of cobordisms, a cornerstone on our way to a category of cobordisms:

\begin{definition}[Equivalence of Cobordisms]
\label{def:cobordism_equivalence}
  Let $\Sigma_0$ and $\Sigma_1$ be two closed $(n-1)$-dimensional manifolds for some $n \in \NN$. Let further be $M$ and $N$ be two cobordisms between them. We call $M$ and $N$ \textbf{equivalent cobordisms}\index{Cobordism!equivalence} or just \textbf{equivalent}, if there exists an orientation-preserving diffeomorphism $\phi: M \atob{} N$ such that the following diagram commutes:
  \begin{center}
  \begin{tikzcd}[sep = normal]
    & N & \\
    \Sigma_0 \arrow[ur] \arrow[dr] & & \arrow[ul] \arrow[dl] \Sigma_1\\
    & M \arrow[uu, "\phi"] \arrow[uu, "\sim"'] &
  \end{tikzcd}
  \end{center}
\end{definition}
When naively constructing a category of cobordisms, one would take closed oriented $(n-1)$-dimensional manifolds as objects and cobordisms between them as objects. As we will soon see, this would lead to an ill-defined composition of cobordisms which would in turn prohibit a categorial treatment. With this in mind, the next subchapters will focus on building the fundaments to define the category of cobordisms.

\newpage
\subsection{Morse Functions and Glueing}
\label{subsec:cobordism_glueing}

This subchapter is devoted to an introduction to \emph{Morse Functions} and the \emph{glueing} of cobordisms. Both will play an important role in constructing and understanding the definition of the \emph{category of cobordisms} in the subsequent chapter \ref{subsec:cobordism_category}. Note that while it is possible to work with cobordisms without mentioning Morse functions or explicitely constructing the glueing of cobordisms, this path offers far more insight into the technical details that make a category of cobordisms possible. We start by giving a compact introduction to Morse functions, roughly following the lines of \cite{FrobAlgebraTQFT}:\\

In the previous chapter we discussed the cylinder construction using the unit interval $I = [0,1]$ \ref{example:unit_interval}, \ref{example:cylinder}. In this discussion it became clear, that for a smooth manifold $\Sigma$ functions of the form $\Sigma \times I \lra I$ that respect the boundaries, can play an important role. Generalizing this idea to general manifolds $M$ yields maps of the form

\begin{center}
\begin{tikzpicture}[tqft/.cd, cobordism/.style={draw},every lower boundary component/.style={draw},every incoming lower boundary component/.style={draw, dashed},every outgoing lower boundary component/.style={draw, solid}]
  \pic[tqft/reverse pair of pants, rotate=90, at={(0,.5)}, name=A];
  \node[pin={[pin distance=1.3cm]30:$M$}] at (.75,1.7) {};
  \draw[|-|] (-.1,-1) -- (1.9,-1) node[below] {$1$} node[below, xshift=-2.cm] {$0$} node[left, xshift=-.21cm] {\tiny{|}} node[left, xshift=-.21cm, yshift=-.3cm] {$t_1$} node[left, xshift=-1.21cm] {\tiny{|}} node[left, xshift=-1.21cm, yshift=-.3cm] {$t_0$};

  \draw[->] (1.5,1.3) -- (1.5,-.75);
  \draw[] (1.5, 1.3) node[left, xshift=.24cm] {$q \bullet$};

  \draw[->] (.51,1.5) -- (.51,-.75);
  \draw[] (1,1) node[left, xshift=-.25cm, yshift=.45cm] {$p \bullet$};
\end{tikzpicture}
\end{center}

It turns out that many of the topological characteristics of the manifold $M$ can be found within such a map. This sketches the rough outline of \emph{Morse theory} in differentiable topology. Since we will only use some of its basic notions, the reader is referred to \cite{Hirsch}
for a more thorough introduction to the field. We begin by remembering some notions of differential geometry:

\begin{definition}
  Let $M$ be a compact smooth manifold of dimension $m$ and $f\colon M \lra I$ a smooth map for $I \subset \RR$ a closed interval. A point $p \in M$ is called a \textbf{critical point}\index{Smooth map!critical point} if $Df(p) = 0$. If a point is not critical, it is called \textbf{regular}\index{Smooth map!critical point}. In the same manner, the image $p$ under $f$ is called a \textbf{critical value}\index{Smooth map!critical value} or \textbf{regular value}\index{Smooth map!regular value} if $p$ is critical or regular.
\end{definition}

In the above picture $p$ would be a critical point with corresponding critical value $t_0$. $q$ on the other hand would be a regular point and $t_1$ a regular value.

\begin{definition}
\label{def:index_nondeg}
  Let $M$ be a smooth manifold, $f \colon M \lra I$ a smooth map onto a closed interval $I \subset \RR$ and $p$ a critical point of $f$. We call $p$ \textbf{nondegenerate}\index{Smooth map!nondegenerate point} if for some coordinate chart $(U,x)$ about $p$ the Hessian of $f$ is nonsingular. Since the Hessian is a real symmetric matrix, we can utilize its real eigenvalues to define the \textbf{index}\index{Smooth map!index} of $f$ at $p$ to be the number of negative eigenvalues.
\end{definition}

The above definition can be better understood using tangible examples:

\begin{example}
\label{example:index_surface}
  Let $M$ be a surface ($2$-dimensional smooth manifold) and $f\colon M \lra I$ a smooth function. In this case, a nondegenerate critical point $p$ has index $0$ iff it is a local minimum. It has index $1$ iff it is a saddle point and index $2$ if it is a local maximum.
\end{example}

Now we have everything needed to define \emph{Morse functions}:

\begin{definition}[Morse functions]
\label{def:cobordism_morse}
  Let $M$ be a smooth manifold and $I$ an interval. We call a smooth map $f\colon M \lra I$ a \textbf{Morse function}\index{Smooth map!Morse function} if all its critical points are nondegenerate. If $M$ is a manifold with boundary, we further require
  $$ f^{-1}(\partial I) = \partial M $$
  and that both boundary points of $I$ are regular values of $f$.
\end{definition}

\begin{remark}
  Note that if $M$ is compact, which holds true for most of our previous and upcoming applications, the number of critical points of a Morse function is finite which makes it possible to rearrange the function such that all those critical points have disjoint images under $f$. For our purposes we will tacitly assume that such a rearrangement was executed.
\end{remark}

While a nice construct, we are not yet sure if Morse functions even exist in every or any of our upcoming applications. The following theorem presented by \textit{Hirsch}, \textit{chapter 6.1.2} of \cite{Hirsch}, does not only solve that concern, it even states that "most" functions are Morse functions:

\begin{theorem}[Existence of Morse Functions]
  For any manifold $M$, Morse functions form a dense subset of $C^s(M, I)$ where $2 \leq s \leq \infty$\index{Smooth map!Existence of Morse functions}.
\end{theorem}

For our purposes the above theorem can be summed up as "Morse functions always exist". Thus in upcoming applications we will simply take them for granted, using them whenever convenient. We will mainly use the notion of \emph{Morse functions} to formulate an ever more specialized definition of \emph{glueing} spaces together and investigating the contributing factors. Since we ultimately want to glue cobordisms together, thus smooth compact oriented manifolds, it is only natural to start at the low level of topological spaces:

\begin{definition}[Glueing of Topological spaces]
  Let $X,Y$ and $Z$ be topological spaces and $f\colon X \lra Y$, $g\colon X \lra Z$ be continous maps between them. We define the \textbf{glueing of $Y$ and $Z$ along $X$}\index{Glueing}, denoted by $Y \cop_{X} Z$ as the quotient space
  $$Y \cop Z / \sim \quad \text{where} \quad y \sim z \quad \text{iff} \quad  \exists x \in X \quad  s.t.  \quad f(x) = y \quad \& \quad g(x) = z$$
  We define the topology on $Y \cop_{X} Z$ using the two functions $f,g$: A set in $Y \cop_{X} Z$ is open, if its inverse image under $f$ and $g$ is open in $Y$ or $Z$ respectively. If the common boundary is understood, we will often denote the glueing $Y \cop_{X} Z$ by $YZ$. This notation will be adapted to smooth manifolds later on.
\end{definition}

One can show that the above construction, sometimes called the "pushout"\index{Glueing!pushout}, is in fact "unique" in a categorical sense: It turns out to be the colimit of the system $Y \longleftarrow X \lra Z$. Using the two evident maps $ Y \lra Y \cop_{X} Z \longleftarrow Z $ this unique property can be expressed as follows: For any commuting diagram of the form
\begin{center}
\begin{tikzcd}
  & P & \\
  Y \arrow[ur, "f^\prime"] & & \arrow[ul, "g^\prime"'] Z \\
  & \arrow[ul, "f"] X \arrow[ur, "g"']&
\end{tikzcd}
\end{center}
there exists a unique continous map $\phi \colon Y \cop_{X} Z \lra P$ such that the following diagram commutes:

\begin{center}
\begin{tikzcd}[sep=large]
  & Y \arrow[dr] \arrow[drr, bend left=30, "f^\prime"] & &  \\
  X \arrow[ur, "g"] \arrow[dr, "f"'] & & Y \cop_{X} Z  \arrow[r, dashed, "\phi"] & P\\
  & Z \arrow[ur] \arrow[urr, bend right=30, "g^\prime"'] & &
\end{tikzcd}
\end{center}

This unique map\index{Glueing!uniqueness} can be interpreted as glueing the two maps $f^\prime$ and $g^\prime$, which in turn allows us to state that any two continous maps that agree on the border between their domains can be glued together. For more material on this, see \textbf{A.3.7} in \cite{FrobAlgebraTQFT}.\\

We now aim to extend this notion of "glueing" to topological manifolds. The extra structure that we need to consider is thus the construction of a $C^0$-atlas on $Y \cop_{X} Z$ from the two $C^0$-atlases of $Y$ and $Z$. First we note that any point of $Y \cop_{X} Z$ that is not on $X$ is trivially covered by one of the two $C^0$-atlases. Now choosing a point in the common boundary $X$ and an open neighbourhood $U$ around it, we want to construct a $C^0$-chart $U \lra \RR^n$ about $p$. The main idea is to use the two intersections $U \cap Y$ and $U \cap Z$ to construct a new glued topological manifold using these open sets. One then constructs two charts
$$ f_0 \colon U \cap Y \lra \RR_+^n, \quad \quad f_1 \colon U \cap Z \lra \RR_-^n $$
to construct a chart $f:U \lra \RR^n$ as needed by using the universal property shown above. In the last steps, one shows that the choices of $f_0$ and $f_1$ do not influence the construction of $f$ by again using the universal property of such a map. A detailed discussion of this can be found in \cite[1.3.2]{FrobAlgebraTQFT}.\\

Our next goal is to generalize the glueing to smooth manifolds and ultimately to cobordisms. Again, we cut this discussion short and provide a compact summary of the results up until the \emph{Regular interval theorem}. An explicit and more extensive version can be found in \cite[1.3.4 - 1.3.6]{FrobAlgebraTQFT}:\\

When defining the glueing of smooth manifolds, it is only natural to first investigate some "toy model" examples. Using $1$-manifolds, one can show that the resulting glueing structure is not unique, that is, one can define multiple distinct smooth structures on the glued manifold. Since given two manifolds $M_0 , M_1$ one aims to define \emph{one} smooth structure on $M = M_0 M_1$, the discussion is expanded to smooth manifolds $X$ that are homeomorphic to $M$ via some map $\phi \colon M \atob{\sim} X $ which restricts to a diffeomorphism on the two manifolds $M_0, M_1$. One can then use $\phi$ to pull the maximal atlas of $X$ to $M$.\\

The thus described process can now be applied to cylinder cobordisms. This will be done in an explicit fashion following the construction seen in \cite[1.3.6]{FrobAlgebraTQFT}. Thus let $M_0 \colon \Sigma_0 \Lra \Sigma_1$ and $M_1 \colon \Sigma_1 \Lra \Sigma_2$ be two cobordisms of dimension $n \in \NN$ which are equivalent to cylinders along $\Sigma_1$ in the sense of \ref{example:cylinder}, that is there exist two homeomorphisms
$$ \phi_0 \colon M_0 \atob{\sim} \Sigma_1 \times [0,1], \quad \quad \phi_1 \colon M_1 \atob{\sim} \Sigma_2 \times [1,2] $$
As proposed before we need to find smooth $n$-manifold $X$ such that there exists a homeomorphisms $\phi \colon M_0 \cop_{\Sigma_1} M_1 \atob{\sim} X$ that restricts to a diffeomorphism on $M_0$ and $M_1$. A natural choice would be the cylinder $\Sigma_1 \times [0,2]$ and indeed, using the universal property, one can take the product of the maps $\phi_0$ and $\phi_1$ to define a homeomorphism
$$ \phi := \phi_0 \cop_{\Sigma_1} \phi_1 : M_0 \cop_{\Sigma_1} M_1 \lra \Sigma_1 \times [0,2] $$
This glued map can be visualized in the following diagram:

\begin{center}
\begin{tikzpicture}[tqft/.cd, cobordism/.style={draw},every lower boundary component/.style={draw},every incoming lower boundary component/.style={draw, dashed},every outgoing lower boundary component/.style={draw, solid}]
  \pic[tqft/cylinder to next, rotate=90, at={(0,2.5)}, name=A];
  \pic[tqft/cylinder to prior, rotate=90, at=(A-outgoing boundary 1)];

  \pic[tqft/cylinder, rotate=90, at={(0,1)}];
  \pic[tqft/cylinder, rotate=90, at={(2,1)}];
  \draw[|-|] (0,0) -- (4,0) node[xshift=-2.cm] {\tiny{|}} node[xshift=-4.cm] {\tiny{|}} node[xshift=-0.cm] {\tiny{|}} node[xshift=-4.cm, yshift=-.4cm]{$0$} node[xshift=-2.cm, yshift=-.4cm]{$1$} node[xshift=-0.cm, yshift=-.4cm]{$2$};

  \node[pin={[pin distance=.75cm]150:$M_0$}] at (1.,3.) {};
  \node[pin={[pin distance=.75cm]30:$M_1$}] at (3.,3.) {};

  \draw (-.5, 2.5) node[xshift=0.0cm, yshift=0.0cm] {$\Sigma_0$};
  \draw (4.5, 2.5) node[xshift=0.0cm, yshift=0.0cm] {$\Sigma_2$};
  \draw (2., 4.2) node[xshift=0.0cm, yshift=0.0cm] {$\Sigma_1$};

  \draw[->] (2,3) -- (2,1.5) node[xshift=.2cm, yshift=.75cm]{$\phi$};
  \draw[->] (.5,2) -- (.5,1.5) node[xshift=-.25cm, yshift=.35cm]{$\phi_0$};
  \draw[->] (3.5,2) -- (3.5,1.5) node[xshift=.25cm, yshift=.35cm]{$\phi_1$};

  \draw (.0, 1) node[xshift=-1.2cm, yshift=0.0cm] {$\Sigma_1 \times [0,2]$};

\end{tikzpicture}
\end{center}

Now since $\Sigma_1 \times [0,2]$ has a smooth structures that agrees with that of $\Sigma_1 \times [0,1]$ and $\Sigma_1 \times [1,2]$ respectively, we can use the pullback by $\phi$ to get a smooth structure on $M = M_0 M_2$. This discussion leads to a very powerful theorem which will allow us to move towards the glueing of general cobordisms:

\begin{theorem}[Regular Interval Theorem]\index{Cobordism!Regular Interval Theorem}
\label{theorem:regular_interval}
  Let $M \colon \Sigma_0 \Lra \Sigma_1$ be an $n$-cobordism and $f\colon M \lra [0,1]$ a Morse function without critical points such that $f^{-1}(0) = \Sigma_0$ and $f^{-1}(1) = \Sigma_1$. Let further $\pi \colon \Sigma_0 \times [0,1] \lra [0,1]$ be the obvious projection. Then there exists a diffeomorphism $\phi \colon \Sigma_0 \times [0,1] \lra M$ such that the following diagram commutes:
  \begin{center}
  \begin{tikzcd}[sep = large]
    \Sigma_0 \times [0,1] \arrow[dr, "\pi"'] \arrow[r, "\phi"] & M \arrow[d, "f"] \\
    & ~[0,1]
  \end{tikzcd}
  \end{center}
  Note that corresponding maps can be constructed for the cylinder $\Sigma_1 \times [0,1]$.
\begin{proof}
  A detailed proof of this statement is presented in \cite[6.2.2]{Hirsch}.
\end{proof}
\end{theorem}

Using the \emph{Regular Interval Theorem} we can further investigate the splitting of an arbitrary cobordism which comes close to their glueing:

\begin{lemma}
\label{lemma:decomposition_cobord}
  Let $M \colon \Sigma_0 \Lra \Sigma_1$ be an $n$-cobordism and $f \colon M \lra [0,1]$ a Morse function such that $f^{-1}(0) = \Sigma_0$ and $f^{-1}(1) = \Sigma_1$. By definition \ref{def:cobordism_morse} $f$ has no critical points on the boundary components $\Sigma_0, \Sigma_1$.
  Then there exists some $\epsilon > 0$ and a decomposition\index{Cobordism!decomposition} $M = f^{-1}([0,\epsilon]) f^{-1}([\epsilon, 1])$ such that $M_0 := f^{-1}([0,\epsilon])$ is diffeomorphic to the cylinder $\Sigma_0 \times [0,1]$. Similarly there exists another such decomposition for the boundary component $\Sigma_1$.

\begin{proof}
  Since $f$ has no critical points on the boundary components, there exists an $\epsilon > 0$ such that $f$ has no critical values on the interval $[0,\epsilon]$. Now if we restrict $f$ to $M_0 := f^{-1}([0,\epsilon])$ it is a Morse function without any critical points thus meeting the requirements of \ref{theorem:regular_interval}.
  Thus $M_0$ is diffeomorphic to the cylinder $\Sigma_0 \times [0,1]$. This induces a decomposition of $M$ into $M = M_0 M_1$. The composition for the boundary component $\Sigma_1$ follows in a similar fashion.
\end{proof}
\end{lemma}

Indeed this result gives us quite some ideas when it comes to glueing arbitrary cobordisms. Taking two cobordisms $M_0 \colon \Sigma_0 \Lra \Sigma_1$ and $M_1 \colon \Sigma_1 \Lra \Sigma_2$ and two Morse functions $f_0 \colon M_0 \lra [0,1], f_1 \colon M_1 \lra [1,2]$ we can choose $\epsilon > 0$ such that the intervals $[1-\epsilon,1]$ and $[1,1+\epsilon]$ contain no critical values of the respective Morse function.
Thus the two inverse images $f^{-1}([1-\epsilon,1])$ and $f^{-1}([1,1+\epsilon])$ are each diffeomorphic to cylinders by rule of \ref{theorem:regular_interval}. This brings us back to the glueing of cylinder cobordisms discussed above which lets us take the smooth structure of the cylinder cobordism. Thus we obtain a "glueing cobordism"\index{Cobordism!glueing cobordism}\index{Glueing!glueing cobordism} $M = M_0 M_1$. Now as discussed before, the thus defined smooth structure need not be unique. The following theorem and its proof, both of which can be found in \cite[8.2.1]{Hirsch}, makes use of isotops and is thus omited from this script. It basically states that the smooth structure defined on glueing cobordisms is unique up to diffeomorphism:

\begin{theorem}
\label{theorem:glueing_diffeo}
  Let $M_0$ and $M_1$ be two cobordisms and $\Sigma$ the out-boundary of $M_0$ and the in-boundary of $M_1$. Let $M_0 M_1 = M_0 \cop_\Sigma M_1$ be the glueing cobordism and $\alpha, \beta$ two smooth structures on $M_0 M_1$ both respecting the smooth structures of the single cobordisms. Then there exists a diffeomorphism
  $$ \phi \colon (M_0 M_1, \alpha) \lra (M_0 M_1, \beta) \quad \quad \text{such that} \quad \quad \phi_\Sigma = \id_\Sigma. $$
\end{theorem}

This theorem puts us in a good place on our way towards a category of cobordisms, since we now know that given any two cobordisms there exists a diffeomorphism class of glueing cobordisms. Further this diffeomorphism class is well-defined. As was hinted at in the closing statement of subsection \ref{subsec:cobordism_def}, the previous constructions of glueing cobordisms present the core of why we \emph{need} to use equivalence classes of cobordisms instead of just cobordisms as arrows for such a category. On this note, we end our technical discussion of the glueing of cobordisms and defer the investigation of composition of cobordism classes into the next subchapter which will thoroughly construct the symmetric monoidal category of $n$-cobordisms, $\cobn{n}$.

\newpage
\subsection{The Category $\cobn{n}$}
\label{subsec:cobordism_category}

This subchapter is dedicated to the construction and classification of the category of $n$-cobordisms, $\cobn{n}$, using the tools presented in \ref{sec:Categorial_Preliminaries}. As has been sufficiently announced, once we use equivalence classes, and thus diffeomorphism classes, of cobordisms as arrows and the previously presented glueing \ref{subsec:cobordism_glueing} as composition, we naturally happen upon the structure of a category.
Using \ref{subsec:Monoidal_categories} and \ref{subsec:braided_categories} we will see that $\cobn{n}$ is indeed a symmetric monoidal category thus carrying a rich algebraic structure. We again roughly follow the outlines of \cite{FrobAlgebraTQFT} while incorporating discussions and insights from \cite{Intro_TQFT}.\\

Let us consider two cobordisms $M_0$ and $M_1$ such that $\Sigma$ is the out-boundary of $M_0$ and the in-boundary of $M_1$. By rule of \ref{theorem:glueing_diffeo} there exists a smooth structure unique up to diffeomorphism. Since we ultimately want to end up with a category where cobordism classes make up the arrows and their in- and out-boundaries the objects, we need to show that this smooth structure does not depend on the particular choice of a representant. Thus consider the following:

\begin{construction}
  Let $M_0\colon \Sigma_0 \Lra \Sigma_1$ and $M_1\colon \Sigma_1 \Lra \Sigma_2$ be two cobordisms. Further let $\phi_0\colon M_0 \atob{\sim} M_0^\prime$ and $\phi_1\colon M_1 \atob{\sim} M_1^\prime$ be two diffeomorphisms respecting the boundary such that we obtain the following diagram:
  \begin{center}
  \begin{tikzcd}[sep=large]
    & M_0^\prime & & M_1^\prime & \\
    \Sigma_0 \arrow[ur] \arrow[dr] & & \arrow[ul] \arrow[dl] \Sigma_1 \arrow[ur] \arrow[dr] & & \arrow[ul] \arrow[dl] \Sigma_2 \\
    & M_0 \arrow[uu,"\phi_0"] \arrow[uu,"\sim"'] & & M_1 \arrow[uu,"\sim"] \arrow[uu,"\phi_1"'] &
  \end{tikzcd}
  \end{center}
  Now for the two glueings $M_0 M_1$ and $M_0^\prime M_1^\prime$ we utilize the previously discussed universal property of glueings adapted to the category of continous maps to obtain a glueing of $\phi_0$ and $\phi_1$ denoted by $\phi$ which is a homeomorphism that restricts to a diffeomorphism on each glued piece. Thus we obtain the following diagram:
  \begin{center}
  \begin{tikzcd}[sep=large]
    & M_0^\prime M_1^\prime & \\
    \Sigma_0 \arrow[ur] \arrow[dr] & & \arrow[ul] \arrow[dl] \Sigma_2 \\
    & M_0 M_1 \arrow[uu, "\phi"] &
  \end{tikzcd}
  \end{center}
  This puts us in a comfortable position: Both $M_0 M_1$ and $M_0^\prime M_1^\prime$ are cobordisms with smooth structures. Using the homeomorphism $\phi$, which is universal, we can define a smooth structure on $M_0^\prime M_1^\prime$ using that of $M_0 M_1$. While this structure might differ from the original one, \emph{Theorem} \ref{theorem:glueing_diffeo} tells us that the two need to be diffeomorphic. This shows that the glueing cobordism does indeed not depend on the choice of representants.
\end{construction}

We use the above construction as a motivation for the following definition:

\begin{definition}[Composition cobordism]
\label{def:cobordism_composition}
  Let $M_0 \colon \Sigma_0 \Lra \Sigma_1$ and $M_1 \colon \Sigma_1 \Lra \Sigma_2$ be two cobordism classes. We define their \textbf{composition cobordism}\index{Cobordism!composition} as the equivalence class of cobordisms $M_0 M_1 \colon \colon \Sigma_0 \Lra \Sigma_2$ defined by the glueing of cobordism classes.
\end{definition}

Given a composition rule for what are to be arrows, we follow the requirements of the definition of a category \ref{def:category_set} and check if the composition is indeed associative.

\begin{lemma}
\label{lemma:cobordism_associative}
  Given three cobordism classes and respective representants denoted by
  \begin{center}
  \begin{tikzcd}[sep=large]
    \Sigma_0 \arrow[r, Rightarrow, "M_0"] & \Sigma_1 \arrow[r, Rightarrow, "M_1"] & \Sigma_2 \arrow[r, Rightarrow, "M_2"] & \Sigma_3
  \end{tikzcd}
  \end{center}
  the order of composition does not matter, i.e. $(M_0 M_1)M_2 = M_0(M_1 M_2)$ in terms of equivalence classes of cobordisms.
\begin{proof}
  First note that there exists a natural isomorphism for the glueing of topological manifolds such that
  $$ \left( M_0 \cop_{\Sigma_1} M_1 \right) \cop_{\Sigma_2} M_2 \cong M_0 \cop_{\Sigma_1} \left( M_1 \cop_{\Sigma_2} M_2 \right) $$
  Details on this can be found in the previously referenced discussion of glueings in \cite{FrobAlgebraTQFT}. Now we defined a smooth structure on such glueings by taking the smooth structure on the cobordisms and appropriately replacing it around the connecting boundaries. This was done by using the \emph{Regular Interval Theorem} \ref{theorem:regular_interval} to use the smooth structure of a cylinder in a small neighbourhood of the boundaries $\Sigma_1$ and $\Sigma_2$. Since these neighbourhoods are necessarily disjoint and the smooth structure is unaffected elsewhere, the order we replace the structure about the boundaries does indeed not matter. This proves the statement.
\end{proof}
\end{lemma}

With associativity of cobordism composition, and thus composition of arrows, established, we only need an identity morphism. To no surprise the cylinder cobordism \ref{example:cylinder} is an excellent choice:

\begin{lemma}
\label{lemma:cobordism_unity}
  Let $M_0: \Sigma_0 \Lra \Sigma_1$ be a cobordism and denote by $C_0$ and $C_1$ a cylinder cobordism over $\Sigma_0$ and $\Sigma_1$ respectively. Then $C_0 M_0$ and $M_0 C_1$ belong to the same cobordism class as $M_0$\index{Cobordism!unity}.
\begin{proof}
  We conduct the proof for $C_0$, the proof for $C_1$ is completely analogous. We have previsouly shown that any cobordism can be decomposed into a cylinder composed with the remaining part \ref{lemma:decomposition_cobord}.
  Thus we decompose $M_0$ into the cobordism $M_{[0,\epsilon]} M_{[\epsilon,1]}$ where $\epsilon > 0$ is a fitting constant for the respective construction in \ref{lemma:decomposition_cobord}. Note by rule of this constructon, $M_{[0,\epsilon]}$ is diffeomorphic to a cylinder over $\Sigma_0$. Now we make us of the associativity \ref{lemma:cobordism_associative} to note
  $$ C_0 M_0 = C_0 (M_{[0,\epsilon]} M_{[\epsilon, 1]}) \cong (C_0 M_{[0,\epsilon]}) M_{[\epsilon, 1]} \cong M_{[0,\epsilon]} M_{[\epsilon, 1]} = M_0$$
  Since the above equation holds modulo diffeomorphisms, $C_0 M_0$ and $M_0$ belong to the same diffeomorphism class. Note that in the second step, we used the fact that composed cylinders again form a cylinder which can be contracted to form $M_{[0,\epsilon]}$ using a diffeomorphism.
\end{proof}
\end{lemma}

Now that we have established that cobordisms qualify as categorial arrows and know that cylinder cobordisms act as identity arrows, we present a short input on \emph{invertible cobordisms} and then proceed to define $\cobn{n}$.

\begin{definition}[Invertible Cobordisms]
\label{def:cobordism_invertible}
    Let $M \colon \Sigma_0 \Lra \Sigma_1$ be a cobordism. We call $M$ an \textbf{invertible cobordism}\index{Cobordism!invertible} if there exists a cobordism $M^{-1} \colon \Sigma_1 \Lra \Sigma_0$ such that $M M^{-1}$ is a cylinder cobordism over $\Sigma_1$ and $M^{-1} M$ is a cylinder cobordism over $\Sigma_0$.
\end{definition}

Now we can finally define the \emph{category of cobordisms}:

\begin{definition}[Category of Cobordisms]
  For any $n \in \NN$ we denote by $\cobn{n}$ the \textbf{category of n-dimensional cobordisms}\index{The category of!n-cobordisms}\index{Cobordism!category} using $(n-1)$-dimensional closed oriented manifolds as objects and diffeomorphism classes of cobordisms between them as arrows. The previously defined composition of cobordisms \ref{def:cobordism_composition} serves as the composition of arrows, associativity is granted by rule of \ref{lemma:cobordism_associative}. Further for any object we obtain an identity arrow by forming the cylinder cobordism \ref{lemma:cobordism_unity}.
\end{definition}

Note that $\cobn{n}$ is named after its arrows since they are of major interest. While the objects are "plain" $(n-1)$-manifolds, we have already seen how much work was needed to fit cobordisms into the suit of a categorial arrow. It will later turn out that the special case of $\cobtwo$ is actually completely classified by a certain set of "generator" cobordisms which underlines their significance.\\

To obtain some fruitful training with the tensor product and we prove some results about invertible cobordisms which will be used in the investigation of $\cobtwo$ \ref{subsec:cobordism_2D}. It turns out that the definition of an invertible cobordism can be used to draw a close connection to topological properties of the involved manifolds:

\begin{lemma}
\label{lemma:cobordism_connected}
  Let $M \colon \Sigma_0 \Lra \Sigma_1$ be an invertible cobordism in $\cobn{n}$ for some $n \in \NN$. If the manifold $M$ is connected, so are $\Sigma_0$ and $\Sigma_1$.
\begin{proof}
  Since $M$ is invertible, the composition of cobordisms $M M^{-1}$ is diffeomorphic to the cylinder $\Sigma_0 \times I$ with $I$ the unit interval. For a cylinder we know that every point lies in the same connected component as some point of the in-boundary which happens to be the in boundary of $M$ itself. Now since $M$ is connected, all points necessarily lie in the same connected component which makes $\Sigma_0 \times I$ connected. But this amounts to saying that the base of the cylinder, $\Sigma_0$, is connected. Repeating the same argument for $M^{-1}M$ diffeomorphic to $\Sigma_1 \times I$ proves the statement.
\end{proof}
\end{lemma}

Further digressing into invertible cobordisms, the disjoint union of cobordisms can be shown to reveal the invertibility of the factors:

\begin{lemma}
\label{lemma:cobordism_disjoint_invertible}
  Let $M \colon \Sigma_0 \Lra \Sigma_1$ and $N \colon \Theta_0 \Lra \Theta_1$ be two cobordisms in $\cobn{n}$ for some $n \in \NN$. If their disjoint union $L := M \coprod N \colon \Sigma_0 \coprod \Theta_0 \Lra \Sigma_1 \coprod \Theta_1$ is invertible, then so are $M$ and $N$.
\begin{proof}
  Since $L$ is invertible, there exists a cobordism $L^{-1}$ such that $L L^{-1}$ is diffeomorphic to the cylinder cobordism over $\Sigma_0 \coprod \Theta_0$ which is diffeomorphic to $(\Sigma_0 \times I) \coprod (\Theta_0 \times I)$. Thus $L L^{-1}$ necessarily has more than one connected component, which means $L^{-1}$ cannot be connected as a manifold by rule of \ref{lemma:cobordism_connected}. Using this, we split
  $$L^{-1} = M^\prime \coprod N^\prime \quad \text{where} \quad M^\prime \colon \Sigma_1 \Lra \Sigma_0 \quad \text{and} \quad N^\prime \colon \Theta_1 \Lra \Theta_0$$
  Now investigate these two factors separately: Consider $M M^\prime \colon \Sigma_0 \Lra \Sigma_0$. This cobordism is necessarily diffeomorphic to the cylinder over $\Sigma_0$ since
  $$L L^{-1} = MM^\prime \coprod NN^\prime \cong (\Sigma_0 \times I) \coprod (\Theta_0 \times I)$$
  The same argument reveals that $NN^\prime$ is diffeomorphic to the cylinder over $\Theta_0$. Now if we repeat this entire process for $L^{-1}L$ we see that $M^\prime$ and $N^\prime$ are indeed two-sided inverses of $M$ and $N$ respectively which proves that both are invertible cobordisms.
\end{proof}
\end{lemma}

As previously mentioned, we can further specialize the categorial type of $\cobn{n}$. Ultimately we want to arrive at a symmetric monoidal category making use of \emph{subsections} \ref{subsec:Monoidal_categories} and \ref{subsec:braided_categories}. We start by defining a monoidal structure on $\cobn{n}$, restricting to a strict monoidal structure by rule of the \textbf{Strictification Theorem} \ref{theo:strictification_theorem}. The symmetric structure will then naturally occur from the definition of the tensor product.\\

First we note that given two oriented manifolds $M_0, M_1$ one can take the disjoint union $M_0 \coprod M_1$ to form another smooth oriented manifold. Now for this "coproduct"\index{Cobordism!coproduct} manifold there exists a unique orientation such that both evident inclusion maps are orientation preserving. Note that the empty manifold $\emptyset$ serves as a "unity" for the coproduct since for any manifold $M_0$ we have $M_0 \coprod \emptyset \cong M_0$.\\

This product, which does in fact form a monoidal category together with the category of oriented manifolds, can be extended to $\cobn{n}$: Given two $n$-cobordisms $M_0 \colon \Sigma_0 \Lra \Sigma_1$ and $N_0 \colon \Theta_0 \Lra \Theta_1$ we can form another cobordism using the disjoint union (i.e. coproduct)
$$M_0 \coprod N_0 \colon \Sigma_0 \coprod \Theta_0 \Lra \Sigma_1 \coprod \Theta_1$$
Again note that the "empty cobordism" $\emptyset_n \colon \emptyset_{n-1} \Lra \emptyset_{n-1}$ acts as a unit with respect to the thus defined product of cobordisms. Now moving from cobordisms to cobordism classes, we note that, since the union is disjoint, choosing different representatives of a cobordism class yields a diffeomorphic cobordism. Thus the disjoint union of cobordism classes is indeed well-defined. Since the empty cobordism $\emptyset_n \colon \emptyset_{n-1} \Lra \emptyset_{n-1}$ still acts as a union, we obtain a strict monoidal category:

\begin{definition}[Monoidal Category of Cobordisms]
  The category $\cobn{n}$ naturally extends to a monoidal category by taking the coproduct $\coprod$ as its tensor product and the empty $(n-1)$-dimensional manifol as a unit element. Thus we obtain $\triple{\cobn{n}}{\coprod}{\emptyset_{n-1}}$\index{The monoidal category of!n-cobordisms}\index{Cobordism!monoidal category}.
\end{definition}

To define a braiding and ultimately a symmetric braiding on $\cobn{n}$ we need to further expand upon the previously discussed cylinder cobordism \ref{example:cylinder}:

\begin{example}[Cylinder Construction]
\label{example:cylinder_construction}
  The main idea of the cylinder construction is very simple: Given two closed $(n-1)$-manifolds $\Sigma_0$ and $\Sigma_1$ and a diffeomorphism $\phi \colon \Sigma_0 \lra \Sigma_1$ between them, we want to construct a cylinder cobordism from $\Sigma_0$ to $\Sigma_1$. We will see however that this assignment naturally embeds into the categorial setting.\\
  Consider the category of $(n-1)$-dimensional smooth manifolds denoted by $\cat{Man}_{n-1}$ with arrows diffeomorphisms between them. Given any such arrow $f \colon \Sigma_0 \lra \Sigma_1$ we can the cylinder $\Sigma_0 \times I$, where $I$ denoted the unit interval, and map $\Sigma_0$ onto its in-boundary via the identity map and $\Sigma_1$ into its out-boundary using $f^{-1}$. The analogous construction is valid for $\Sigma_1 \times I$ yielding two equivalent cobordisms:
  \begin{center}
  \begin{tikzcd}[]
    & \Sigma_1 \times I &  \\
    \Sigma_0 \arrow[dr, "\id"'] \arrow[ur, "f"]& & \arrow[ul, "\id"'] \arrow[dl, "f^{-1}"] \Sigma_1 \\
    & \Sigma_0 \times I \arrow[uu, "f \times \id"'] &
  \end{tikzcd}
  \end{center}
  We denote the thus defined cylinder cobordism \index{Cobordism!cylinder} by $C_f$. Since $f$ is a diffeomorphism, one immediately sees that $C_f$ is indeed an invertible cobordism in the sense of \ref{def:cobordism_invertible}. Further note that for any object $\Sigma$ of $\cat{Man}_{n-1}$ the identity diffeomorphism $\id_\Sigma \colon \Sigma \lra \Sigma$ induces the identity cobordism over $\Sigma$. Carefully glancing at the definition of a functor \ref{def:functors} we turn to the composition of induced cylinder cobordisms:\\

  Given two arrows in $\cat{Man}_{n-1}$ of the form $\Sigma_0 \atob{f} \Sigma_1 \atob{g} \Sigma_2$ we investigate the composition of the two induced cylinder cobordisms $C_f$ and $C_g$. Here we use the cobordisms defined via cylinders over $\Sigma_0$ and $\Sigma_1$ respectively. Thus note that
  \begin{align*}
    C_f C_g &= (\Sigma_0 \times [0,1], \id, f^{-1}) \ (\Sigma_1 \times [0,1], \id, g^{-1}) \\
    &\cong \ (\Sigma_0 \times [0,1], \id, f^{-1}) \ (\Sigma_0 \times [0,1], f^{-1}, f^{-1} \circ g^{-1})\\
    &= (\Sigma_0 \times [0,1], \id, f^{-1} \circ g^{-1}) = C_{g \circ f}
  \end{align*}
  Putting all the properties of the \emph{cylinder construction} together we obtain a contravariant functor $\cat{C}_n \colon \cat{Man}_{n-1} \lra \cobn{n}$ that assigns to each diffeomorphism a cylinder cobordism.
\end{example}

Since we work with equivalence classes of cobordisms, it is important to know which diffeomorphisms are mapped onto the same object in $\cobn{n}$ by the functor $\cat{C}_n$. The following result follows the lines of \textit{proposition 1.3.23} of \cite{FrobAlgebraTQFT}:

\begin{prop}
\label{prop:cobordism_homotopy}
  Two diffeomorphisms $\phi \colon \Sigma_0 \lra \Sigma_1 $ and $ \psi \colon \Sigma_0 \lra \Sigma_1$ induce the same cobordism $\Sigma_0 \Lra \Sigma_1$ iff they are smoothly homotopic.
\begin{proof}
  Now if we have $\phi$ and $\psi$ inducing the same cobordism class, there exists an equivalence of cobordisms $M \colon \Sigma_0 \times I \lra \Sigma_1 \times I$ which we can compose with the evident projection $\pr_{\Sigma_1}$ to obtain a homotopy between $\phi$ and $\psi$. The converse direction requires a bit more thought:\\
  $\phi, \psi$ being homotopic means that there exists a smooth map $\Phi \colon \Sigma_0 \times I \lra \Sigma_1$ such that $\phi = \Phi(\cdot,0)$ and $\psi = \Phi(\cdot, 1)$. This amounts to saying that the following diagram commutes:
  \begin{center}
  \begin{tikzcd}[]
    & \Sigma_0 \times I \arrow[dd, "\Phi \times \pr_I"] &  \\
    \Sigma_0 \arrow[dr, "\phi"'] \arrow[ur, "\id"] & & \arrow[ul, "\id"'] \arrow[dl, "\psi"] \Sigma_0 \\
    & \Sigma_1 \times I &
  \end{tikzcd}
  \end{center}
  where we understand all maps but $\Phi$ as diffeomorphisms onto the in- and out-boundary respectively. To see that this diagram does indeed induce two equivalenc cobordisms for $\phi$ and $\psi$, compose the right part with $\psi^{-1} \colon \Sigma_1 \lra \Sigma_0$ to obtain
  \begin{center}
  \begin{tikzcd}[]
    & \Sigma_0 \times I \arrow[dd, "\Phi \times \pr_I"] &  \\
    \Sigma_0 \arrow[dr, "\phi"'] \arrow[ur, "\id_{\Sigma_0}"] & & \arrow[ul, "\psi^{-1}"'] \arrow[dl, "\id_{\Sigma_1}"] \Sigma_1 \\
    & \Sigma_1 \times I &
  \end{tikzcd}
  \end{center}
  But this is nothing else but stating that $\Phi \times \pr_I$ is an equivalence of the two resulting cobordisms $(\Sigma_0 \times I,\phi, \id_{\Sigma_1})$ and $(\Sigma_1 \times I, \id_{\Sigma_0}, \psi^{-1})$. Thus homotopic maps indeed induce representants of the same cobordism class which proves the statement.
\end{proof}
\end{prop}

The following discussion of the braiding that naturally occurs on $\cobn{n}$ will explicitely make use of the \emph{cylinder construction}:\\

Currently we work with the monoidal category $\triple{\cobn{n}}{\coprod}{\emptyset_n}$ that uses the adapted coproduct on cobordisms as a tensor product. Taking one step back and working with smooth $(n-1)$-dimensional manifolds $\Sigma, \Theta$ we immediately see that there exists a "twist" map
$$\gamma_{\Sigma \Theta} \colon \Sigma \coprod \Theta \lra \Theta \coprod \Sigma$$
defined by swapping the two factors of the coproduct. As such it clearly defines a diffeomorphism. Now using the cylinder construction \ref{example:cylinder_construction} we can assign an $n$-cobordism to this diffeomorphism denoted by
$$ \cat{T}_{\Sigma \Theta} := \cat{C}_{\gamma{\Sigma \Theta}} \colon \Sigma \coprod \Theta \Lra \Theta \coprod \Sigma $$
When the involved manifolds are understood, or once we refer to the map assigning to two manifolds their twist cobordism, we simply write $\cat{T}$ for the \textbf{twist cobordism}\index{Cobordism!twist}. The above secenario put into two dimensions can be used to visualize the twist cobordism as
\begin{center}
\begin{tikzpicture}[tqft/.cd, cobordism/.style={draw},every lower boundary component/.style={draw},every incoming lower boundary component/.style={draw, dashed},every outgoing lower boundary component/.style={draw, solid}]

  \pic[tqft/cylinder to next, rotate=90, at={(0,0)}];
  \pic[tqft/cylinder to prior, rotate=90, at={(0,1)}];

  \draw (-.5,-.1) node[xshift=0.0cm, yshift=0.0cm] {$\Theta$};
  \draw (-.5,1.1) node[xshift=0.0cm, yshift=0.0cm] {$\Sigma$};

  \draw (2.5,1.1) node[xshift=0.0cm, yshift=0.0cm] {$\Theta$};
  \draw (2.5,-.1) node[xshift=0.0cm, yshift=0.0cm] {$\Sigma$};

\end{tikzpicture}
\end{center}
Note that the crossing of the two cylinders is intended to avoid any resemblance of cylinders passing under or over each other. First of all, this clearly gives us a braided monoidal category \ref{def:monoidal_braided} $\triple{\cobn{n}}{\coprod}{\emptyset_n, \cat{T}}$\index{Cobordism!braided monoidal category}\index{The braided monoidal category of!n-cobordisms}, naturality of the twist stemming from the vivid examples in $2D$. Since all cobordisms defined via the cylinder construction are in fact invertible, the twist map $\cat{T}$ does in particular classify as a symmetric braiding \ref{def:monoidal_symmetric}:

\begin{definition}[Symmetric Monoidal Category of Cobordisms]
  The monoidal category of cobordisms given by $\triple{\cobn{n}}{\coprod}{\emptyset_n}$ extends to a symmetric monoidal category by using the twist cobordism $\cat{T}$ as its symmetric braiding. We thus obtain the symmetric monoidal category $\pair{\cobn{n}}{\cat{T}}$\index{Cobordism!symmetric monoidal category}\index{The symmetric monoidal category of!n-cobordisms}.
\end{definition}

\begin{example}
  Consider $\cobn{n}$ and let $k \in \NN$ with $k < n$. Now let $\Sigma$ be a closed compact oriented manifold of dimension $k$. Then the cartesian product of manifolds given by $(\ins) \times \Sigma$ induces a symmetric monoidal functor
  $$ (\ins) \times \Sigma \colon \cobn{n-k} \lra \cobn{n} $$
  This can be quickly verified by checking the properties given in \ref{def:monoidal_braided_functor}.
\end{example}

This concludes the chapter on the category of cobordisms. We could indeed use the full extent of the presented introductory material on categories discovering that cobordisms form a symmetric monoidal category, thus bearing a lot of additional structure. The precise classification will reappear when defining and discussing TQFTs in chapter \ref{sec:TQFT}. Using the notions presented in the previous pages the next subsection investigates the category $\cobtwo$ which will serve as groundwork for $2$-dimensional TQFTs.

\newpage
\subsection{Cobordisms in two Dimensions}
\label{subsec:cobordism_2D}

There are quite some reasons to put an increased focus on the category of cobordisms $\cobn{n}$ in two dimensions. This is no particular surprise, since looking at "small" dimensions is a frequently used and, in many cases, very powerful tool to learn about higher, often more complex dimensions. It is indeed very hard to describe $\cobn{n}$ for $n \geq 3$:\\

Considering the amount of possible $2$-dimensional closed oriented and non-diffeomorphic manifolds, that is oriented surfaces, this may already seem evident. However the classification of $n$-dimensional oriented surfaces plays a more important role looking at the arrows instead of the objects of $\cobn{n}$. Note that closed $2$-dimensional surfaces can indeed be completely classified while higher dimensions did not yet reveal similar properties. With this in mind, the chapter will work towards a complete classification, of $\cobtwo$ using a finite number of generator cobordisms. This amounts to finding a set of generating arrows in the sense of \ref{def:generating_set}.\\

Note that our restriction to "generating" any possible cobordism stems from another previously mentioned advantage of $\cobtwo$\index{Cobordism!2-cobordisms}. All of its objects are diffeomorphic to disjoint unions of the circle. This in particular enables the stunningly powerful visual representation of operations like composition and disjoint union that has been previously used. Note that while this can be captured in the categorial concept of the \emph{skeleton of a category} \ref{def:skeleton}, we are not \emph{yet} in the position to identify it. Since we again roughly follow the lines of \cite{FrobAlgebraTQFT}, we begin by providing central theorem whose proof and consequences serve as the raison d'être for the entire chapter.

\begin{theorem}[Classification of $\cobtwo$]
\label{theorem:cob2_classification}
  Using the composition of cobordisms and the coproduct, the following six cobordisms form a set of generators \ref{def:generating_set}\index{Cobordism!2-cobordisms!generators} for the (skeleton of the) symmetric monoidal category $\cobtwo$:
  \begin{center}
  \begin{tikzpicture}[tqft/.cd, cobordism/.style={draw},every lower boundary component/.style={draw},every incoming lower boundary component/.style={draw, dashed},every outgoing lower boundary component/.style={draw, solid}]

    \pic[tqft/cup, rotate=90, at={(2,.5)}];
    \pic[tqft/cap, rotate=90, at={(2,.5)}];
    \pic[tqft/cylinder, rotate=90, at={(5,.5)}];
    \pic[tqft/reverse pair of pants, rotate=90, at={(8,-.5)}];
    \pic[tqft/pair of pants, rotate=90, at={(11,.5)}];
    \pic[tqft/cylinder to next, rotate=90, at={(14,0)}];
    \pic[tqft/cylinder to prior, rotate=90, at={(14,1)}];

  \end{tikzpicture}
  \end{center}
  These cobordisms will further be called, by order of appearance in the above picture, the \textbf{cup}, \textbf{cap}, \textbf{cylinder} or \textbf{identity}, \textbf{copants}, \textbf{pants} and \textbf{twist} cobordism.
\end{theorem}

An interesting bit we have already proved for any cobordism is the possible decomposition into two pieces where one is diffeomorphic to a cylinder cobordism \ref{lemma:decomposition_cobord}. Note that this is not really a strong statement in the categorial context since the identity arrow is somewhat "understood" or rather trivially included; strictly speaking one could even exclude it from the above theorem. In this investigation however we will explicitely include the identity cobordism.\\

The proof of the classification theorem \ref{theorem:cob2_classification} has two popular approaches both of which can be found in \textit{chapter 1.4} of \cite{FrobAlgebraTQFT}. In this chapter we have chosen the ansatz using \emph{Morse Theory} to stay aligned with previous results like \ref{theorem:regular_interval} and \ref{lemma:decomposition_cobord}. This also significantly reduces the preliminary overhead. Note that from here on, we refer to objects and arrows of $\cobtwo$ if not noted otherwise.\\

The main idea of the Morse-theoretic proof is to assign to every cobordism a Morse-function, split its codomain into pieces with at most one critical value and then investigate the inverse of these pieces. It will turn out that the critical points on the cobordisms together with their index \ref{def:index_nondeg} will fully classify the type at hand. Again, we have already seen a manifestation of this principle in \ref{theorem:regular_interval}. Here we had cobordisms with Morse functions such that they had no critical points. This induced a diffeomorphism relation of \emph{one} of the boundaries which prevents us from concluding that such a cobordism is the cylinder (i.e. identity) cobordism. To fully classify this case, we need more groundwork:

\begin{lemma}
\label{lemma:cobordism_same_connected}
  Let $M \colon \Sigma_0 \Lra \Sigma_1$ be an arrow of $\cobn{n}$ for some $n \in \NN$. If $M$ is invertible, $\Sigma_0$ and $\Sigma_1$ have the same number of connected components.
\begin{proof}
  If $M$ is connected, we can simply use \ref{lemma:cobordism_connected} which proves the claim. Now if $M$ is not connected, we would like to "reorder" its boundaries in such a way that it can be represented as a disjoint union of cobordisms. This would throw us into the waters of \ref{lemma:cobordism_disjoint_invertible} thus proving that every factor is invertible. Inductively repeating this process, using \ref{lemma:cobordism_connected} when ariving at connected cobordisms, would prove that $\Sigma_0$ and $\Sigma_1$ have the same number of connected components. Now this "reordering" of boundaries is nothing but a fitting repeated application of the twist cobordisms $\cat{T}$ which serves to prove the statement.
\end{proof}
\end{lemma}

We will from now on refer to cobordisms composed of identities and twist cobordisms as \textbf{permutation cobordisms}\index{Cobordism!2-cobordisms!permutation}. Turning back to $\cobtwo$ we apply this more general result to arrive at a powerful conclusion for its objects.

\begin{prop}
\label{prop:cobordism_diffeo}
  Any two closed oriented $1$-manifold $\Sigma_0$ and $\Sigma_1$ are diffeomorphic iff there is an invertible cobordism between them.
\begin{proof}
  Let $\phi \colon \Sigma_0 \lra \Sigma_1$ be a diffeomorphism. We can use the cylinder construction \ref{example:cylinder_construction} to form a cobordism between them. Since the same construction is valid for the inverse map, we obtain an invertible cobordism in $\cobtwo$. If on the other hand we have an invertible cobordism $M \colon \Sigma_0 \Lra \Sigma_1$, we remind ourselves that $\Sigma_0 \cong \Sigma_1$ amounts to both having the same number of connected components since closed oriented $1$-manifolds are formed by disjoint unions of circles. However this holds true due to the previously proven lemma \ref{lemma:cobordism_same_connected} which proves the statement.
\end{proof}
\end{prop}

This tells us something very intricate about $\cobtwo$: Two objects of $\cobtwo$ lie in the same isomorphism class of objects, i.e. objects for which there exist isomorphism arrows \ref{def:object_iso} between them, iff they have the same number of connected components. This finally enables us to unambigously define the skeleton of $\cobtwo$ which will serve for our further discussions:

\begin{notation}
\label{nota:CobSkeleton}
  We now know that we can uniquely classify an isomorphism (in this case diffeomorphism) class of objects in $\cobtwo$ with an integer $n \geq 0$ indicating the number of connected components, i.e. the number of disjoint circles. One can thus define a \textbf{skeleton} of $\cobtwo$\index{Cobordism!2-cobordisms!skeleton} via the full subcategory $\{0,1,2, ...\}$ with arrows all arrows between these objects. We will henceforth use the skeleton of $\cobtwo$ instead of the entire category while, in an abuse of notation, denoting it by $\cobtwo$. Note that the set of generators in the \emph{Classification Theorem} \ref{theorem:cob2_classification} technically refers to the skeleton of $\cobtwo$ which of course also generates all cobordisms in the full category.
\end{notation}

This also makes intuitive sense if we inspect the pictorial representation of $\cobtwo$ used thus far. As an example, the following cobordism is clearly invertible:
\begin{center}
\begin{tikzpicture}[tqft/.cd, cobordism/.style={draw},every lower boundary component/.style={draw},every incoming lower boundary component/.style={draw, dashed},every outgoing lower boundary component/.style={draw, solid}]

  \pic[tqft/cylinder to next, rotate=90, at={(0,0)}];
  \pic[tqft/cylinder to prior, rotate=90, at={(0,1)}];

  \draw (-.5,.5) node[xshift=0.0cm, yshift=0.0cm] {$\Sigma_0$};
  \draw (2.5,.5) node[xshift=0.0cm, yshift=0.0cm] {$\Sigma_1$};

  \draw (3.5,.5) node[xshift=0.0cm, yshift=0.0cm] {$\Lra$};

  \draw (4.5,.5) node[xshift=0.0cm, yshift=0.0cm] {$\Sigma_0$};
  \draw (9.5,.5) node[xshift=0.0cm, yshift=0.0cm] {$\Sigma_0$};

  \pic[tqft/cylinder to next, rotate=90, at={(5,0)}];
  \pic[tqft/cylinder to prior, rotate=90, at={(5,1)}];
  \pic[tqft/cylinder to next, rotate=90, at={(7,0)}];
  \pic[tqft/cylinder to prior, rotate=90, at={(7,1)}];

  \draw (10.5,.5) node[xshift=0.0cm, yshift=0.0cm] {$\amp$};

  \draw (11.5,.5) node[xshift=0.0cm, yshift=0.0cm] {$\Sigma_1$};
  \draw (16.5,.5) node[xshift=0.0cm, yshift=0.0cm] {$\Sigma_1$};

  \pic[tqft/cylinder to next, rotate=90, at={(12,0)}];
  \pic[tqft/cylinder to prior, rotate=90, at={(12,1)}];
  \pic[tqft/cylinder to next, rotate=90, at={(14,0)}];
  \pic[tqft/cylinder to prior, rotate=90, at={(14,1)}];

\end{tikzpicture}
\end{center}

while this one clearly does not induce two cylinder cobordisms when composed with its "naive inverse":

\begin{center}
\begin{tikzpicture}[tqft/.cd, cobordism/.style={draw},every lower boundary component/.style={draw},every incoming lower boundary component/.style={draw, dashed},every outgoing lower boundary component/.style={draw, solid}]

  \pic[tqft/cylinder to next, rotate=90, at={(0,0)}];
  \pic[tqft/cap, rotate=90, at={(0,0)}];

  \draw (-.5,.0) node[xshift=0.0cm, yshift=0.0cm] {$\Sigma_0$};
  \draw (2.5,.5) node[xshift=0.0cm, yshift=0.0cm] {$\Sigma_1$};

  \draw (3.5,.5) node[xshift=0.0cm, yshift=0.0cm] {$\Lra$};

  \draw (4.5,.0) node[xshift=0.0cm, yshift=0.0cm] {$\Sigma_0$};
  \draw (9.5,.0) node[xshift=0.0cm, yshift=0.0cm] {$\Sigma_0$};

  \pic[tqft/cylinder to next, rotate=90, at={(5,0)}];
  \pic[tqft/cap, rotate=90, at={(5,0)}];
  \pic[tqft/cylinder to prior, rotate=90, at={(7,1)}];
  \pic[tqft/cup, rotate=90, at={(7,0)}];

  \draw (10.5,.5) node[xshift=0.0cm, yshift=0.0cm] {$\amp$};

  \draw (11.5,.5) node[xshift=0.0cm, yshift=0.0cm] {$\Sigma_1$};
  \draw (16.5,.5) node[xshift=0.0cm, yshift=0.0cm] {$\Sigma_1$};

  \pic[tqft/cylinder to prior, rotate=90, at={(12,1)}];
  \pic[tqft/cup, rotate=90, at={(12,0)}];
  \pic[tqft/cylinder to next, rotate=90, at={(14,0)}];
  \pic[tqft/cap, rotate=90, at={(14,0)}];

\end{tikzpicture}
\end{center}

We can in fact prove the lingering assumption that only disjoint unions of the twist $\cat{T}$ and the identity cobordism form invertible cobordisms:

\begin{lemma}
\label{lemma:cobordism_invertible_2_diffeo}
  The only invertible cobordisms in $\cobtwo$ are induced by diffeomorphisms. Further every invertible cobordism is a  \textit{permutation cobordisms}.
\begin{proof}
  Putting together \ref{lemma:cobordism_same_connected} and \ref{prop:cobordism_diffeo} we note that given two objects $\Sigma_0, \Sigma_1 \in \cobtwo$ we have the following equality: If there exists an invertible cobordism between them, $\Sigma_0$ and $\Sigma_1$ have the same number of connected components. Thus they are diffeomorphic since all closed oriented $1$-manifolds are disjoint unions of a number of circles equal to their connected components. But this diffeomorphism leads back to an invertible cobordism between them. Thus all invertible cobordisms are induced by diffeomorphisms.\\

  Taking this one step further, we want to utilize homotopy \ref{prop:cobordism_homotopy} to classify the different types of possible cobordisms. However every orientation-preserving diffeomorphism between circles is necessarily smoothly homotopic to the identity. Thus, up to homotopy, every diffeomorphism between objects of $\Sigma_0$ and $\Sigma_1$ needs to be a permutation of their respective connected components. Since such diffeomorphisms induce the invertible \textit{permutation cobordisms}, every invertible cobordism is of this type. This includes the identity cobordism that is induced by the identity permutation.
\end{proof}
\end{lemma}

Now, finally, we can inject all of the above results into a profound treatment of the decomposition of arbitrary cobordisms in $\cobtwo$ as was promised after stating the \emph{Classification Theorem} \ref{theorem:cob2_classification}. First, we combine \ref{theorem:regular_interval} and \ref{lemma:cobordism_invertible_2_diffeo} to state the following corollary:

\begin{corollary}
\label{corollary:no_crit_perm}
  If a cobordism in $\cobtwo$ admits a Morse function without critical points, it is equivalent to a \emph{permutation cobordism}.
\end{corollary}

Next, we again refer to \emph{Hirsch} for a Morse theoretic tool:

\begin{lemma}[Hirsch \ref{Hirsch}, 9.3.3]
\label{lemma:two_discs_missing}
  Let $M$ be a compact connected orientable surface admitting a Morse function having only one critical point of index $1$. Then $M$ is a disk with two holes.
\end{lemma}

The above lemma can be rephrased and visualized for our setting. A connected cobordism in $\cobtwo$ which admits a Morse functions with exactly one critical point $p$ of index $1$ is diffeomorphic to one of the following two cobordisms:

\begin{center}
\begin{tikzpicture}[tqft/.cd, cobordism/.style={draw},every lower boundary component/.style={draw},every incoming lower boundary component/.style={draw, dashed},every outgoing lower boundary component/.style={draw, solid}]

  \pic[tqft/reverse pair of pants, rotate=90, at={(0,.5)}, name=A];
  \draw[|-|] (-.1,-.5) -- (1.9,-.5) node[below] {$1$} node[below, xshift=-2.cm] {$0$} node[left, xshift=-1.21cm] {\tiny{|}} node[left, xshift=-1.21cm, yshift=-.3cm] {$t$};
  \draw[] (1,1) node[left, xshift=-.25cm, yshift=.45cm] {$p \bullet$};

  \draw[] (4,1) node[]{or};

  \pic[tqft/pair of pants, rotate=90, at={(6,1.5)}, name=B];
  \draw[|-|] (6,-.5) -- (8,-.5) node[below] {$1$} node[below, xshift=-2.cm] {$0$} node[left, xshift=-.21cm] {\tiny{|}} node[left, xshift=-.21cm, yshift=-.3cm] {$t$};
  \draw[] (8,1) node[left, xshift=-.03cm, yshift=.45cm] {$\bullet p$};
\end{tikzpicture}
\end{center}

This in turn allows us to finally conclude the proof of the \emph{Classification Theorem}:

\begin{proofof}[Classification Theorem, \ref{theorem:cob2_classification}]
  Let $M \colon \Sigma_0 \Lra \Sigma_1$ be a cobordism in $\cobtwo$ and $f \colon M \lra I$ a Morse function such that $f^{-1}(0) = \Sigma_0$ and $f^{-1}(1) = \Sigma_1$. We required in \ref{def:cobordism_morse} that both boundary points of $I$ are regular values of $f$ and that all critical points of $f$ have disjoint images under $f$.
  Further, the number of critical points of $f$ is finite, say $n$, and let $k \geq n$. This allows us to subdivide the interval $I$ into subintervals using $k+1$ strictly increasing regular values $x_0, ..., x_k$ such that for any $i = 0, ..., k-1$ the subinterval $[x_i, x_{i+1}]$ contains \emph{at most} one critical value. Using these subintervals, we also decompose our cobordism $M$ using the glueing of cobordisms into
  $$ M = \coprod_{i=0}^{k-1} M_{[x_i, x_{i+1}]} := \coprod_{i=0}^{k-1} f^{-1} \left([x_i, x_{i+1}]\right) $$
  This puts us in the comfortable position to inspect every piece by itself. Thus for $j$ arbitrary, we inspect $M_{[x_j, x_{j+1}]}$. This piece of $M$ contains at most one critical point $p$ contained in only one of the possibly many connected components of $M$. Since all other components have no critical points, they are permutation cobordisms by rule of \ref{corollary:no_crit_perm} and thus particularly generated by the twist and identity cobordism. This completely exhausts those pieces that contain no critical point.\\

  Now assume that there is a critical point $p$. We take only the connected component that contains and enote it by $C_p$. Comparing with \ref{example:index_surface} we see that if $p$ has index $0$, it is a local minimum, thus $C_p$ is diffeomorphic to

  \begin{center}
  \begin{tikzpicture}[tqft/.cd, cobordism/.style={draw},every lower boundary component/.style={draw},every incoming lower boundary component/.style={draw, dashed},every outgoing lower boundary component/.style={draw, solid}]
    \pic[tqft/cap, rotate=90, at={(0,0)}, name=A];
    \draw[] (0,0) node[left, xshift=1.75cm, yshift=-.05cm] {$p \ \bullet$};
  \end{tikzpicture}
  \end{center}

  If it has index $2$, $p$ is a local maximum and $C_p$ is diffeomorphic to

  \begin{center}
  \begin{tikzpicture}[tqft/.cd, cobordism/.style={draw},every lower boundary component/.style={draw},every incoming lower boundary component/.style={draw, dashed},every outgoing lower boundary component/.style={draw, solid}]
    \pic[tqft/cup, rotate=90, at={(0,0)}, name=A];
    \draw[] (0,0) node[left, xshift=1.09cm, yshift=-.05cm] {$\bullet \ p$};
  \end{tikzpicture}
  \end{center}

  The last remaining case is the saddle with index $1$ which, using \ref{lemma:two_discs_missing}, corresponds to $C_p$ being diffeomorphic to

  \begin{center}
  \begin{tikzpicture}[tqft/.cd, cobordism/.style={draw},every lower boundary component/.style={draw},every incoming lower boundary component/.style={draw, dashed},every outgoing lower boundary component/.style={draw, solid}]

    \pic[tqft/reverse pair of pants, rotate=90, at={(0,.5)}, name=A];
    \draw[] (1,1) node[left, xshift=-.25cm, yshift=.45cm] {$p \bullet$};

    \draw[] (4,1) node[]{or};

    \pic[tqft/pair of pants, rotate=90, at={(6,1.5)}, name=B];
    \draw[] (8,1) node[left, xshift=-.03cm, yshift=.45cm] {$\bullet p$};
  \end{tikzpicture}
  \end{center}

  Since $j$ was arbitrary, we have completely decomposed all pieces and thus any possible cobordism in $\cobtwo$ into the generators given in $\ref{theorem:cob2_classification}$ which proves the statement. \qed
\end{proofof}

Now that we have proven \ref{theorem:cob2_classification}, we can turn towards relations between different connected cobordisms which can be trivially generalized to non-connected ones using the disjoint union. Most of these relations capture the symmetric and monoidal property of $\cobtwo$. Again we will represent them pictorially to simplify notation and support intuition.

\newpage
\begin{theorem}[Generator relations]\index{Cobordism!2-cobordisms!generator relations}
\label{theorem:cobordism_relations}
The following relations hold true:
\begin{itemize}
  \item \textbf{The identity relations}
  \begin{equation}\label{eq:identity_relations}
    \begin{tikzpicture}[scale=.5, every tqft/.style={transform shape},tqft/.cd, cobordism/.style={draw},every lower boundary component/.style={draw},every incoming lower boundary component/.style={draw, dashed},every outgoing lower boundary component/.style={draw, solid}]

      \pic[tqft/cap, rotate=90, at={(9-5.5,.0)}];
      \pic[tqft/cylinder, rotate=90, at={(11-5.5,.0)}];
      \draw (14-5.5,0) node[xshift=0.0cm, yshift=0.0cm] {$=$};
      \pic[tqft/cap, rotate=90, at={(13-5,.0)}];

      \pic[tqft/cup, rotate=90, at={(20-2,.0)}];
      \pic[tqft/cylinder, rotate=90, at={(18-2,.0)}];
      \draw (21-6.5,0) node[xshift=0.0cm, yshift=0.0cm] {$=$};
      \pic[tqft/cup, rotate=90, at={(22-9,.0)}];


      \pic[tqft/cylinder, rotate=90, at={(-1.0-3,-2.0)}];
      \pic[tqft/cylinder, rotate=90, at={(-1.0-3,-4.0)}];
      \pic[tqft/reverse pair of pants, rotate=90, at={(1.-3, -4)}];
      \draw (4-3,-3) node[xshift=0.0cm, yshift=0.0cm] {$=$};
      \pic[tqft/reverse pair of pants, rotate=90, at={(5.-3, -4)}];
      \draw (8-3,-3) node[xshift=0.0cm, yshift=0.0cm] {$=$};
      \pic[tqft/reverse pair of pants, rotate=90, at={(9.-3, -4)}];
      \pic[tqft/cylinder, rotate=90, at={(11.0-3,-3.0)}];

      \pic[tqft/pair of pants, rotate=90, at={(23, -3)}];
      \pic[tqft/cylinder, rotate=90, at={(25,-2.0)}];
      \pic[tqft/cylinder, rotate=90, at={(25,-4.0)}];
      \draw (18,-3) node[xshift=0.0cm, yshift=0.0cm] {$=$};
      \pic[tqft/pair of pants, rotate=90, at={(19, -3)}];
      \draw (22,-3) node[xshift=0.0cm, yshift=0.0cm] {$=$};
      \pic[tqft/cylinder, rotate=90, at={(13,-3.0)}];
      \pic[tqft/pair of pants, rotate=90, at={(15, -3)}];

    \end{tikzpicture}
  \end{equation}

  \item \textbf{Unit and Counit relations}
  \begin{equation}\label{eq:unit_counit}
  \begin{tikzpicture}[scale=.5, every tqft/.style={transform shape},tqft/.cd, cobordism/.style={draw},every lower boundary component/.style={draw},every incoming lower boundary component/.style={draw, dashed},every outgoing lower boundary component/.style={draw, solid}]

    \pic[tqft/cylinder, rotate=90, at={(.0,1)}];
    \pic[tqft/cap, rotate=90, at={(0,-1)}];
    \pic[tqft/reverse pair of pants, rotate=90, at={(2,-1)}];
    \draw (5,0) node[xshift=0.0cm, yshift=0.0cm] {$=$};
    \pic[tqft/cylinder, rotate=90, at={(6,.0)}];
    \draw (9,0) node[xshift=0.0cm, yshift=0.0cm] {$=$};
    \pic[tqft/cylinder, rotate=90, at={(10,-1)}];
    \pic[tqft/cap, rotate=90, at={(10,1)}];
    \pic[tqft/reverse pair of pants, rotate=90, at={(12,-1)}];

    \pic[tqft/cylinder, rotate=90, at={(29,1)}];
    \pic[tqft/cup, rotate=90, at={(29,-1)}];
    \pic[tqft/pair of pants, rotate=90, at={(27,0)}];
    \draw (22,0) node[xshift=0.0cm, yshift=0.0cm] {$=$};
    \pic[tqft/cylinder, rotate=90, at={(23,.0)}];
    \draw (26,0) node[xshift=0.0cm, yshift=0.0cm] {$=$};
    \pic[tqft/cylinder, rotate=90, at={(19,-1)}];
    \pic[tqft/cup, rotate=90, at={(19,1)}];
    \pic[tqft/pair of pants, rotate=90, at={(17,0)}];

  \end{tikzpicture}
  \end{equation}

  \item \textbf{Associativity and Coassociativity relations}
  \begin{equation}\label{eq:asso_coasso}
  \begin{tikzpicture}[scale=.5, every tqft/.style={transform shape},tqft/.cd, cobordism/.style={draw},every lower boundary component/.style={draw},every incoming lower boundary component/.style={draw, dashed},every outgoing lower boundary component/.style={draw, solid}]

    \pic[tqft/reverse pair of pants, rotate=90, at={(0,-2)}];
    \pic[tqft/reverse pair of pants, name=A, rotate=90, at={(2,-1)}];
    \pic[tqft/cylinder to prior, rotate=90, anchor=outgoing boundary 1, at=(A-incoming boundary 2)];
    \draw (5,0) node[xshift=0.0cm, yshift=0.0cm] {$=$};
    \pic[tqft/reverse pair of pants, name=B, rotate=90, at={(8,-1)}];
    \pic[tqft/reverse pair of pants, rotate=90, anchor=outgoing boundary 1, at=(B-incoming boundary 2)];
    \pic[tqft/cylinder to next, rotate=90, anchor=outgoing boundary 1, at=(B-incoming boundary 1)];

    \pic[tqft/pair of pants, rotate=90, name=C, at={(21-2,0)}];
    \pic[tqft/cylinder to next, rotate=90, anchor=incoming boundary 1, at=(C-outgoing boundary 2)];
    \pic[tqft/pair of pants, rotate=90, anchor=incoming boundary 1, at=(C-outgoing boundary 1)];
    \draw (20-2,0) node[xshift=0.0cm, yshift=0.0cm] {$=$};
    \pic[tqft/pair of pants, rotate=90, name=D, at={(15-2,0)}];
    \pic[tqft/cylinder to prior, rotate=90, anchor=incoming boundary 1, at=(D-outgoing boundary 1)];
    \pic[tqft/pair of pants, rotate=90, anchor=incoming boundary 1, at=(D-outgoing boundary 2)];

  \end{tikzpicture}
  \end{equation}

  \item \textbf{Commutativity and Cocomutativity relations}
  \begin{equation}\label{eq:commu_cocommu}
  \begin{tikzpicture}[scale=.5, every tqft/.style={transform shape},tqft/.cd, cobordism/.style={draw},every lower boundary component/.style={draw},every incoming lower boundary component/.style={draw, dashed},every outgoing lower boundary component/.style={draw, solid}]

    \pic[tqft/reverse pair of pants, name=A, rotate=90, at={(2,-1)}];
    \pic[tqft,incoming boundary components=1, outgoing boundary components=1, offset=-1, rotate=90, anchor=outgoing boundary 1, at=(A-incoming boundary 1)];
    \pic[tqft,incoming boundary components=1, outgoing boundary components=1, offset=1, rotate=90, anchor=outgoing boundary 1, at=(A-incoming boundary 2)];
    \draw (5,0) node[xshift=0.0cm, yshift=0.0cm] {$=$};
    \pic[tqft/reverse pair of pants, rotate=90, at={(6,-1)}];

    \pic[tqft/pair of pants, rotate=90, at={(11,0)}];
    \draw (14,0) node[xshift=0.0cm, yshift=0.0cm] {$=$};
    \pic[tqft/pair of pants, name=B, rotate=90, at={(15,0)}];
    \pic[tqft,incoming boundary components=1, outgoing boundary components=1, offset=-1, rotate=90, anchor=incoming boundary 1, at=(B-outgoing boundary 2)];
    \pic[tqft,incoming boundary components=1, outgoing boundary components=1, offset=1, rotate=90, anchor=incoming boundary 1, at=(B-outgoing boundary 1)];

  \end{tikzpicture}
  \end{equation}

  \item \textbf{The Frobenius relations}\index{Cobordism!2-cobordisms!Frobenius relations}
  \begin{equation}\label{eq:frobenius_relation}
  \begin{tikzpicture}[scale=.5, every tqft/.style={transform shape},tqft/.cd, cobordism/.style={draw},every lower boundary component/.style={draw},every incoming lower boundary component/.style={draw, dashed},every outgoing lower boundary component/.style={draw, solid}]

    \pic[tqft/reverse pair of pants, name=A, rotate=90, at={(2,-2)}];
    \pic[tqft/cylinder to prior, name=B, rotate=90, at={(2,2)}];
    \pic[tqft/pair of pants, anchor=outgoing boundary 1, rotate=90, at=(A-incoming boundary 2)];
    \pic[tqft/cylinder to prior, anchor=outgoing boundary 1, at=(A-incoming boundary 1), rotate=90];

    \draw (5,0) node[xshift=0.0cm, yshift=0.0cm] {$=$};

    \pic[tqft/reverse pair of pants, rotate=90, at={(6,-1)}];
    \pic[tqft/pair of pants, rotate=90, at={(8,0)}];

    \draw (11,0) node[xshift=0.0cm, yshift=0.0cm] {$=$};

    \pic[tqft/pair of pants, rotate=90, name=C, at={(12,-1)}];
    \pic[tqft/cylinder to next, name=D, rotate=90, at={(12,1)}];
    \pic[tqft/reverse pair of pants, rotate=90, anchor=incoming boundary 1, at=(C-outgoing boundary 2)];
    \pic[tqft/cylinder to next, rotate=90, anchor=incoming boundary 1, at=(C-outgoing boundary 1)];

  \end{tikzpicture}
  \end{equation}

\end{itemize}
\begin{proof}
  All of the above relations follow from \ref{lemma:cobordism_same_connected} together with \ref{prop:cobordism_diffeo}. It is evident that for every set of relations, the cobobordisms have the same number of in- and out-boundaries and genus $0$.
\end{proof}
\end{theorem}

\newpage

Before closing this chapter, a short discussion and resulting lemma on the completeness, uniqueness and role of the presented set of relations is at order. First of all, note that the choice of relations is convenient but in no way unique or minimal. In fact, there does not exist a unique minimal set of such relations. Even between only the Frobenius relation \eqref{eq:frobenius_relation}, the Unit/Counit relations \eqref{eq:unit_counit} and the Associativity/Coassociativity relations \eqref{eq:asso_coasso} there exist superfluous relations. We will show an excerpt of this in the following lemma:

\begin{lemma}
\label{lemma:ImplyAssoCoasso}
  The Frobenius relation \eqref{eq:frobenius_relation} together with the Unit/Counit relations \eqref{eq:unit_counit} imply the Associativity/Coassociativity relations \eqref{eq:asso_coasso}.
\begin{proof}
  We will show that the Frobenius relation together with the Unit/Counit relations imply the Associativity relation. The respective proof for the Coassociativity is completely analogous. To make clear, what happens in each step, changes will be highlighted in green colour:
  \begin{center}
  \begin{tikzpicture}[scale=.5, every tqft/.style={transform shape},tqft/.cd, cobordism/.style={draw},every lower boundary component/.style={draw},every incoming lower boundary component/.style={draw, dashed},every outgoing lower boundary component/.style={draw, solid}]

    \pic[tqft/reverse pair of pants, rotate=90, at={(0,-3)}];
    \pic[tqft/reverse pair of pants, name=A, rotate=90, at={(2,-2)}];
    \pic[tqft/cylinder to prior, rotate=90, anchor=outgoing boundary 1, at=(A-incoming boundary 2)];

    \draw (5,-1) node[xshift=0.0cm, yshift=0.0cm] {$=$};

    \pic[tqft/reverse pair of pants, name=B, rotate=90, at={(10,-1)}];
    \pic[tqft/reverse pair of pants, name=B1, rotate=90, anchor=outgoing boundary 1, at=(B-incoming boundary 1)];
    \pic[tqft/cylinder to prior, name=B2, rotate=90, anchor=outgoing boundary 1, at=(B-incoming boundary 2)];
    \pic[tqft/cylinder, fill=green!30, rotate=90, anchor=outgoing boundary 1, at=(B2-incoming boundary 1)];
    \pic[tqft/cylinder, fill=green!30, rotate=90, anchor=outgoing boundary 1, at=(B1-incoming boundary 2)];
    \pic[tqft/pair of pants, fill=green!30, rotate=90, name=B3, anchor=outgoing boundary 2, at=(B1-incoming boundary 1)];
    \pic[tqft/cylinder, fill=green!30, rotate=90, name=B4, anchor=incoming boundary 1, at=(B3-outgoing boundary 1)];
    \pic[tqft/cup, fill=green!30, rotate=90, anchor=incoming boundary 1, at=(B4-outgoing boundary 1)];

    \draw (13,-1) node[xshift=0.0cm, yshift=0.0cm] {$=$};

    \pic[tqft/reverse pair of pants, name=B, rotate=90, at={(18,0)}];
    \pic[tqft/cylinder, fill=green!30, rotate=90, name=B4, anchor=outgoing boundary 1, at=(B-incoming boundary 1)];
    \pic[tqft/pair of pants, fill=green!30, rotate=90, name=B3, anchor=outgoing boundary 2, at=(B4-incoming boundary 1)];
    \pic[tqft/reverse pair of pants, fill=green!30, name=B1, rotate=90, anchor=incoming boundary 2, at=(B3-outgoing boundary 1)];
    \pic[tqft/cylinder, name=B2, rotate=90, anchor=outgoing boundary 1, at=(B-incoming boundary 2)];
    \pic[tqft/cylinder, rotate=90, anchor=outgoing boundary 1, at=(B2-incoming boundary 1)];
    \pic[tqft/cup, rotate=90, anchor=incoming boundary 1, at=(B1-outgoing boundary 1)];
    \pic[tqft/cylinder, fill=green!30, rotate=90, anchor=outgoing boundary 1, at=(B1-incoming boundary 1)];

    \draw (21,-1) node[xshift=0.0cm, yshift=0.0cm] {$=$};

    \pic[tqft/reverse pair of pants, name=C, rotate=90, at={(24,0)}];
    \pic[tqft/pair of pants, rotate=90, name=C3, anchor=outgoing boundary 2, at=(C-incoming boundary 1)];
    \pic[tqft/cylinder, fill=green!30, rotate=90, name=C4, anchor=incoming boundary 1, at=(C3-outgoing boundary 1)];
    \pic[tqft/reverse pair of pants, name=C1, rotate=90, anchor=incoming boundary 2, at=(C4-outgoing boundary 1)];
    \pic[tqft/cylinder, name=C2, rotate=90, anchor=outgoing boundary 1, at=(C-incoming boundary 2)];
    \pic[tqft/cup, rotate=90, anchor=incoming boundary 1, at=(C1-outgoing boundary 1)];
    \pic[tqft/cylinder, name=C5, rotate=90, anchor=outgoing boundary 1, at=(C1-incoming boundary 1)];
    \pic[tqft/cylinder, fill=green!30, rotate=90, anchor=outgoing boundary 1, at=(C5-incoming boundary 1)];
    \pic[tqft/cylinder, fill=green!30, rotate=90, anchor=incoming boundary 1, at=(C-outgoing boundary 1)];

  \end{tikzpicture}
  \end{center}
  So far we have used the Counit relation in the first step to create a situation in which we can use the Frobenius relation in the second step. The third step just "inserted" and deleted some identity cobordisms to again prepare the stage for the Frobenius relation, this time in the upper branch of the cobordism. From here on, we note:
  \begin{center}
  \begin{tikzpicture}[scale=.5, every tqft/.style={transform shape},tqft/.cd, cobordism/.style={draw},every lower boundary component/.style={draw},every incoming lower boundary component/.style={draw, dashed},every outgoing lower boundary component/.style={draw, solid}]

    \pic[tqft/reverse pair of pants, name=C, rotate=90, at={(2,0)}];
    \pic[tqft/pair of pants, rotate=90, name=C3, anchor=outgoing boundary 2, at=(C-incoming boundary 1)];
    \pic[tqft/cylinder, rotate=90, name=C4, anchor=incoming boundary 1, at=(C3-outgoing boundary 1)];
    \pic[tqft/reverse pair of pants, name=C1, rotate=90, anchor=incoming boundary 2, at=(C4-outgoing boundary 1)];
    \pic[tqft/cylinder, name=C2, rotate=90, anchor=outgoing boundary 1, at=(C-incoming boundary 2)];
    \pic[tqft/cup, rotate=90, anchor=incoming boundary 1, at=(C1-outgoing boundary 1)];
    \pic[tqft/cylinder, name=C5, rotate=90, anchor=outgoing boundary 1, at=(C1-incoming boundary 1)];
    \pic[tqft/cylinder, rotate=90, anchor=outgoing boundary 1, at=(C5-incoming boundary 1)];
    \pic[tqft/cylinder, rotate=90, anchor=incoming boundary 1, at=(C-outgoing boundary 1)];

    \draw (7,-1) node[xshift=0.0cm, yshift=0.0cm] {$=$};

    \pic[tqft/reverse pair of pants, fill=green!30, name=D, rotate=90, at={(8,-1)}];
    \pic[tqft/pair of pants, fill=green!30, rotate=90, name=D3, anchor=incoming boundary 1, at=(D-outgoing boundary 1)];
    \pic[tqft/reverse pair of pants, name=D1, rotate=90, anchor=incoming boundary 2, at=(D3-outgoing boundary 1)];
    \pic[tqft/cylinder, rotate=90, name=D4, anchor=incoming boundary 1, at=(D3-outgoing boundary 2)];
    \pic[tqft/cup, rotate=90, anchor=incoming boundary 1, at=(D1-outgoing boundary 1)];
    \pic[tqft/cylinder, name=D5, rotate=90, anchor=outgoing boundary 1, at=(D1-incoming boundary 1)];
    \pic[tqft/cylinder, rotate=90, anchor=outgoing boundary 1, at=(D5-incoming boundary 1)];

    \draw (15,-1) node[xshift=0.0cm, yshift=0.0cm] {$=$};

    \pic[tqft/reverse pair of pants, name=D, rotate=90, at={(16,-1)}];
    \pic[tqft/reverse pair of pants, name=D1, fill=green!30, rotate=90, anchor=incoming boundary 2, at=(D-outgoing boundary 1)];
    \pic[tqft/pair of pants, fill=green!30, rotate=90, name=D3, anchor=incoming boundary 1, at=(D1-outgoing boundary 1)];
    \pic[tqft/cylinder, rotate=90, name=D4, anchor=incoming boundary 1, at=(D3-outgoing boundary 2)];
    \pic[tqft/cup, rotate=90, anchor=incoming boundary 1, at=(D3-outgoing boundary 1)];
    \pic[tqft/cylinder to next, name=D5, rotate=90, anchor=outgoing boundary 1, at=(D1-incoming boundary 1)];

    \draw (25,-1) node[xshift=0.0cm, yshift=0.0cm] {$=$};

    \pic[tqft/reverse pair of pants, name=A, rotate=90, at={(26,-1)}];
    \pic[tqft/reverse pair of pants, name=A1, fill=green!30, rotate=90, anchor=incoming boundary 2, at=(A-outgoing boundary 1)];
    \pic[tqft/cylinder to next, rotate=90, anchor=outgoing boundary 1, at=(A1-incoming boundary 1)];

  \end{tikzpicture}
  \end{center}
  This proves the statement. Note that we could have introduced parts of this proof as lemmas but since no further application is intended at this point, they have been omitted resulting in a slightly longer proof.
\end{proof}
\end{lemma}

The above proof clearly shows, that we based our choice on convenience rather than strict criteria. And indeed, every of the above relations reflects a certain aspect of previous and/or upcoming notions. First note that the \emph{Identity relations} \eqref{eq:identity_relations} basically embody the idea of the cylinder cobordism serving as an identity, hence the name. \emph{The Unit/Counit relations} \eqref{eq:unit_counit} show that the \textbf{cup cobordism} $\emptyset \Lra 1$ and the \textbf{cap cobordism} $1 \Lra \emptyset$ act as unit and counit respectively. This corresponds to the empty manifold being the unit/counit of the tensor product. Upon closer inspection of the \emph{Associativity and Coassociativity relations} \eqref{eq:asso_coasso} we recognize the associativity of the tensor product in a monoidal category. Meanwhile the \emph{Commutativity and Cocomutativity relations} represent the invariance of the tensor product and the coproduct under the braiding isomorphism, namely the twist cobordism. The \emph{Frobenius relations} \eqref{eq:frobenius_relation} will become apparent in \emph{chapter} \ref{subsec:2DTQFT} when we explore the connection between certain Frobenius algebras $\comk$ and $\cobtwo$.\\

This concludes the chapter on $2$-dimensional cobordisms. We will again encounter $2$-cobordisms in \emph{chapter} \ref{subsec:CatComfrob} where we present a first pictorial link between the category of commutative Frobenius algebras over a field $\kk$, $\comk$, and $\cobtwo$ which will greatly help when treating $2$-dimensional TQFTs in \emph{subchapter} \ref{subsec:2DTQFT}.

\newpage

\newpage
\section{The Category of Commutative Frobenius Algebras}
\label{sec:Comfrob}

At first, this chapter might seem a little off. After all we are working towards a description of \emph{Topological Quantum Field Theory}. So far we have mainly been concerned with general category theory, particularly symmetric monoidal categories and smooth manifolds qualifying as such a category, namely $\cobn{n}$. At least the latter ties directly to a geometric formulation of physics and while the precise definition will have to wait until chapter \ref{sec:TQFT}, we can "read" the objects of $\cobn{n}$ as physical states and the cobordisms between them as evolutions embodying various types of interactions. On the other hand introducing a special type of algebra does not have such an adhoc analogy in physics.\\

It will soon become apparent that Frobenius Algebras over a field $\kk$ do not only classify as a symmetric monoidal category, they also exhibit some striking resemblances to $2$-dimensional cobordisms, that is $\cobtwo$. This resemblance will in a later subchapter be formulated into a full theorem (see \ref{theorem:Equivalence2TQFT}) unveiling a deep link between $2$-dimensional TQFTs and the category $\comk$. It is due to this theorem, which allows for a full classification of $2$-dimensional TQFTs, that an exhibition of Frobenius Algebras takes place in this chapter. A first pictorial representation of the analogy between the two categories will be shown in \ref{subsec:CatComfrob}.\\

Note that due to the rather special role of Frobenius Algebras in this project, the introduction will remain rather focused on the main results needed for later applications. While this won't preclude some interesting proofs and examples, the interested reader is referred to \cite[Chapter 2]{FrobAlgebraTQFT} where a far more extensive treatment using linear algebra and multiple classifications of Frobenius Algebras are discussed in detail. Meanwhile this project will take a very categorial point of view. Needless to say we will also incorporate ideas from \cite{Intro_TQFT} especially when looking towards applications in $TQFTs$.

\subsection{Algebraic Preliminaries}
\label{subsec:AlgebraicPreliminaries}

To introduce Frobenius Algebras, we need some background on vector spaces and algebras in general. Frobenius Algebras are algebras over a field equipped with extra structure in the form of so-called \emph{nondegenerate pairings}, maps that preserve certain algebraic properties of the involved vector spaces, both of which will be defined in this subsection. Note that the material does also include a few basic concepts taught in introductory courses on linear algebra to stay as self contained as possible and fix important notation.\\

Note that when talking about vector spaces, we mean vector spaces over a field $\kk$, thus objects of the symmetric monoidal category $\vectk$. When working with multiple vector spaces, they share the same field unless stated otherwise. We start by defining arbitrary pairings on vector spaces.

\begin{definition}[Pairings]
  Let $V,W$ be two vector spaces. A \textbf{pairing}\index{Pairing} of $V$ and $W$ is a linear map $\beta \colon V \tensor W \lra \kk$. We denote its action on elements of the two spaces by
  \begin{align*}
    \beta \colon V \tensor W &\lra \kk, \\
    v \tensor w & \lra \pair{v}{w}
  \end{align*}
\end{definition}

What we need next is a notion of non-degeneracy of such pairings. That is a classification of such pairings that are in some sense invertible.

\begin{definition}[Non-degenerate pairings]
\label{def:NdegPairing}
  A pairing $\beta \colon V \tensor W \lra \kk$ for two vector spaces $V,W$ is called \textbf{nondegenerate in $V$}\index{Pairing!nondegenerate}, if there exists a map $\gamma_V \colon \kk \lra W \tensor V$ such that the following diagram commutes
  \begin{center}
  \begin{tikzcd}[sep = huge]
      V \tensor (W \tensor V) \arrow[r, "\simeq"] & (V \tensor W) \tensor V \arrow[d, "\beta \tensor \id_V"] \\
      V = V \tensor \kk \arrow[u, "\id_V \tensor \gamma_V"] \arrow[r, leftrightarrow, "\id_V"'] & V = \kk \tensor V
  \end{tikzcd}
  \end{center}
  $\gamma$ is called a \textbf{copairing in $V$}\index{Pairing!copairing}. We call $\beta$ \textbf{nondegenerate in $W$} if there exists a $\gamma_W$ that the respective diagram for $W$ commutes. Accordingly, this $\gamma_W$ would be called a \textbf{copairing in $W$}. $\beta$ is simply called \textbf{nondegenerate} if it is nondegenerate in both variables.
\end{definition}

\begin{remark}
  Note that the above definition of a nondegenerate pairing does indeed coincide with either the usual conditions
  \begin{enumerate}
    \item $\pair{v}{w} = 0 \ \forall v \in V \ \quad \Lra w = 0$
    \item $\pair{v}{w} = 0 \ \forall w \in W \quad \Lra v = 0$
  \end{enumerate}
  once we know that the two spaces are of the same \emph{finite} dimension. We can then recover injectivity and formulate nondegeneracy in the usual way. The interested reader is referred to \cite[Lemma 2.1.15]{FrobAlgebraTQFT} for a deeper analysis of this fact. In this project, we will only make use of the above, admitedly unusual, form of nondegeneracy since its abstract form will prove quite helpful in some of the upcoming proofs.
\end{remark}

A first interesting result is that the copairing of a nondegenerate pairing is in fact unique meaning that we can unambigously talk about "the copairing" of a nondegenerate pairing.

\begin{lemma}
  If a pairing $\beta \colon V \tensor W \lra \kk$ is nondegenerate, the two copairings in either variable agree, namely $\gamma_V = \gamma_W$. In fact any pair of such maps need to agree and thus the copairing is unique.
\begin{proof}
  First note that due to the degeneracy of $\beta$, we can factor the identities on $V$ and $W$ as
  $$ \id_V = (\beta \tensor \id_V) \circ (\id_V \tensor \gamma_V) \quad \text{and} \quad \id_W = (\id_W \tensor \beta) \circ (\gamma_W \tensor \id_W) $$
  This lets us express the two copairings as
  \begin{align*}
    \gamma_V = (\id_W \tensor \id_V) \circ \gamma_V &= (\id_W \tensor \beta \tensor \id_V) \circ (\gamma_W \tensor \id_W \tensor \id_V)  \circ \gamma_V\\
    &= (\id_W \tensor \beta \tensor \id_V) \circ (\gamma_W \tensor \gamma_V) \\
    &= (\id_W \tensor \beta \tensor \id_V) \circ (\id_W \tensor \id_V \tensor \gamma_V)  \circ \gamma_W = (\id_W \tensor \id_V) \circ \gamma_W = \gamma_W
  \end{align*}
  Which proves that $\gamma_V = \gamma_W$. Since we did not make any asumption about the choice of possible copairings, this is a general result for any possible pair of copairings in different variable. This proves that the copairing is in fact unique.
\end{proof}
\end{lemma}

Now given any pairing $\beta \colon V \tensor W \lra \kk$ we can define two naturally arising maps onto the dual vector spaces by setting
\begin{align*}
  \beta_l \colon W &\lra V^*, \quad w \lmap \pair{\ \cdot \ }{w} \\
  \beta_r \colon V &\lra W^*, \quad v \lmap \pair{v}{\ \cdot \ }
\end{align*}
Suggestively named, they stand for the left and right insertion respectively. An interesting result about these maps and the underlying vector spaces is the following:

\begin{lemma}
\label{lemma:NdegPairingFinite}
  A pairing $\beta \colon V \tensor W \lra \kk$ of vector spaces is nondegenerate in $V$ iff $V$ is finite dimensional and the map $\beta_l$ is injective. Accordingly it is nondegenerate in $W$ iff $W$ is finite dimensional and the map $\beta_r$ is injective.
\begin{proof}
  See \cite[Lemma 2.1.12]{FrobAlgebraTQFT}.
\end{proof}
\end{lemma}

Using these definitions, we find that for finite-dimensional vector spaces the two maps are closely related via the dual functor:

\begin{lemma}
\label{lemma:PairingFiniteDual}
  For a pairing $\beta \colon V \tensor W \lra \kk$ of two vector spaces of finite dimension, the two maps $\beta_l$ and $\beta_r$ are the dual of each other.
\begin{proof}
  We inspect the dual of $\beta_l$, the same arguments can be drawn for $\beta_r$. The dual is given by
  \begin{align*}
    \Hom(V^*, \kk) & \lra \Hom(W, \kk), \quad \quad T \lmap T \circ \beta_l
  \end{align*}
  Since we can identify $V \simeq \Hom(V^*, \kk)$ by mapping $v \in V$ to $T_v := [\lambda \lmap \lambda(v)]$, we can inspect the dual of $\beta_l$ as the map
  \begin{align*}
    V & \lra \Hom(W, \kk), \quad \quad v \lmap T_v \circ \beta_l
  \end{align*}
  But now the action of $T_v \circ \beta_l$ on an element $w \in W$ is given by
  $$ T_v \circ \beta_l(w) = T_v (\pair{\ \cdot \ }{w}) = \pair{v}{w} = \beta_r(v) (w)$$
  But this just shows that the dual map is indeed the same as $\beta_r$.
\end{proof}
\end{lemma}

Combining these two results, we first see that nondegeneracy of a pairing means that both $\beta_l$ and $\beta_r$ are injective and both underlying vector spaces are finite dimensional. In this case the existence of these two isomorphisms uniquely determines $\beta$ by setting $\beta := \beta_r(\ins) (\ins) = \beta_r(\ins)(\ins)$. Using the duality between the two maps, we arrive at the following statement:

\begin{lemma}
\label{lemma:NdegPairingIso}
  Given a pairing $\beta \colon V \tensor W \lra \kk$ between finite-dimensional vector spaces, the following statements are equivalent:
  \begin{enumerate}[label=(\roman*), itemsep = -.5ex]
    \item $\beta$ is nondegenerate.
    \item $\beta_l$ is an isomorphism between $W$ and $V^*$.
    \item $\beta_r$ is an isomorphism between $V$ and $W^*$.
  \end{enumerate}
  Note that this also shows that $V$ and $W$ are of the same dimension, once there exists a nondegenerate pairing between them.
\begin{proof}
  To show that $(ii) \Leftrightarrow (iii)$ we simply note that the dual functor is an equivalence of categories and thus in particular preserves invertibility. With this established, it is also immediate that each of the two statements result in $(i)$ being true. The last piece is to start with $(i)$: If $\beta$ is nondegenerate, both $\beta_l$ and $\beta_r$ are injective. Thus their duals are surjective which amounts to saying that both are bijective and hence isomorphisms of their respective domain and codomain.
\end{proof}
\end{lemma}

With a well-defined notion of nondegenerate pairings at hand, we define the algebras which we will later endow with such pairings. Note that while we first give a definition of $\kk$-algebras in $\vectk$, the concept is directly embedded into the categorial preliminaries.

\begin{definition}[$\kk$-Algebras – Vector Space Version]
\label{def:kAlgebra}
  A $\kk$-vector space $A$ together with two $\kk$-linear maps
  \begin{align*}
    \mu \colon A \tensor A \lra A, \quad \quad \eta \colon \kk \lra A
  \end{align*}
  such that the following diagrams commute
  \begin{center}
  \begin{tikzcd}
    A \otimes (A \otimes A) \arrow[r, "\sim"] \arrow[d, "\id_A \otimes \mu"'] & (A \otimes A) \otimes A \arrow[r, "\mu \otimes \id_A"] & A\otimes A \arrow[d, "\mu"] & \kk \otimes A \arrow[dr, "\sim"'] \arrow[r, "\eta \otimes \id_A"]& A \otimes A \arrow[d, "\mu"] & A \otimes \kk \arrow[l, "\id_A \otimes \eta"'] \arrow[dl, "\sim"]\\
    A \otimes A \arrow[rr, "\mu"] & & A & & A &
  \end{tikzcd}
  \end{center}
  is called a \textbf{$\kk$-algebra}. We call $\mu$ the \textbf{multiplication map}\index{$\kk$-Algebra!multiplication} and $\eta$ the \textbf{unit map}\index{$\kk$-Algebra!unit}. We denote the action of the multiplication map by
  $$ \mu \colon A \tensor A \lra A, \quad \quad x \tensor y \lmap xy $$
\end{definition}

As promised, this will not be the definitive definition of a $\kk$-algebra. Looking at the commuting diagrams in its above definition, we see a striking resemblance to the ones in the definition of a monoid in a monoidal category \ref{def:monoidal_monoids}. And in fact a $\kk$-algebra will turn out to be just that!\\

However to draw this beautiful parallel, we need to have an underlying monoidal category. Of course the natural candidate is $\vectk$ since every $\kk$-algebra particularly is a $\kk$-vector space. We already know that $\vectk$ is a monoidal category regarding it as the triple $\triple{\vectk}{\tensor}{\kk}$, whose properties we have quietly exploited throughout this introduction. Without further ado we can now formulate a beautiful compact categorial definition of $\kk$-algebras:

\begin{definition}[$\kk$-Algebras – Categorial Version]
\label{def:kAlgebraCoalgebraCat}
  A $\kk$-algebra is a monoid in $\vectk$. Without giving the conrete definition as in \ref{def:kAlgebra}\index{$\kk$-Algebra} note that \emph{comonoids} in $\vectk$ correspond to $\kk$-Coalgebras\index{$\kk$-Algebra!$\kk$-Coalgebra}.
\end{definition}

Since monoids of a category forms a category itself, we denote the category $\cat{Mon}_{\vectk}$ by $\alg$ where the arrows are $\kk$-algebra morphisms defined as in \ref{def:monoidal_monoid_arrows}. From here on out, we can think of $\kk$-algebras as objects of $\alg$. Following this reading we will subsequently deviate from the usual formulation of and leading up to Frobenius Algebras by choosing a more categorial vocabulary. In this spirit it is the natural next step to look at the left and right actions \ref{def:monoid_action} of $\alg$ on $\vectk$:

\begin{definition}[Left and Right Modules]
  Given $A \in \alg$ and $V \in \vectk$, we call a right action of $A$ on $V$ together with $V$ a \textbf{right A-module}\index{$\kk$-Algebra!action}\index{$\kk$-Algebra!module}, for left actions we use the term \textbf{left A-module}. In a slight abuse of notation we also call $V$ itself a left or right A-module implicitely assuming a fitting left or right action.
\end{definition}

Speaking in terms of vector spaces a right $A$-module is a vector space $V$ together with a $\kk$-linear map
$$ \rho \colon V \tensor A \lra V, \quad \quad x \tensor a \lmap \rho(x,a) $$
such that the according diagrams for an action of a monoid commute which can be expressed for arbitrary $x\in V$ and $a,b \in A$ as
$$ \rho(\rho(x,a),b) = \rho(x,ab) \quad \quad \text{\&} \quad \quad \rho(x,e) = x $$
Accordingly a left $A$-module for $V$ is the tuple $(V, \lambda)$ where $\lambda$ is $\kk$-linear and maps
$$ \lambda \colon A \tensor V \lra V, \quad \quad a \tensor x \lmap \lambda(a,x) $$
such that
$$ \lambda(a,\lambda(b,x)) = \lambda(ab,x) $$

Note that following this definition, every element of $\alg$ in particular forms a right module with itself, i.e. every $A \in \alg$ is a right $A$-module with the right action given by the multiplication map. Closely following the categorial treatment of actions of monoids, we aim to formulate morphisms of left and right actions \ref{def:monoid_action_morphisms} in the context of algebra modules.

\begin{definition}[$A$-Module Homomorphisms]
  Given $A \in \alg$ and two left or two right $A$-modules $V,W$ we call a morphism of their respective left or right actions a \textbf{left} or \textbf{right A-homomorphism}\index{$\kk$-Algebra!module morphism}. To further adapt common notation, we call such an homomorphism \textbf{left} or \textbf{right A-linear}\index{$\kk$-Algebra!module morphism!left/right linear}, a name stemming from the suggestive commuting diagrams in \ref{def:monoid_action_morphisms}.
\end{definition}

Note that for $\kk$-algebras, we denote the category of left and right actions, thus of left and right modules, now given adapted morphisms in the form of left and right $A$-module homomorphisms, by $\lmod{A}$ and $\rmod{A}$. Note that the two are just more suggestive names for $\lact{A}$ and $\ract{A}$ in the monoidal category $\vectk$. Using these terms we can formulate a very interesting application of the concept of duality:

\begin{lemma}[Duals of Module Homomorphisms]
  Given $A \in \alg$ the natural dual operation on left or right modules and their homomorphisms is a contravariant functor $\lmod{A} \lra \rmod{A}$ or $\rmod{A} \lra \lmod{A}$ respectively.
\begin{proof}
  Again we will prove the result only for left actions, the right variant is completely analogous. Let $V, W \in \vectk$ be two left $A$-modules and $\phi \colon V \lra W$ be a morphism of their left actions $\sigma, \lambda$. Our first goal is to construct a natural right $A$-module structure on $V^*$. To this end, note that $V^* = \Hom(V, \kk)$ and thus we can define
  $$ \sigma^* \colon V^* \tensor A \lra V^*, \quad \quad \chi \tensor a \lmap \sigma^*(\chi, a) := [x \in V \lmap \chi(\sigma(a, x))] $$
  Now since for arbitrary $\chi \in V^*$, $x \in V$ and $a,b \in A$
  $$ \sigma^*(\sigma^*(\chi, a), b)(x) = \sigma^*(\chi, a)(\sigma(b,x)) = \chi(\sigma(a, \sigma(b,x))) = \chi(\sigma(ab,x))) = \sigma^*(\chi, ab)(x)$$
  the constructed map is indeed a right action making $(V^*, \sigma^*)$ a right $A$-module which shows that the dualisation works on objects. Now to combat the left $A$-module homomorphism $\phi$, we define its dual morphism as
  $$ \phi^* \colon W^* \lra V^*, \quad \quad \theta \lmap \theta(\phi(\ins)) $$
  Now all we need to check is if the following diagram commutes to verify that $\phi^*$ is a right $A$-module homomorphism:
  \begin{center}
  \begin{tikzcd}[sep=huge]
    W^* \tensor A \arrow[d, "\lambda^*"'] \arrow[r, "\phi^* \tensor \id_A"]  & V^* \tensor A \arrow[d, "\sigma^*"] \\
    W^* \arrow[r, "\phi^*"'] & V^*
  \end{tikzcd}
  \end{center}
  However we can simply use our definitions together with the properties of module homomorphisms to state for arbitrary $a \in A$, $\theta \in W^*$ and $x \in V$:
  \begin{align*}
    \sigma^*( \phi^* \tensor \id_A (\theta \tensor a)) &= \sigma^*( \theta(\phi(\ins)), a) = \theta(\phi(\sigma(a, \ins))) \\
    &= \theta(\lambda(a,\phi(\ins))) = \phi^*(\theta(\lambda(a, \ins))) = \phi^*(\lambda^*(\theta, a))
  \end{align*}
  This shows that the two defined dual mappings on modules and their homomorphisms do indeed act like a functor between $\lmod{A}$ and $\rmod{A}$. The reverse statement follows analogously.
\end{proof}
\end{lemma}

Now having a dual map at hand, it is only natural to try to recover the idea of a double dual of an object being isomorphic to the object itself. And indeed for finite-dimensional vector spaces, we find the analogous statement:

\begin{lemma}
  If for any $A \in \alg$ we restrict to modules in the full subcategory of finite-dimensional vector spaces, $\cat{FinVect}_\kk$, the above dual functor is an equivalence of categories between the full subcategories $\cat{Finlmod}_A$ of $\lmod{A}$ and $\cat{Finrmod}_A$ of $\rmod{A}$.
\begin{proof}
  Let $V,W \in \cat{FinVect}_\kk$ be left $A$-modules with actions $\sigma, \lambda$ and $\phi \colon V \lra W$ a left $A$-module homomorphism of those actions. We already know that $\psi^* \colon W^* \lra V^*$ is a fitting right $A$-module homomorphism between the right actions $\sigma^*$ and $\lambda^*$. Now we can exploit our previous result about the dual of objects on $\finvect$ (see \ref{lemma:double_dual}) to state $V^{**} \simeq V$ and $W^{**} \simeq W$. Thus we investigate
  $$ \phi^{**} := (\phi^*)^* \colon V \lra W, \quad \quad x \lmap \phi^*(\ins)(x) $$
  But for any $\theta \in W^*$ we have
  $$ \theta(\phi^{**}(x)) = \phi^*(\theta)(x) = \theta(\phi(x)) $$
  and thus $\phi^{**} = \phi$. Now investigating the action $\sigma$ we need to show that $\sigma^{**} = \sigma$ to conclude the proof. To this end we explicitely write down its definition:
  $$ \sigma^{**} := (\sigma^*)^* \colon V \tensor A \lra V, \quad \quad x \tensor a  \lmap \sigma^{**}(a, x) := [ \chi \in V^* \lmap \sigma^*(\chi, a)(x)] $$
  Thus for any $a \in A$, $x \in V$ and $\chi \in V^*$ we get
  $$ \chi(\sigma^{**}(a, x)) =  \sigma^*(\chi, a)(x) = \chi(\sigma(a,x)) $$
  which proves that the defined dual operation is indeed an equivalence of categories between $\cat{Finlmod}_A$ and $\cat{Finrmod}_A$.
\end{proof}
\end{lemma}

Turning back to pairings of vector spaces of the form $\beta \colon V \tensor W \lra \kk$, we quickly see that, regardless of $V,W$ being left or right $A$-modules, the field $\kk$ is not an $A$-module and thus $\beta$ can not be $A$-linear. However for the case that $V$ is a right $A$-module and $W$ is a left $A$-module, we can investigate a particular diagram:

\begin{definition}[Associative Pairings]
\label{def:AssociativePairings}
  A pairing $\beta \colon V \tensor W \lra \kk$ where $V$ is a right $A$-module with action $\rho$ and $W$ is a left $A$-module with action $\lambda$ is called an \textbf{associative pairing}\index{Pairing!associative} if the following diagram commutes:
  \begin{center}
  \begin{tikzcd}[sep=huge]
    V \tensor A \tensor W \arrow[d, "\id_V \tensor \lambda"'] \arrow[r, "\rho \tensor \id_W"] & V \tensor W \arrow[d, "\beta"] \\
    V \tensor W \arrow[r, "\beta"'] & \kk
  \end{tikzcd}
  \end{center}
  The reason for the "associative" prefix is more apparent when writing the above diagram for arbitrary $x,y \in V$ and $a \in A$ as
  $$ \beta(x \tensor \lambda(a,y)) = \beta(\rho(x,a) \tensor y) $$
\end{definition}

Using this definition leads to an interesting consequence for the mappings $\beta_l$ and $\beta_r$:

\begin{lemma}
\label{lemma:PairingAssociativeLinear}
  For a pairing $\beta \colon V \tensor W \lra \kk$ where $V$ is a right $A$-module with action $\rho$ and $W$ is a left $A$-module with action $\lambda$ the following three statements are equivalent:
  \begin{enumerate}[label=(\roman*), itemsep = -.5ex]
    \item $\beta$ is associative
    \item $\beta_l \colon W \lra V^*$ is left $A$-linear
    \item $\beta_r \colon V \lra W^*$ is right $A$-linear
  \end{enumerate}
\begin{proof}
  Left linearity of $\beta_l$ can be expressed diagrammatically using the properties of morphisms of actions of monoids as
  \begin{center}
  \begin{tikzcd}[sep=huge]
    A \tensor W \arrow[d, "\lambda"'] \arrow[r, "\id_A \tensor \beta_l"]  & A \tensor V^* \arrow[d, "\rho^*"] \\
    W \arrow[r, "\beta_l"'] & V^*
  \end{tikzcd}
  \end{center}
  Now this diagram is equivalent to the following equation being true for any $x \in V$, $t \in W$ and $a \in A$:
  \begin{align*}
    \pair{x}{\lambda(a,t)} = \pair{\ins}{\lambda(a,t)}(x) = \beta_l(\lambda(a,t))(x) = \rho^*(\id_A \tensor \beta_l(a \tensor t)) (x) = \rho^*(a, \pair{\ins}{t}) (x) = \pair{\rho(x,a)}{t}
  \end{align*}
  which amounts to $\pair{x}{\lambda(a,t)} = \pair{\rho(x,a)}{t}$ i.e. the associativity condition. The same reasoning holds for $\beta_r$ which concludes the proof.
\end{proof}
\end{lemma}

This concludes the short introduction to pairings, $\kk$-algebras and modules. In the next chapter we will use many of the above concepts and results for a streamlined yet thorough introduction to Frobenius Algebras.

\newpage
\subsection{Frobenius Algebras}
\label{subsec:Frob}

This subchapter is dedicated to introducing the central notions and definitions of Frobenius Algebras and their commutative versions. Note that the ansatz chosen here is, akin to the previous chapter, rather categorial following the main definition from \cite{FrobeniusPaper} to give a top-down view on Frobenius Algebras. Of course we will constrast the categorial introduction with a concrete discussion of Frobenius Algebras in $\vectk$ following the lines of \cite{FrobAlgebraTQFT} and proving our categorial definition to be equivalent to a more standard approach. This will indeed be the main result of this subchapter whose goal it is to introduce general constructs in category theory and have them descend onto the monoidal category $\vectk$ to enable both an abstract and a concrete treatment of $2D$-$TQFT$ in \ref{subsec:2DTQFT}.\\

Without further ado, let us define Frobenius Algebras in a monoidal category. Note that it draws upon the concepts of monoids \ref{def:monoidal_comonoids} and comonoids \ref{def:monoidal_comonoids}:

\begin{definition}[Frobenius Algebra – Category Version]
\label{def:FrobeniusCat}
  Let $\triple{\cat{C}}{\otimes}{e}$ be a monoidal category and $A \in \cat{C}$. A tuple $\pair{A, \mu, \eta, \delta}{\epsilon}$ is called a \textbf{Frobenius Algebra}\index{Frobenius Algebra} if
  \begin{enumerate}[itemsep = 0ex]
    \item $\triple{A}{\mu}{\eta}$ is a monoid,

    \item $\triple{A}{\delta}{\epsilon}$ is a comonoid,

    \item the following diagram, called the \textbf{Frobenius Relation}\index{Frobenius Algebra!Frobenius Relation}, commutes:
    \begin{center}
    \begin{tikzcd}[sep = huge]\label{eq:FrobeniusRelationFrob}
      (A \tensor A) \tensor A \arrow[dd, "\simeq"'] & \arrow[r, "\id_A \tensor \delta"] A \arrow[l, "\delta \tensor \id_A"'] \tensor A \arrow[d, "\mu"] & A \tensor (A \tensor A) \arrow[dd, "\simeq"] \\
      & A \arrow[d, "\delta"] & \\
      A \tensor (A \tensor A) \arrow[r, "\id_A \tensor \mu"'] & A \tensor A & \arrow[l, "\mu \tensor \id_A"] (A \tensor A) \tensor A
    \end{tikzcd}
    \end{center}
  \end{enumerate}
\end{definition}

\begin{example}
  An immediate example is the circle $1 \in \cobtwo$. Looking at the set of generators \ref{theorem:cob2_classification} and the set of relations they satisfy \ref{theorem:cobordism_relations} we can immediately recover all of the above properties. Note that the additional relations involving the twist cobordism even make $1$ a \emph{commutative Frobenius Algebra} which will be defined later on.
\end{example}

Looking at \ref{lemma:RigidComonoid}
we see that the underlying object of any Frobenius Algebra is automatically rigid. Thus we can already assert that any Frobenius Algebra over the monoidal category $\vectk$ is in particular defined over a finite-dimensional vector space! Another immediate result is that any monoidal functor preserves Frobenius Algebras in a monoidal category, that is it maps Frobenius Algebras to Frobenius Algebras. The next result shows that the counit of a Frobenius Algebra induces a "pairing" morphism:

\begin{lemma}
\label{lemma:PairingFrob}
  Let $\pair{A, \mu, \eta, \delta}{\epsilon}$ be a Frobenius Algebra in a monoidal category $\triple{\cat{C}}{\otimes}{e}$. Then there exist
  morphisms $\beta \colon A \tensor A \lra e$ and $\theta \colon e \lra A \tensor A$ such that the following diagrams commute:
  \begin{center}
  \begin{tikzcd}[sep=huge]
    A \tensor A \tensor A \arrow[d, "\id_A \tensor \mu"'] \arrow[r, "\mu \tensor \id_A"] & A \tensor A \arrow[d, "\beta"] & A \tensor A \tensor A & \arrow[l, "\id \tensor \delta"'] A \tensor A \\
    A \tensor A \arrow[r, "\beta"'] & e & A \tensor A \arrow[u, "\delta \tensor \id"] & e\arrow[l, "\theta"] \arrow[u, "\theta"']
  \end{tikzcd}
  \end{center}
  To draw a parallel to pairings in the context of $\kk$-Algebras, we call such morphisms \textbf{pairings}\index{Frobenius Algebra!pairing} and \textbf{copairings}\index{Frobenius Algebra!copairing} respectively. Note that the diagrams are categorial mirrors of \ref{def:AssociativePairings}.
\begin{proof}
  To prove this, simply define $\beta := \epsilon \circ \mu$ and $\theta:= \delta \circ \eta$. Now note
  \begin{align*}
    \beta \circ ( \id_A \tensor \mu) &= \epsilon \circ \mu \circ (\id_A \tensor \mu) = \epsilon \circ \mu \circ (\mu \tensor \id_A) = \beta \circ (\mu \tensor \id_A) \\
    (\delta \tensor \id_A) \circ \theta &= (\delta \tensor \id_A) \circ (\delta \circ \eta) = (\id_A \tensor \delta) \circ (\delta \circ \eta) = (\id_A \tensor \delta) \circ \theta
  \end{align*}
  Thus the above diagrams do indeed commute.
\end{proof}
\end{lemma}

\begin{remark}
  Note that our choices for $\beta$ and $\theta$ were in no way unique since we did not use any commutation properties of $\epsilon$ or $\eta$ respectively. Any morphisms $A \lra e$ and $e \lra A$ could take their place. We are however obliged to use $\mu$ and $\delta$ respectively. This can in fact be incorporated into various equivalent definitions of Frobenius Algebras of $\finvect$ which we won't treat in this project. The interested reader is referred to \cite[Lemma 2.2.8]{FrobAlgebraTQFT}.
\end{remark}

We can however classify certain types of such morphisms which will allow for an interesting parallel between them and nondegenerate pairings \ref{def:NdegPairing} which we discussed previously:

\begin{definition}
\label{def:NdegPairingFrob}
  We call a morphism $\beta \colon A \tensor A \lra e$ with the same commutation properties as in \ref{lemma:PairingFrob} \textbf{nondegenerate}\index{Frobenius Algebra!nondegenerate pairing/copairing} iff there exists a copairing $\theta \colon e \lra A \tensor A$ such that the following diagram commutes:
  \begin{center}
  \begin{tikzcd}[sep=huge]
    A \tensor A \tensor A \arrow[r, "\beta \tensor \id_A"] & e \tensor A \simeq A \\
    A \simeq A \tensor \arrow[u, "\id_A \tensor \theta"] \arrow[ur, "\id_A"'] e & ~
  \end{tikzcd}
  \end{center}
  Subsequently we also call the copairing \textbf{nondegenerate}.
\end{definition}

Indeed this leads us to the following interesting result about nondegenerate pairings:

\begin{lemma}
\label{lemma:CopairingUnique}
  Let $\beta$ be a nondegenerate pairing and $\theta$ the respective copairing such that the diagrams in \ref{def:NdegPairingFrob} commute. Any other copairing $\xi$ to which this applies is isomorphic to $\theta$. Moreover $\beta := \epsilon \circ \mu$ is a nondegenerate pairing with copairing $\theta := \delta \circ \eta$.
\begin{proof}
  Let $\beta, \theta $ and $\xi$ be as above. Then in particular
  $$ (\id_A \tensor \id_A) \circ \theta \simeq (\id_A \tensor \beta \tensor \id_A) \circ (\xi \tensor \theta) \simeq (\id_A \tensor \id_A) \circ \xi $$
  which implies $\xi \simeq \theta$. Thus the copairing to a nondegenerate pairing is unique. Now inspecting $\beta := \epsilon \circ \mu$ and $\theta := \delta \circ \eta$ we can immediately verify that
  \begin{align*}
    (\beta \tensor \id_A) \circ (\id_A \tensor \theta) &= ((\epsilon \circ \mu) \tensor \id_A) \circ (\id_A \tensor (\delta \circ \eta)) \\
    &= (\epsilon \tensor \id_A) \circ (\mu \tensor \id_A) \circ (\id_A \tensor \delta) \circ (\id_A \tensor \eta) \\
    &= (\epsilon \tensor \id_A) \circ \delta \circ \mu \circ (\id_A \tensor \eta) \simeq \id_A
  \end{align*}
  where the last isomorphy was derived from the fact that $\triple{A}{\mu}{\eta}$ is a monoida and $\triple{A}{\delta}{\epsilon}$ is a comonoid. This proves the statement.
\end{proof}
\end{lemma}

On the other hand, given any pairing $\beta$ as in \ref{lemma:PairingFrob} determines \textbf{a} possible counit morphism by defining $\epsilon := \beta(e, \ins ) \colon A \lra e$. Thus we could equivalently give $\pair{A, \mu, \eta }{\delta}$ together with a morphism $\beta \colon A \otimes A \lra e$ as above. The analogous result holds for $\theta$ and the unit morphism. Note that for $\beta = \epsilon \circ \mu$ this is the canonical way to retrieve $\epsilon$.\\

Summing up, we can indeed assign a special role to $\beta := \epsilon \circ \mu$ and $\theta := \delta \circ \eta$! They are the unique nondegenerate pairing and copairing that utilize only the given morphisms of the ambient Frobenius Algebra and also such that $\theta$ has the induced morphism $\epsilon := \beta(e, \ins )$ as counit. We will subsequently call them \emph{the nondegenerate pairing} and \emph{the nondegenerate copairing}\index{Frobenius Algebra!unique pairing/copairing}.

\begin{lemma}
\label{lemma:ComultiplicationUnique}
  For a given counit $\epsilon$ of a Frobenius Algebra, $\delta$ is the unique comultiplication\index{Frobenius Algebra!unique comultiplication} such that $\triple{A}{\delta}{\epsilon}$ is a comonoid and such that the Frobenius Relation in \ref{def:FrobeniusCat} commutes.
\begin{proof}
  Let $\omega$ be another comultiplication satisfying all of the above. Thus we inspect $\omega \circ \eta$. Indeed we get
  $$ (\id_A \tensor \epsilon) \circ (\id_A \tensor \mu) \circ (\omega \tensor \id_A) \circ (\eta \tensor \id_A) = (\id_A \tensor \epsilon) \circ \omega \circ \mu \circ (\eta \tensor \id_A)  \simeq \id_A$$
  which qualifies $\omega$ as a copairing for the nondegenerate pairing $\beta$. From Lemma \ref{lemma:CopairingUnique} we know that the copairing is unique and thus $\omega \circ \eta \simeq \theta = \delta \circ \eta$ which shows that $\omega$ is indeed isomorphic to the given comultiplication $\delta$.
\end{proof}
\end{lemma}

Since Frobenius Algebras are defined over monoidal categories, the next natural step is to investigate their behaviour when defined over \emph{braided monoidal categories} \ref{def:monoidal_braided} and \emph{symmetric monoidal categories} \ref{def:monoidal_symmetric}.

\begin{definition}[Symmetric Frobenius Algebras]
\label{def:SymmetricFrobenius}
  Given a braided monoidal category $\triple{\cat{C}}{\tensor}{e, \gamma}$ a Frobenius Algebra $\pair{A, \mu, \eta, \delta}{\epsilon}$ therein is \textbf{symmetric}\index{Frobenius Algebra!symmetric}, if the following diagram commutes:
  \begin{center}
  \begin{tikzcd}[sep=huge]
    A \tensor A \arrow[d, "\mu"'] \arrow[r, "\gamma"] & A \tensor A \arrow[r, "\mu"] & A \arrow[d, "\epsilon"] \\
    A \arrow[rr, "\epsilon"] & & e
  \end{tikzcd}
  \end{center}
\end{definition}

Using the morphism $\beta = \epsilon \circ \mu$ which we investigated in \ref{lemma:PairingFrob} we can sum the above definition up as follows: In a symmetric Frobenius Algebra, the "pairing" $\beta$ is invariant under a twist of its factors. While this seems like a "fair enough" requirement, it has a major problem: There is no a priori reason to assign a preferred role to the multiplication $\mu$ instead of the comultiplication $\eta$. Indeed we find the following insightful result:

\begin{lemma}
\label{lemma:SymmetricComultiplication}
  A Frobenius Algebra $\pair{A, \mu, \eta, \delta}{\epsilon}$ in a braided monoidal category $\triple{\cat{C}}{\tensor}{e, \gamma}$ is symmetric iff $\gamma^{-1} \circ \delta \circ \eta \simeq \delta \circ \eta$.
\begin{proof}
  This prove utilises the uniqueness of the nondegenerate copairing/pairing association proven in \ref{lemma:CopairingUnique}. Thus let the Frobenius Algebra be symmetric and remember that $\beta := \epsilon \circ \mu$ is the nondegenerate pairing with copairing $\theta := \delta \circ \eta$. Also note that
  $$ (\id_A \tensor \gamma) \circ (\id_A \tensor \gamma) = (\id_A \tensor \id_A \tensor \id_A) $$
  Further, since $\gamma$ is natural in both of its factors,
  $$(\id_A \tensor \mu) \circ (\id_A \tensor \gamma) \circ (\gamma^{-1} \tensor \id_A) = \gamma \circ (\mu \tensor \id_A)$$
  as well as
  $$(\id_A \tensor \gamma) \circ (\gamma^{-1} \tensor \id_A) \circ (\delta \tensor  \id_A) = (\id_A \tensor \delta) \circ \gamma$$
  This let's us state the following chain of equalities/isomorphisms:
  \begin{align*}
    &~ \ \ \ (\beta \tensor \id_A) \circ (\id_A \tensor (\gamma^{-1} \circ \theta))\\
    &= (\epsilon \tensor \id_A) \circ (\mu \tensor \id_A) \circ (\id_A \tensor \gamma^{-1}) \circ (\id_A \tensor \delta) \circ (\id_A \tensor \eta) \\
    &= (\epsilon \tensor \id_A) \circ (\mu \tensor \id_A)  \circ (\id_A \tensor \gamma) \circ (\gamma^{-1} \tensor \id_A) \circ (\id_A \tensor \gamma) \circ (\gamma^{-1} \tensor \id_A) \circ (\id_A \tensor \delta) \circ (\id_A \tensor \eta) \\
    &= (\epsilon \tensor \id_A) \circ \gamma \circ \delta \circ \mu \circ \gamma \circ (\id_A \tensor \eta) \simeq (\id_A \tensor \epsilon) \circ \delta \circ \mu \circ (\eta \tensor \id_A) \simeq \id_A
  \end{align*}
  Thus $\gamma^{-1} \circ \theta$ is a possible nondegenerate copairing for $\beta$ which results in $\gamma^{-1} \circ \delta \circ \eta \simeq \theta = \delta \circ \eta$. This proves one direction of the statement. The reverse is just a rephrasing of the above chain of equations/isomorphisms and is thus not done explicitely.
\end{proof}
\end{lemma}

This shows that we could equivalently define a symmetric Frobenius Algebra using the comultiplication. Note that the braiding needed us to resort to the uniqueness of the nondegenerate pairing/copairing instead of that of the comultiplication. We can however restrict to a symmetric monoidal category which yields an even richer algebraic structure to Frobenius Algebras:

\begin{definition}[Commutative Frobenius Algebras]
\label{def:CommutativeFrobenius}
  Given a braided monoidal category $\triple{\cat{C}}{\tensor}{e, \gamma}$ a Frobenius Algebra $\pair{A, \mu, \eta, \delta}{\epsilon}$ therein is \textbf{commutative}\index{Frobenius Algebra!commutative}, if the following diagram commutes:
  \begin{center}
  \begin{tikzcd}[sep=huge]
    A \tensor A \arrow[dr, "\mu"'] \arrow[r, "\gamma"] & A \tensor A \arrow[d, "\mu"] \\
    ~ & A
  \end{tikzcd}
  \end{center}
\end{definition}

\begin{remark}
  It goes without saying that a commutative Frobenius Algebra in particular classifies as a symmetric Frobenius Algebra. Meanwhile the contrary does not need to hold. This manifests in the category of commutative Frobenius Algebras being a full subcategory of the category of symmetric Frobenius Algebras which in turn is a full subcategory of the category of Frobenius Algebra in a symmetric monoidal category. We will define and discuss these categories after a few important points.
\end{remark}

The above definition can be summed up as follows: In a commutative Frobenius Algebra, the multiplication respects the symmetric structure of the ambient symmetric monoidal category. But a priori, there is once again no reason why we should prefer the multiplication over the comultiplication. Indeed we obtain the following result:

\begin{lemma}\label{lemma:FrobAlgComm}
  A Frobenius Algebra $\pair{A, \mu, \eta, \delta}{\epsilon}$ in a symmetric monoidal category $\triple{\cat{C}}{\tensor}{e, \gamma}$ is commutative iff $\delta \circ \gamma = \delta$.
\begin{proof}
  Remember that
  $$ (\id_A \tensor \gamma) \circ (\id_A \tensor \gamma) = (\id_A \tensor \id_A \tensor \id_A) $$
  And further, since $\gamma$ is natural in both of its factors,
  $$(\id_A \tensor \mu) \circ (\gamma \tensor \id_A) \circ (\id_A \tensor \gamma) = \gamma \circ (\mu \tensor \id_A)$$
  as well as
  $$(\gamma \tensor \id_A) \circ (\id_A \tensor \gamma) \circ (\delta \tensor  \id_A) = (\id_A \tensor \delta) \circ \gamma$$
  The main strategy of this proof is to use the symmetry of the Frobenius Algebra to show that $\gamma \circ \delta$ behaves like a comultiplication for the counit $\epsilon$ which allows us to invoke \ref{lemma:ComultiplicationUnique} to state that $\gamma \circ \delta \simeq \delta$. The converse statement, namely the same for the multiplication $\mu$, follows analogously and is thus omitted. Since $\gamma$ is an isomorphism, we immediately see that $\triple{A}{\gamma \circ \delta}{\epsilon}$ is in fact a comonoid. All we need to show is that it additionally satisfies the Frobenius Relation, hence we ultimately need to show
  $$ (\id_A \tensor \mu) \circ (\gamma \tensor \id_A) \circ (\delta \tensor \id_A) = (\gamma \circ \delta) \circ \mu $$
  Thus observe
  \begin{align*}
    &~ \ \ \ (\id_A \tensor \mu) \circ (\gamma \tensor \id_A) \circ (\delta \tensor \id_A) \\
    &= (\id_A \tensor \mu) \circ (\id_A \tensor \gamma) \circ (\gamma \tensor \id_A) \circ (\id_A \tensor \gamma) \circ (\id_A \tensor \gamma) \circ (\delta \tensor \id_A) \\
    &= (\id_A \tensor \mu) \circ (\gamma \tensor \id_A) \circ (\id_A \tensor \gamma) \circ (\gamma \tensor \id_A) \circ (\id_A \tensor \gamma) \circ (\delta \tensor \id_A) \\
    &= \gamma \circ (\mu \tensor \id_A) \circ (\id_A \tensor \delta) \circ \gamma = \gamma \circ \delta \circ \mu \circ \gamma = \gamma \circ \delta \circ \mu
  \end{align*}
  In the first step, we used the commutativity of $\mu$ and inserted a double twist equaling an identity. In the second and third step the naturality of the braiding was used to permute them and then commute with $\mu$ and $\delta$ respectively. The fourth step used the Frobenius Relation diagram to commute $\mu$ and $\delta$ and the last step is just the definition of commutativity of $\mu$. Thus we have shown that $\gamma \circ \delta$ indeed behaves like a comultiplication with counit $\epsilon$ making it isomorphic to $\delta$ itself which proves the statement. Again note that the reverse follows completely analogous.
\end{proof}
\end{lemma}

The last two definitions mainly introduced specialized Frobenius Algebras "respecting" the braided or symmetric structure of the ambient monoidal category. The respective lemmas we proved can be seen as "coherence" results, essentially stating that both definitions do not discriminate between the monoid and the comonoid structure of a Frobenius Algebra. Next we finally turn to a category of Frobenius Algebras the reason for the delay being the preceding definitions of symmetic and commutative Frobenius Algebras which now allow for a concurrent discussion of all three arising categories. With the objects thoroughly introduced, we define morphisms between them:

\begin{definition}[Frobenius Algebra Morphisms]
\label{def:FrobAlgMorph}
  Given two Frobenius Algebras $\pair{A, \mu, \eta, \delta}{\epsilon}$ and $\pair{A^\prime, \mu^\prime, \nu^\prime, \delta^\prime}{\epsilon^\prime}$ in a symmetric monoidal category $\triple{\cat{C}}{\tensor}{e, \gamma}$ a \textbf{Frobenius Algebra morphism}\index{Frobenius Algebra!morphisms} between them is a morphism
  $f \colon A \lra A^\prime$ that is both a monoid morphism (see \ref{def:monoidal_monoid_arrows}) $f \colon \pair{A, \mu }{\nu} \lra \pair{A^\prime, \mu^\prime }{\nu^\prime}$ and a comonoid morphism $f \colon \pair{A, \delta}{\epsilon} \lra \pair{A^\prime, \delta^\prime}{\epsilon^\prime}$.
\end{definition}

Now we \emph{could} easily define respective versions of Frobenius Algebra morphisms for their commutative and/or symmetric delegates but there is no real merrit in doing so for our further investigations. We can be content with knowing that the symmetric morphisms respect and preserve the braiding structure while the commutative ones do so with the symmetric braiding. Ultimately we have already achieved our interim goal in the form of the following definition:

\begin{definition}[Category of Frobenius Algebras]
\label{def:FrobAlgCat}
  Given a symmetric monoidal category $\triple{\cat{C}}{\tensor}{e, \gamma}$ we define the \textbf{Category of Frobenius Algebras}\index{Frobenius Algebra!category}\index{The category of!Frobenius Algebras} of $\cat{C}$ to be the category with Frobenius Algebras in $\cat{C}$ as its objects and Frobenius Algebra morphisms between them as arrows. Composition is given by regular composition of arrows in $\cat{C}$. We denote this category by $\cat{Frob}_\cat{C}$. Further denote the full subcategory of symmetric Frobenius Algebras by $\symfrob_\cat{C}$\index{Frobenius Algebra!symmetric subcategory}\index{The category of!symmetric Frobenius Algebras} and the full subcategory of commutative Frobenius Algebras by $\comfrob_\cat{C}$\index{Frobenius Algebra!commutative subcategory}\index{The category of!commutative Frobenius Algebras}.
\end{definition}
~

Of course the above definition is not general since we demanded a \emph{symmetric monoidal category} while we might just ask for a monoidal category. Given that the existence of even a braiding is not guaranteed here, $\symfrob$ and $\comfrob$ might as well be empty. However since we ultimately work with symmetric monoidal categories like $\cobn{n}$ or $\vectk$, this generalisation does not really provide extra insights. Far more thrilling is the following result about the descend of a braided, symmetric or just purely monoidal structure of the ambient category to the category of Frobenius Algebras:

\begin{corollary}
\label{corollary:FrobCatSym}
  It is an interesting albeit cumbersome task to show that $\triple{\cat{Frob}_\cat{C}}{\tensor}{e}$\index{The monoidal category of!Frobenius Algebras} is in fact itself a symmetric monoidal category with the twist given by the twist of factors of the product. Moreover $\triple{\symfrob_\cat{C}}{\tensor}{e}$\index{The braided monoidal category of!symmetric Frobenius Algebras} is a braided monoidal category and $\triple{\comfrob_\cat{C}}{\tensor}{e}$\index{The symmetric monoidal category of!commutative Frobenius Algebras} a symmetric monoidal category.
\end{corollary}

To verify all three of the above statements would require an excessive amount of tedious writing and is much more efficiently conducted in the pictorial representation akin to that of cobordisms. For a proof of the first statement see \cite[2.4.6 – 2.4.8]{FrobAlgebraTQFT}. We will comment on the pictorial representation used throughout his exposition of the matter in the next chapter.\\

Note that throughout this chapter any examples, especially those relating to $\vectk$, were deliberately avoided whenever not strictly needed for further application. While this might be debatable from a structural perspective, it enables the following workflow: \ref{subsec:AlgebraicPreliminaries} introduced fairly general notions from linear algebra leaning towards a more categorial language as the chapter progressed. \ref{subsec:Frob} then presented Frobenius Algebras purely formulated within category theory while referencing familiar terms of the previous chapter to ultimately foreshadow the inevitable superposition of the two. Unlike custom practise in literature, see \cite[chapter 2.5]{FrobAlgebraTQFT} or \cite[chapter 3.3]{Intro_TQFT}, the superposition in the following chapter can peer through the kaleidoscope of category theory rather than the seasoned binoculars of linear algebra.\\

With the category of Frobenius Algebras at hand, we conclude the top-down introduction to the categorial concept of Frobenius Algebras. The next subsection is dedicated to and an amalgamation of the results of \ref{subsec:AlgebraicPreliminaries} and $\ref{subsec:Frob}$, showing that many of the usual algebraic properties of Frobenius Algebras over a field $\kk$ are in fact just manifestations of the respective general concept from category theory.

\newpage
\subsection{The Category $\comk$}
\label{subsec:CatComfrob}
This entire chapter can be seen as an application of chapter \ref{subsec:Frob} to the algebraic preliminaries given in \ref{subsec:AlgebraicPreliminaries}. Our focus will be the emergence of Frobenius Algebras in the category $\vectk$, a notion that will turn out to have deep ties to $2D$-TQFTs in \ref{sec:TQFT}.
We have already seen that for a symmetric monoidal category $\triple{\cat{C}}{\tensor}{e}$ the category $\comfrob_\cat{C}$ is symmetric monoidal as well. Thus it belongs to the same type of category as $\cobn{n}$, $\vectk$ and in particular $\cobtwo$. The latter relation will serve as the common background for identifying $2$-dimensional TQFTs and the category of commutative Frobenius Algebras in $\vectk$ denoted by $\comk$. Tacitly reading over this major spoiler of the results of this chapter and the next section, we turn towards a rather evident definition of Frobenius Algebras in $\vectk$:

\begin{definition}[Frobenius Algebras over $\kk$]
\label{def:FrobField}
  We call a Frobenius Algebra in the category of vector spaces over a field $\kk$, namely $\vectk$, a \textbf{Frobenius Algebra over a field}\index{Frobenius Algebra!over a field} $\kk$. In a slight abuse of notation we denote the category of these algebras by $\cat{Frob}_\kk$. Its symmetric and commutative subcategories are denoted by $\symfrob_\kk$ and $\comk$ respectively.
\end{definition}

Let us unpack this definition: First of all note that a Frobenius Algebra $\pair{A, \mu, \eta, \delta}{\epsilon}$ forms a monoid $\pair{A, \mu }{\eta}$ and a comonoid $\pair{A, \delta}{\epsilon}$. Thus $\pair{A, \mu }{\eta}$ defines a $\kk$-Algebra and $\pair{A, \delta}{\epsilon}$ a $\kk$-Coalgebra \ref{def:kAlgebraCoalgebraCat}.

This already has some interesting consequences. By rule of \ref{lemma:CopairingUnique} there exists a nondegenerate pairing $\beta := \epsilon \circ \mu$ with nondegenerate copairing $\theta := \delta \circ \eta$ where they. Since the categorial concept of a nondegenerate pairing/copairing \ref{def:NdegPairingFrob} is generalisation of the definition of a nondegenerate pairing for vector spaces \ref{def:NdegPairing}, we can use \ref{lemma:NdegPairingFinite} to conclude:

\begin{corollary}
  The underlying vector space $A$ of any Frobenius Algebra in $\vectk$ is finite-dimensional.
\end{corollary}

Further note that the nondegenerate pairing and coparing are by definition of the categorial term \ref{lemma:PairingFrob} in particular associative \ref{def:AssociativePairings} where we interpreted $A$ as a left and right $A$-module. From here on, we will state different yet equivalent definitions of Frobenius Algebras in $\vectk$ taken from \cite[chapter 2]{FrobAlgebraTQFT} to show how they \emph{and} their equivalence directly emerge from their categorial classification and the many proofs given in \ref{subsec:Frob}. The first and easiest one is the following:

\begin{definition}[Found in {{\cite[3.7]{Intro_TQFT}, \cite[2.3.24]{FrobAlgebraTQFT}}}]
\label{def:FrobKockFirst}
  A Frobenius Algebra over $\kk$ is a $\kk$-vector space $A$ with maps $\mu,  \eta, \delta, \epsilon$ such that $\triple{A}{\mu}{\eta}$ is a $\kk$-Algebra, $\triple{A}{\delta}{\epsilon}$ is a $\kk$-Coalgebra and such that
  $$ (\mu \tensor \id_A) \circ (\id_A \tensor \delta) = \delta \circ \mu = (\id_A \tensor \mu) \circ (\delta \tensor \id_A) $$
\end{definition}

This definition, respectively uniquely defining proposition in \cite[2.3.22 – 2.3.24]{FrobAlgebraTQFT}, is nothing but a \emph{literal} translation of the previously given categorial definition \ref{def:FrobAlgCat} into a concrete version of \ref{def:FrobField}. Since there is not much to conduct from here on, let us investigate the next definition:

\begin{definition}[Found in {{\cite[2.2.5]{FrobAlgebraTQFT}}}]
\label{def:FrobKockSecond}
  A Frobenius Algebra is a $\kk$-Algebra of finite dimension, equipped with an associative nondegenerate pairing $\beta \colon A \tensor A \lra \kk$.
\end{definition}

Again we can quickly recover this definition from our categorial one: We already stated that there exists the associative nondegenerate pairing $\beta := \epsilon \circ \mu$ with copairing $\theta := \delta \circ \eta$ for any Frobenius Algebra stemming from \ref{def:FrobField}. Now going the reverse way we know that by rule of \ref{lemma:CopairingUnique} the copairing for $\beta$ is unique. Since further $\beta$ uniquely determines the counit $\epsilon$ for its unique nondegenerate copairing $\theta$, we can invoke  \ref{lemma:ComultiplicationUnique} to see that we retrieve the full structure of a Frobenius Algebra in the categorial sense. The last definition we will touch here is again a bit easier.

\begin{definition}[Found in {{\cite[2.2.6]{FrobAlgebraTQFT}}}]
\label{def:FrobKockThird}
  A Frobenius Algebra is a ﬁnite-dimensional $\kk$-algebra A equipped with a left A-isomorphism to its dual. Alternatively (and equivalently) $A$ is equipped with a right $A$-isomorphism to its dual.
\end{definition}

From \ref{lemma:NdegPairingIso} we know that these requirements are not only equivalent to each other but also equivalent to the existence of a nondegenerate pairing $\beta$. Since the isomorphisms are left or right $A$-linear respectively, we can use \ref{lemma:PairingAssociativeLinear} to see that $\beta$ is associative. This brings us back to \ref{def:FrobKockSecond} and thus the definition leads to the categorial version. Since the categorial version in particular implies the existence of a associative nondegenerate pairing $\beta$ for which, using \ref{lemma:NdegPairingIso}, we obtain two equivalent isomorphisms as requested in \ref{def:FrobKockThird}, the two definitions are equivalent.\\

Next we turn to the precise classification of subcategories of $\cat{Frob}_\kk$. To this end remember that $\triple{\vectk}{\tensor}{\kk}$ is in fact a \emph{symmetric} monoidal category. Thus we can consider the full subcategories $\symfrob_\kk$ and $\comk$ which are braided/symmetric monoidal categories in their own right. Further investigating the latter let's us draw a baffling conclusion that inspired the very idea to compare $2D$-TQFTs and $\comk$. We will exhibit this discussion by sketching a construction rather than in a strictly formulated paragraph to allow for some historical and intuitive digressions.

\begin{construction}
\label{construction:comFrobToCob2}\index{Frobenius Algebra!pictorial representation}
  As previously mentioned, the author of \cite{FrobAlgebraTQFT} chooses a very different ansatz to discuss Frobenius Algebras. Not only does he start in the classical setting of linear algebra to define the category $\comk$ in the very end (compare \cite[chapter 2]{FrobAlgebraTQFT}), he also provides a purely pictorial representation of many proofs and concepts. This can be seen starting \cite[2.3.5]{FrobAlgebraTQFT}. While the given form of representation would have failed to support many of the interesting and more general proofs of \ref{subsec:Frob}, it is sufficient to provide the specialised variants for Frobenius Algebras in $\vectk$.\\

  While the intention is not to rewrite all or any of these proofs, there is another striking reason to investigate the pictorial representation of the morphisms of a Frobenius Algebra in $\comk$: It looks exactly like the set of generators of $\cobtwo$ we introduced in the classification theorem \ref{theorem:cob2_classification}! Of course this is no mere coincidence but a deep rooted result which we will come to discuss in detail in \ref{subsec:2DTQFT}. Historically the likeness of $2D$-TQFTs and $\comk$ was first discovered by the theoretical physicist \textit{Robbert Dijkgraaf} in his Ph.D. thesis \cite{DijkgraafOne} and later formally proven, albeit without strict treatment of symmetric structures, by \textit{Abrams} \cite{Abrams} and \textit{Quinn} \cite{Quinn}.
  A large part of the foundation of the proof presented in \ref{subsec:2DTQFT} will be the matching pictorial representations of $\cobtwo$ and $\comk$ providing more than enough reason to investigate its roots right here.\\

  Without further ado, let us give a "dictionary" for the morphisms of a Frobenius Algebra $\pair{A, \mu, \eta, \delta}{\epsilon}$ within a symmetric monoidal category $\pair{\cat{C}, \tensor, e}{\gamma}$:

  \begin{center}
  \begin{tikzpicture}[tqft/.cd, cobordism/.style={draw},every lower boundary component/.style={draw},every incoming lower boundary component/.style={draw, dashed},every outgoing lower boundary component/.style={draw, solid}]

    \pic[tqft/cup, rotate=90, at={(2,.5)}];
    \pic[tqft/cap, rotate=90, at={(2,.5)}];
    \pic[tqft/cylinder, rotate=90, at={(5,.5)}];
    \pic[tqft/reverse pair of pants, rotate=90, at={(8,-.5)}];
    \pic[tqft/pair of pants, rotate=90, at={(11,.5)}];
    \pic[tqft/cylinder to next, rotate=90, at={(14,0)}];
    \pic[tqft/cylinder to prior, rotate=90, at={(14,1)}];

    \draw[] (2.1,-.65) node[]{$\updownarrow$};
    \draw[] (3.85,-.65) node[]{$\updownarrow$};
    \draw[] (6.,-.65) node[]{$\updownarrow$};
    \draw[] (9.,-.9) node[]{$\updownarrow$};
    \draw[] (12.,-.9) node[]{$\updownarrow$};
    \draw[] (15.,-.75) node[]{$\updownarrow$};

    \draw[] (2.1,-1.5) node[]{$\epsilon$};
    \draw[] (3.85,-1.5) node[]{$\eta$};
    \draw[] (6.,-1.5) node[]{$\id_A$};
    \draw[] (9.,-1.5) node[]{$\mu$};
    \draw[] (12.,-1.5) node[]{$\delta$};
    \draw[] (15.,-1.5) node[]{$\gamma$};

  \end{tikzpicture}
  \end{center}

  Thus we assign to the

  \begin{center}
  \begin{tabular}{ l l l l l }
    \textbf{unit} & $\epsilon \colon$ & $ A \lra \kk$ & the \textbf{cup} generator &  $1 \lra 0$\\
    \textbf{counit} & $\eta \colon$ & $ \kk \lra A$ & the \textbf{cocup} generator &  $0 \lra 1$\\
    \textbf{identity} & $\id_A \colon$ & $ A \lra A$ & the \textbf{identity} generator &  $1 \lra 1$\\
    \textbf{multiplication} & $\mu \colon$ & $ A \tensor A \lra A$ & the \textbf{copair of pants} generator &  $2 \lra 1$\\
    \textbf{comultiplication} & $\delta \colon$ & $ A \lra A \tensor A$ & the \textbf{pair of pants} generator &  $1 \lra 2$\\
    \textbf{twist} & $\gamma \colon$ & $ A \tensor A \atob{\sim} A \tensor A$ & the \textbf{twist} generator & $2 \atob{\sim} 2$
  \end{tabular}
  \end{center}

  The above comparison immediately shows that the assignment makes sense in terms of in- and outgoing objects. Where the cobordisms have units $0$, the morphisms of the Algebra have the unit $\kk$ and every instance of $A$ is matched by exactly one copy of $S^1$ at the cobordism level. While we will get to a strict formulation of this obvious link between the two structures in the next chapter, let us for now look at the possible applications of this "pictionary" as \textit{Kock} calls it \cite[2.3.7]{FrobAlgebraTQFT}. For example we can now give a pretty representation of the nondegenerate pairing $\beta$:
  \begin{center}
  \begin{tikzpicture}[tqft/.cd, cobordism/.style={draw},every lower boundary component/.style={draw},every incoming lower boundary component/.style={draw, dashed},every outgoing lower boundary component/.style={draw, solid}]

    \pic[tqft/cup, rotate=90, at={(2,.5)}];
    \pic[tqft/reverse pair of pants, rotate=90, at={(0,-.5)}];
    \pic[tqft, incoming boundary components = 2, outgoing boundary components = 0, rotate = 90, at={(3.75,-.5)}];

    \draw[] (-1.5,.5) node[]{$\beta = \epsilon \circ \mu =$};
    \draw[] (3,.5) node[]{$\simeq$};

  \end{tikzpicture}
  \end{center}

  Of course this is not all we can do in this representation. While redrawing completed proofs would be a rather mundane task, we still have an open question regarding the name of the Frobenius Relations for cobordisms \eqref{eq:frobenius_relation}. Looking at the suspiciously similar named equations for Frobenius Algebras \eqref{eq:FrobeniusRelationFrob} we come to the following conclusion:

  \begin{center}
  \begin{tikzpicture}[scale=.5, every tqft/.style={transform shape},tqft/.cd, cobordism/.style={draw},every lower boundary component/.style={draw},every incoming lower boundary component/.style={draw, dashed},every outgoing lower boundary component/.style={draw, solid}]

    \draw[] (8.,-4.5) node[]{$(\id_A \tensor \mu) \circ (\delta \tensor \id_A) = \delta \circ \mu = (\mu \tensor \id_A) \circ (\id_A \tensor \delta)$};
    \draw[] (8.,-3.) node[]{$\Longleftrightarrow$};

    \pic[tqft/reverse pair of pants, name=A, rotate=90, at={(2,-2)}];
    \pic[tqft/cylinder to prior, name=B, rotate=90, at={(2,2)}];
    \pic[tqft/pair of pants, anchor=outgoing boundary 1, rotate=90, at=(A-incoming boundary 2)];
    \pic[tqft/cylinder to prior, anchor=outgoing boundary 1, at=(A-incoming boundary 1), rotate=90];

    \draw (5,0) node[xshift=0.0cm, yshift=0.0cm] {$=$};

    \pic[tqft/reverse pair of pants, rotate=90, at={(6,-1)}];
    \pic[tqft/pair of pants, rotate=90, at={(8,0)}];

    \draw (11,0) node[xshift=0.0cm, yshift=0.0cm] {$=$};

    \pic[tqft/pair of pants, rotate=90, name=C, at={(12,-1)}];
    \pic[tqft/cylinder to next, name=D, rotate=90, at={(12,1)}];
    \pic[tqft/reverse pair of pants, rotate=90, anchor=incoming boundary 1, at=(C-outgoing boundary 2)];
    \pic[tqft/cylinder to next, rotate=90, anchor=incoming boundary 1, at=(C-outgoing boundary 1)];

  \end{tikzpicture}
  \end{center}

  Thus the name "Frobenius Relations" in \eqref{eq:frobenius_relation} is indeed justified. It is but a pictorial representation of the same relations for Frobenius Algebras. Note that the other relations of course have their counterparts in Frobenius Algebras as well. The Unit and Counit relations \eqref{eq:unit_counit} are nothing else than the respective commutative diagrams in the definition of monoids \ref{def:monoidal_monoids} and comonoids \ref{def:monoidal_comonoids}.
  The same goes for the Associativity and Coassociativity relations \eqref{eq:asso_coasso}. Meanwhile Identity relations \eqref{eq:identity_relations} are even more fundamental depicting the behaviour of identity morphisms in a category. The symmetric structure of the ambient monoidal category is recovered in form of \emph{commutative} Frobenius Algebras where multiplication and comultiplication are invariant under $\gamma$. The same property is depicted in the Commutativity and Cocomutativity relations \eqref{eq:commu_cocommu}.\\

  All in all we see that there is a complete mapping between the defining structures of the two a priori very different categories. While so far this is just a "hunch" with strong evidence, we will dedicate an entire chapter to its formalisation and proof. It is no far stretch to say that this will pose the main result of the entire project.
\end{construction}

This concludes our chapter on $\comk$ and thus our section on Frobenius Algebras in general. At the same time it ends the background given in this project, we are now finally in a comfortable place to investigate TQFTs with all the tools at hand to derive many of their fasicnating properties and instances.

\newpage

\newpage
\section{Topological Quantum Field Theory}
\label{sec:TQFT}

In this last section we will explore the definition, properties and even some insightful examples of \emph{Topological Quantum Field Theories}. It draws heavily from the three previous ones and, when looking at $2$-dimensional TQFT, unites them all. While the previous sections have been rather devoid of any physical examples and analogies, this one attempts to weave in a physical perspective on the matter thus trying to appeal to both mathematicians and physicists.\\

In the first chapter we will unfold the general ideas of TQFT, go on to give the mathematical definition and follow up with further discussions about its form and properties. In the following chapter all of the constructions and statements accumulated thus far will be incorporated into the discussion of $2$-dimensional TQFTs and ultimately their equivalence to $\comk$. While the previous chapters will contain previously unseen amounts of physical input, the last chapter is dedicated to the discussion of occurences and applications of TQFTs in physics, again with a particular focus on their $2$-dimensional representatives.\\

Due to its various combined goals this chapter draws from quite a few sources. While we will again use \cite{FrobAlgebraTQFT} and increasingly more parts of \cite{Intro_TQFT}, they both work with the axiomatic definition of TQFTs given by \textit{Atiyah} \cite{Atiyah_1} which we will take a closer look at due to its invaluable physical insights. When talking about $2$-dimensional TQFTs, we will again return to the first two sources while keeping a very unique perspective incorporating the categorial version of Frobenius Algebras presented in \ref{subsec:CatComfrob}.

\subsection{Motivation, Definition and properties}
\label{subsec:TQFTdef}

Classical physical theories like Classical Mechanics, General Relativity, Electrodynamics and more generally classical field theories are mainly administered by the theory of differential equations and their geometric formulations. As such they are particularly concerned with local properties of physical systems. This property is even incorporated as an axiom at the heart of all classical theories, the "Principle of Locality". In rather vague terms it can be formulated as follows:

\begin{principle}[Principle of Locality]
  Any influence of one state on another needs to be mediated by a field traveling through the space between the two.
\end{principle}

This is further determined by the "speed of light" denoted by $c$ which is introduced in Special Relativity. It poses a limit for how fast any such field, hence any interaction, can travel, thus strictly limiting the distance from which a state can be influenced within a certain timeframe. With the advent of Quantum Mechanics and the subsequent verification of its principles and their realisation in nature, a seemingly nonlocal theory took its rightful place next to the established classical theories. Quite fascinatingly Quantum Mechanics exhibits a close connection to low-dimensional topology 
thus further establishing its link to global properties of physical systems rather than local ones. According to Atiyah, this is no particular surprise since "both quantum theory and topology are characterized by discrete phenomena emerging from a continuous background" \cite[p. 175]{Atiyah_1}. While one can investigate this connection in far more rigorous terms, we are mainly concerned with the following arising problem: How can we unify classical theories, which are inherently local, with quantum theories exhibiting clearly nonlocal properties?\\

The field addressing this unification is \emph{Quantum Field Theory}. Many of its popular formulations start from a field theoretic framework of classical mechanics, prominently using jet bundles and either infinite-dimensional manifolds or a multisymplectic formalism, 
to then perform a "Quantisation" of the acquired theory. These theories come in many different flavours embodying different approaches and requirements. We however will only be concerned with a particular subset of so-called "functorial" QFTs, namely those that can be formulated as a functor from one category to another:\\

The main idea of functorial Quantum Field Theory is to encapsulate the assignment that such a theory represents in a functorial prescription between certain categories. Since QFTs describe the transition from a classical theory to a quantised version of the same, it is only natural to take some cobordism category as the "base" category for such a prescription; As previously discussed cobordisms can, allowing a vague analogy, be seen as "evolving slices of a spacetime" and thus leave more than enough freedom to formulate classical geometric theories using further restrictions. On the domain of our functor we want to describe a quantum theory, which is commonly done using hilbert spaces. We catch these objects using the rather general category of vector spaces over a field, $\vectk$, where usually $\kk = \CC$ or $\kk = \RR$. This already brings us to the first important remark:

\begin{remark}
  While in this project oriented and compact manifolds were used to describe cobordisms \ref{def:cobordism_oriented}, they could be endowed with various different additional structures. For example they could have a conformal structure, a spin structure, a framing or work with homotopy classes of maps into some classifying space. This results in conformal QFT (CFT), spin, framed or homotopy TQFT \cite[2.3]{Intro_TQFT}. The prefix "topological" stems from the absence of any notion of a metric or a conformal structure. Technically we are thus working with so-called "oriented closed TQFTs" which we will, in a slight abuse of notation, only call "TQFTs".
\end{remark}

Without further ado, let us take a look at the definition of a Topological Quantum Field Theory to then discuss the implications of the different axioms. Note that while the idea of the definition and the "axioms" (i.e. requirements) are taken from \cite{Atiyah_1}, the presented version is much closer to that given in \cite[1.2.23]{FrobAlgebraTQFT} which enables closer physical analogies:

\begin{definition}[Topological Quantum Field Theory]
\label{def:TQFT}
  An $n$-dimensional \textbf{Topological Quantum Field Theory}\index{TQFT} over a field $\kk$ is a prescription $\ZC$ that associates to each closed oriented $(n-1)$-manifold $\Sigma$ a vector space $\ZC(\Sigma) \in \vectk$ and to each $M \in \cobn{n}$ with $M \colon \Sigma_0 \Lra \Sigma_1$ a linear map $\ZC(M) \colon \ZC(\Sigma_0) \lra \ZC(\Sigma_1)$. This prescription is subject to the following "axioms":
  \begin{enumerate}[label=\alph*)]
    \item The image of two equivalent cobordisms $M \cong N$ is the same, namely $\ZC(M) = \ZC(N)$.

    \item Any cylinder cobordism $\Sigma \times I \colon \Sigma \Lra \Sigma$ is mapped to the identity
    $$\ZC(\Sigma \times I) = \id_{\ZC(\Sigma)} \colon \ZC(\Sigma) \lra \ZC(\Sigma)$$

    \item If a cobordism $M \in \cobn{n}$ has a decomposition of the form $M = M_0 M_1$, $\ZC$ acts as follows:
    $$ \ZC(M) = \ZC(M_1) \circ \ZC(M_0) $$

    \item The disjoint union of $(n-1)$-manifolds $\Sigma_0, \Sigma_1$ is sent to the tensor product of their vector spaces, namely
    $$ \ZC\left(\Sigma_0 \coprod \Sigma_1\right) = \ZC(\Sigma_0) \tensor \ZC(\Sigma_1) $$
    In the same manner if $M = M_0 \coprod M_1$ for $M_0, M_1 \in \cobn{n}$, we get
    $$ \ZC(M) = \ZC(M_0) \tensor \ZC(M_1) $$

    \item The empty manifold, namely $\emptyset_{n-1}$, is sent to the unit of $\vectk$, namely $\kk$. Thus $\ZC(\emptyset_{n-1}) = \kk$.
  \end{enumerate}
\end{definition}

Before delving into the physical aspects of the above definition, let us make some purely mathematical observations. The first three properties simply state that $\ZC$ is a functor in the sense of \ref{def:functors}. Together with the last two properties, this functor is promoted to a monoidal functor. Since it additionally respects the symmetric structure of both the base and the domain category, we can reformulate the definition in the following streamlined and modern form:

\begin{definition}[Topological Quantum Field Theory]
\label{def:FunctorTQFT}
  An $n$-dimensional \textbf{TQFT}\index{TQFT!as a functor} over a field $\kk$ is a symmetric monoidal functor between $\cobn{n}$ and $\vectk$
  $$ \ZC \colon \cobn{n} \lra \vectk $$
  Note that both categories denote the respective symmetric monoidal category.
\end{definition}

Now for the physical perspective: First of all we are working with manifolds of fixed dimension, a very suitable inquiry for any physical setting. Property $a)$ together with $b)$ is yet another manifestation of the "topological" nature of our theory since we work only with diffeomorphism classes of cobordisms without any concern about their possible curvature, metric or other structures. They also embody "relativistic invariance" in that they tell us that our physical theory is invariant under the precise form the transition between the observed "spacetime slices" takes.\\

Property $c)$ might seem a little bland at first sight. However it encapsulates locality by stating that the "evolution", if splitable on the level of cobordisms, is also splitable on the level of linear maps. Hence we obtain an "independence on the history" on the quantum level in that the previous evolution does not explicitely influence the next apart from the output slice it creates.\\
Properties $d)$ and $e)$ are directly linked to the quantum properties of our theory. Just like independent physical systems are described by the tensor product of Hilbert spaces, we want our theory to map cobordisms that do not touch, thus their coproduct, to the tensor product of their respective vector spaces under the TQFT. Meanwhile $e)$ hints at the fact that quantum systems are invariant under global (complex) phases. While this is a simple rewriting of the fact that $\kk$ acts as a unit, we can also read it as the "vacuum" $\emptyset_{n-1}$ being mapped to $\kk$ which expresses the fact that the vacuum of a quantum theory is invariant under complex phases.\\

In practise we will mostly work with the more modern definition of a TQFT given in \ref{def:FunctorTQFT} to provide organic and compact proofs and discussions. However keeping the interpretation of the equivalent definiton \ref{def:TQFT} in mind will prove crucial to understand the physical implications of the upcoming properties and specialised examples like $2$-dimensional TQFT. The first result shows why TQFTs are somewhat "manageable" regarding the involved dimensions:

\begin{theorem}
\label{theorem:FiniteImageNdegPairing}
  Given an $n$-dimensional TQFT $\ZC$ over a field $\kk$, any vector space in its image is equipped with a nondegenerate pairing and is thus in particular of finite dimension.
\begin{proof}
  Let $\Sigma \in \cobn{n}$ and denote by $V := \ZC(\Sigma)$ its image under the TQFT. Now take $\Sigma$ with the opposite orientation, denote it by $\overline{\Sigma}$ and its image by $W := \ZC(\overline{\Sigma})$. Using these two we can define a "pairing-like" and a "copairing-like" cobordism
  $$M \colon \Sigma \coprod \overline{\Sigma} \Lra \emptyset_{n-1}, \quad \quad N \colon \emptyset_{n-1} \Lra \overline{\Sigma} \coprod \Sigma$$
  we denote the induced linear maps by $\beta \colon V \tensor W \lra \kk$ and $\gamma \colon \kk \lra W \tensor V$. Now we know that the cylinder cobordism $\Sigma \times I$ is diffeomorphic to the composition cobordism
  $$ \left(\Sigma \times I \coprod N \right) \left(M \coprod \Sigma \times I \right) \ \cong \ \Sigma \times I$$
  Applying the TQFT and using its properties, this leads to
  $$ \id_{\ZC(\Sigma)} = (\beta \tensor \id_{\ZC(\Sigma)}) \circ (\id_{\ZC(\Sigma)} \tensor \gamma ) $$
  which is nothing but the condition for $\beta$ to be a nondegenerate pairing \ref{def:NdegPairing} with copairing $\gamma$. Thus by \ref{lemma:NdegPairingFinite} $V$ and $W$ are finite-dimensional. Since we did not make any particular choice for $\Sigma$, any image vector space of a TQFT is necessarily finite-dimensional and comes with a nondegenerate pairing.
\end{proof}
\end{theorem}

A direct consequence of the above proof is the following corollary resulting from the isomorphism property seen in \ref{lemma:NdegPairingIso}. It can even be found as a defining axiom in Atiyah's definition of a TQFT \cite[p.178]{Atiyah_1}:

\begin{corollary}
  For any $\Sigma \in \cobn{n}$ one has $\ZC(\overline{\Sigma}) = \ZC(\Sigma)^*$.
\end{corollary}

From here on out, let $\ZC$ be an $n$-dimensional TQFT over $\kk$. The above results give us a direct look at the descend of purely topological properties of the involved manifolds to the algebraic properties of the vector spaces in the image of $\ZC$. Further we now know that any TQFT restricts to a functor onto $\finvect$ which simplifies a lot of constructions.

\begin{remark}
  Note that while finite-dimensional vector spaces suffice for some physical applications like the effective theory of Quantum Computation, there are other models that demand infinite-dimensional vector spaces. This can indeed be achieved by digressing from the purely topological QFTs and introduce metric dependencies or other geometric structures thus altering the cobordism category at hand.
\end{remark}

Further looking into the algebraic properties leads us to an interesting conclusion when combining property $b)$ and $e)$ in \ref{def:TQFT}. Namely if $M$ is an $n$-dimensional manifold without boundary and $\ZC$ an $n$-dimensional TQFT over $\kk$, then we can interpret $M$ as the cobordism
$$ M \colon \emptyset_{n-1} \Lra \emptyset_{n-1} $$
and see that its image under $\ZC$ is a linear map $\ZC(M) \colon \kk \lra \kk$, namely a constant! This is yet another reason why TQFTs are interesting for both mathematicians and physicists. Since they produce constants, i.e. topological invariants, of the manifolds at hand, one can use them as a topological tool in mathematics and as an analytic tool in physics. When encountering a highly perturbative QFT, one could work towards its "topological skeleton", i.e. a TQFT, to compute analytic invariants of the theory and match them with experimental or approximate data.\\

Now using the pairing and copairing we can utilise the above construction to prove the following result:

\begin{lemma}
  For any $\Sigma \in \cobn{n}$ one has $\ZC(\Sigma \times S^1) = \dim(\ZC(\Sigma))$.
\begin{proof}
  First note that $\Sigma \times S^1$ does not have a boundary. Using the short construction above we conclude that the induced cobordism has a constant as its image under $\ZC$. To find this constant, we decompose the "donut" $\Sigma \times S^1$ into the "copairing cobordism" $N \colon \emptyset_{n-1} \Lra \overline{\Sigma} \coprod \Sigma$, the twist cobordism $T \colon \overline{\Sigma} \coprod \Sigma \Lra \Sigma \coprod \overline{\Sigma}$ and the "pairing cobordism" $M \colon \Sigma \coprod \overline{\Sigma} \Lra \emptyset_{n-1}$. A pictorial representation in $2$ dimensions is given by:
  \begin{center}
  \begin{tikzpicture}[tqft/.cd, cobordism/.style={draw},every lower boundary component/.style={draw},every incoming lower boundary component/.style={draw, dashed},every outgoing lower boundary component/.style={draw, solid}]

    \draw[] (3.5,.5) node[]{$N$};
    \draw[] (8.5,.5) node[]{$M$};
    \draw[] (6,-1.) node[]{$T$};

    \pic[tqft, incoming boundary components = 2, outgoing boundary components = 0, rotate = 90, at={(7,-.5)}, name=A];
    \pic[tqft, incoming boundary components = 0, outgoing boundary components = 2, rotate = 90, at={(3,-.5)}, name=B];

    \pic[tqft, incoming boundary components = 1, outgoing boundary components = 1, rotate = 90, offset = 1, at={(5,-.5)}];
    \pic[tqft, incoming boundary components = 1, outgoing boundary components = 1, rotate = 90, offset = -1, at={(5,1.5)}];

  \end{tikzpicture}
  \end{center}

  Applying $\ZC$ yields the composite $\beta \circ \gamma$. For finite vector spaces this is nothing else but the trace over the identity map and thus equal to the dimension of $\ZC(\Sigma)$.
\end{proof}
\end{lemma}

Thus we have a concrete tool to calculate the dimension of the involved vector spaces (or Hilbert spaces in physical applications). Another interesting fact is the following statement about the descend of an $n$-dimensional TQFT to lower dimensions:

\begin{lemma}
  Let $\Theta$ be a closed compact oriented manifold of dimension $k$ where $k < n$. Then the functor
  $$ \ZC_{n-k} \colon \cobn{n-k} \lra \vectk, \quad \quad \left(\Sigma_0 \overset{M}{\Lra} \Sigma_1 \right) \lmap \ZC\left(\Sigma_0 \times \Theta \overset{M \times \Theta}{\Lra} \Sigma_1 \times \Theta\right) $$
  is an $(n-k)$-dimensional TQFT over $\kk$\index{TQFT!dimensional reduction}.
\begin{proof}
  Simply note that
  $$ \ZC_{n-k} = \ZC \circ ((\ins) \times \Theta) $$
  Now since both $\ZC$ and $(\ins) \times \Theta$ are symmetric monoidal functors, their composite is too and thus $\ZC_{n-k}$ is indeed an $(n-k)$-dimensional TQFT over $\kk$ by definition \ref{def:FunctorTQFT}.
\end{proof}
\end{lemma}

This result tells us that given a TQFT of any dimension $n$, we can restrict it to a TQFT of lower dimension by using closed compact oriented manifolds. An instance of such behaviour in physics is the compactification process in Bosonic String Theory where one restricts the $(25+1)$-dimensional space by imposing certain restrictions on some of the dimensions. An easy example is the compactification using a circle where one considers $\RR^{24+1} \times S^1$ to observe the implications of the circle compactification on the dynamics in $\RR^{24+1}$.\\

Another example that might come to the physicists mind when looking at $1$-dimensional cobordisms seen in \ref{example:unit_interval} is Quantum Computing. And apart from the usual pictorial representation of Quantum Computing showing a close analogy to TQFTs, there exists an entire field of studies dedicated to harnessing the stability that topological properties of quantum systems display using TQFT. This approach provides an interesting alternative to stabilize the error prone reality of Quantum Computing. A comprehensible and insightful survey on the field of Topological Quantum Computing can be found in \cite{RowellWang}.\\

The next natural step is to form a category of TQFTs. Since TQFTs are themselves functors between categories, we can utilize natural transformations that act as arrows between functors to form such a category:

\begin{definition}[The Category $\tqft$]
\label{def:CategoryTQFT}
  Define the \textbf{category of $n$-dimensional TQFTs over a field $\kk$}\index{The category of!n-dimensional TQFTs over $\kk$}\index{TQFT!category}, denoted by $\tqft$, to be the category whose objects are $n$-dimensional TQFTs over $\kk$ and whose arrows are monoidal natural transformations \ref{def:monoidal_transformation}. Thus the composition of arrows is given by the usual composition of natural transformation.
\end{definition}

In the upcoming discussions, we will usually refer to $\tqft$ to further establish the categorial nature of the topic at hand. Note that this will also enable a comparison between $\tqft$ and other categories, a cornerstone of the upcoming identification $2\cat{TQFT}_\kk \cong \comk$. Thus having set up the motivation, definition and main properties of general TQFTs, we turn our heads towards the interesting case of their $2$-dimensional representatives.

\newpage
\subsection{$2$-dimensional TQFTs}
\label{subsec:2DTQFT}

In this chapter we will at long last investigate the special properties of $2$-dimensional TQFTs. To achieve a well-defined and rigorous statement about the identification of those TQFTs and commutative Frobenius Algebras, extensive use of the previous three sections and many of the proven results will be needed. Note that a first pictorial link between $2$-dimensional TQFTs and Frobenius Algebras was given in \ref{construction:comFrobToCob2}. This will serve as a base for the upcoming constructions leading to a rigorous formulation and proof of the \emph{Equivalence Theorem} which serves as the main result of this project and provides a deep link between two a priori separate constructs.\\

The exposition will start by taking a closer look at elements of $\twoft$. Due to its role as a generator of the skeleton of $\cobtwo$, their action on $S^1$ will be of particular interest and inspire a more general perspective on TQFTs and "generated" monoidal categories. This investigation together with the fundamental properties of TQFTs from the previous chapter will allow for a rigid formulation of the pictorially "plausible" equivalence of categories $2\cat{TQFT}_\kk \cong \comk$ as well as its proof using everything we know about Frobenius Algebras and monoidal categories. However the result will be proven in far greater generality to then restrict to the special case of $\cobtwo$ and $\vectk$.\\

While the main line of thought again follows the ideas of \cite{FrobAlgebraTQFT} and \cite{Intro_TQFT}, the focus lies on a detailed and general proof of the main theorem with many backlinks to previous results and background. In this spirit, the text will deviate from the classical form of the proof/sketch given in the two sources and adopt a description using the general Frobenius Algebras in a monoidal category introduced in the background \ref{subsec:Frob}. This approach will be a less broad version of the postponed generalisation conducted in \cite[chapter 3.6]{FrobAlgebraTQFT} which contains far more extra topics omitted here to shorten the discussion. Altogether using categorial Frobenius Algebras again allows for a more elegant and general description of the very concept of TQFTs to finally restrict to the special applications of $\tqft$ and discuss its implications.

\subsubsection{Free Monoidal Categories over a Frobenius Algebra}

We have already discussed the identification of $\cobtwo$ with its skeleton \ref{nota:CobSkeleton} generated by the circle $S^1$. It can be summarised by noting that any element of $\cobtwo$ is uniquely determined by a positive integer $n \in \NN$, since we can write it as $\coprod^n S^1$. Now take any $\ZC \in \twoft$ and consider its action on such an element. Since $\ZC$ is a symmetric monoidal functor, we get
$$ \ZC \left( \coprod^n S^1 \right) = \bigotimes^n \ZC(S^1) $$
This calls for a special investigation of $\ZC(S^1)$ since understanding it means understanding the image of any element of $\cobtwo$ under any element of $\twoft$.

\begin{lemma}
\label{lemma:S1FrobAlg}
  For any $\ZC \in \twoft$, $\ZC(S^1)$ is a commutative Frobenius Algebra over $\kk$.
\begin{proof}
  First, let us denote $\ZC(S^1)$ by $F$. So far we only know that $F \in \finvect$ and that there exists a nondegenerate pairing induced by the cobordism $1 \coprod 1 \Lra 0$ with a copairing introduced by $0 \Lra 1 \coprod 1$. Using the second definition of a Frobenius Algebra over a field \ref{def:FrobKockSecond}, we see that this is \emph{almost} enough to state that $F$, or up until now really any image of any TQFT in any dimension \ref{theorem:FiniteImageNdegPairing}, is a Frobenius Algebra.\\

  Not particularly enlightening so far. The missing link is the associativity of the pairing and copairing. For $\cobtwo$ we have already shown that the above cobordisms do indeed satisfy the \emph{Associativity and Coassociativity relations} \eqref{eq:asso_coasso} and thus induce an associative/coassociative pairing/copairing respectively. This is where $\cobtwo$ contributes special structure that does not necessarily exist for an arbitrary cobordism category. Together with \ref{lemma:FrobAlgComm} we further see that $F$ is in particular commutative since the pairing and copairing also satisfy the defining relations \eqref{eq:commu_cocommu}.
\end{proof}
\end{lemma}

Now using the above proof together with the fact that $\cobtwo$ is generated by $S^1$ we could slowly embark onto our journey towards the equivalence of $\cobtwo$ and $\comk$. However the above result can also inspire an even more general result. To this end, let us review the proof:\\

First we used the very general fact that every cobordism is mapped onto an element of $\finvect$ and that there exists a nondegenerate pairing and copairing. The first crucial step was the form of the pairing and copairing. Since $\cobtwo$ is generated by $S^1$, the pairing and copairing did not include reverse orientations but only copies of $S^1$ and thus the images only copies of $F$, not its dual. The next crucial step were the relations $\cobtwo$ happens to satisfy. The Associativity/Coassociativity relations induced the associativity of the nondegenerate pairing and thus made $F$ a Frobenius Algebra. The Commutativity/Cocommutativity relations made $F$ into a commutative Frobenius Algebra.\\

Looking at these properties, we can easily construct a more general type of monoidal category displaying the same properties. Using a generalized notion of a TQFT, we can then recover the above result and an even deeper link between them and general Frobenius Algebras (see \ref{subsec:Frob}). The following definition is a rewriting of the one given in \cite[3.6.16]{FrobAlgebraTQFT}:

\begin{definition}[Free monoidal category over a Frobenius Algebra]
\label{def:FreeFrobAlgCat}
  Let $\pair{\chi, \tensor}{0}$ be a monoidal category whose skeleton is generated by $1$ and $\tensor$. Namely all of its objects are given by an integer $n \in \NN$ and are of the form
  $$ n := 1 \tensor 1 \tensor ... \tensor 1 = \bigotimes^n 1 $$
  Further its set of arrows is generated by the following four arrows
  $$ \mu \colon 2 \lra 1, \quad \quad \delta \colon 1 \lra 2, \quad \quad \eta \colon 0 \lra 1, \quad \quad \epsilon \colon 1 \lra 0 $$
  which satisfy the following relations:
  \begin{enumerate}[label=\alph*)]
    \item $ \mu \circ (\id \tensor \eta) = \id = \mu \circ (\eta \tensor \id) $ \hfill (Commutativity)
    \item $ (\id \tensor \epsilon) \circ \delta = \id = (\epsilon \tensor \id) \circ \delta $ \hfill (Cocommutativity)
    \item $ (\id \tensor \mu) \circ (\delta \tensor \id) = \delta \circ \mu = (\mu \tensor \id) \circ (\id \tensor \delta) $ \hfill (Frobenius)
  \end{enumerate}
  Note that this makes $1 \in \chi$ a Frobenius Algebra in the sense of \ref{def:FrobeniusCat}. Accordingly we call $\chi$ a \textbf{free monoidal category over a Frobenius Algebra}\index{Monoidal category!free}.
\end{definition}

This definition sparks some direct analogies. First of all, a subset of the relations in $\cobtwo$ \ref{theorem:cobordism_relations} holds true. Also the object $1 \in \chi$ takes the same role as $S^1$ in $\cobtwo$ since both generate the respective skeleton via the respective tensor product. Furthermore, like in $\cobtwo$, we could give a pictorial representation of the above objects and arrows and repeat the pictorial proof of \ref{lemma:ImplyAssoCoasso} to derive the following more general result:

\begin{corollary}
  The three relations in \ref{def:FreeFrobAlgCat} imply the following two relations:
  $$ \mu \circ (\id \tensor \mu) = \mu \circ (\mu \tensor \id), \quad \quad (\delta \tensor \id) \circ \delta = (\id \tensor \delta) \circ \delta $$
  They obviously correspond to the Associativity/Coassociativity relations \eqref{eq:asso_coasso}.
\end{corollary}

Altogether we the definition recovers quite a few of the important properties that $S^1$ carries in $\cobtwo$. What hasn't been touched just yet, is the notion of a twist morphism. This can be rectified by defining the following specialisation:

\begin{definition}[Free Symmetric Monoidal Categories over a Frobenius Algebra]
\label{def:FreeSymFrobAlgCat}
  Let $\pair{\chi, \tensor}{0}$ be a free monoidal category over a Frobenius Algebra. Let further $\gamma$ be a braiding and element of the generators of the arrow set such that $\pair{\chi, \tensor}{0, \gamma}$ is a symmetric monoidal category and such that the following additional relations hold:
  $$ \mu \circ \gamma = \mu, \quad \quad \quad \quad \delta = \gamma \circ \delta $$
  Then we call $\chi$ a \textbf{free symmetric monoidal category over a Frobenius Algebra}\index{Monoidal category!free symmetric}. Note that in this case $1$ is a commutative Frobenius Algebra.
\end{definition}

Closely inspecting the above definition and comparing it to the skeleton of $\cobtwo$ \ref{nota:CobSkeleton} one can directly see that the two are essentially the same. Namely whenever working with a free symmetric monoidal category over a Frobenius Algebra, one can equivalently work with $\cobtwo$. Note that this establishes an equivalence of categories, which will not be formalised or used in this project.

\subsubsection{Functor Categories}

$\tqft$ is a prime example of a functor category. Namely its objects are symmetric monoidal functors and its arrows monoidal natural transformations between them. Just like the presented generalisation from $\cobtwo$ to free symmetric monoidal categories over a Frobenius Algebra, we can take another step starting with TQFTs.

\begin{definition}[The Category of Symmetric Monoidal Functors]
  Let $\triple{\cat{C}}{\tensor}{e}$ and $\triple{\cat{D}}{\widetilde{\tensor}}{i}$ be two monoidal categories. We define the \textbf{category of monoidal functors}\index{The category of!monoidal functors}\index{Functor!monoidal functor category} as the category whose objects are monoidal functors from $\cat{C}$ to $\cat{D}$ and whose arrows are monoidal natural transformations between them. We denote it by $\mon(\cat{C}, \cat{D})$.\\
  If the two categories further carry symmetric structures, we define the \textbf{category of symmetric monoidal functors}\index{The category of!symmetric monoidal functors}\index{Functor!symmetric monoidal functor category} by taking symmetric monoidal functors from $\cat{C}$ to $\cat{D}$ as objects and monoidal natural transformations as arrows between them.
\end{definition}

An immediate example is the following equivalence:
$$ \tqft \cong \smon(\cobn{n}, \vectk) $$
This in turn inspires the generalisation of the idea of a TQFT; Namely for any two symmetric monoidal categories $\cat{C}, \cat{D}$ the category $\smon(\cat{C}, \cat{D})$. Note however that the category $\tqft$ allows for the previously discussed physical analogies while other such categories might not be of any physical interest.\\

The first interesting result using this generalisation involves the generating element of a free monoidal category over a Frobenius Algebra denoted by $\chi$. As shown for $S^1$ in $\cobtwo$ we can now take $1 \in \chi$ which in particular is a Frobenius Algebra and note that monoidal functors map Frobenius Algebras onto Frobenius Algebras and respect commutativity. This proves the following lemma:

\begin{lemma}
\label{lemma:FrobToFrob}
  Let $\chi$ be a free monoidal category over a Frobenius Algebra $1$ and let $\triple{\cat{C}}{\ttensor}{e}$ be any monoidal category. Then the image of $1$ under any object in $\mon(\chi, \cat{C})$ is a Frobenius Algebra. Should $\chi$ be symmetric, then the image of $1$, now a commutative Frobenius Algebra, is also commutative.
\end{lemma}

Note that the above statement makes \ref{lemma:S1FrobAlg} obsolete. While we could have used the fact that monoidal functors preserve the structure of Frobenius Algebras to shorten said lemma's proof the goal was rather to get insight into the special role $S^1$ takes apart from being a Frobenius Algebra.

\newpage
\subsubsection{Equivalence Theorems}

With these results and constructions at hand, we can finally prove the deep rooted connection between certain functor categories and Frobenius Algebras that is beginning to surface. Without using the symmetric structure, thus falling back to free monoidal categories over a Frobenius Algebra, we can prove the first core result about the equivalence to Frobenius Algebras:

\begin{theorem}[Equivalence Theorem – Monoidal Functor Categories]
\label{theorem:EquivalenceFrob}\index{Functor!monoidal functor category!Equivalence Theorem}
  Let $\chi$ be a free monoidal category over the Frobenius Algebra $1$ and $\triple{\cat{C}}{\tensor}{e}$ any monoidal category. Then there exists a natural isomorphism
  $$ \mon(\chi, \cat{C}) \cong \cat{Frob}_\cat{C} $$
\begin{proof}
  We start by looking at the objects of $\mon(\chi, \cat{C})$:
  \begin{itemize}
    \item[$\Lra$]
      Since the skeleton of $\chi$ is generated by a finite set of objects and arrows, we can completely determine any functor by its values on these. Thus starting with any monoidal functor let $F$ be the image of $1 \in \chi$. We already know that $F$ inherits the structure of a Frobenius Algebra. Since we work with monoidal functors, we further know that the image of $n$ is $\bigotimes^n F$. The first map we take care of is the identity arrow $\id \colon 1 \lra 1$ which is sent to the identity $\id_F \colon F \lra F$.
      Looking at the very suggestively named arrows in \ref{def:FreeFrobAlgCat} the arbitrarily chosen functor maps $\mu$ to $F \tensor F \lra F$, $\delta$ to $F \lra F \tensor F$, $\eta$ to $e \lra F$ and $\epsilon$ to $F \lra e$.
      Using the relations these arrows satisfy in $\chi$ we immediately see that their images satisfy the relations for a monoid/comonoid such that $\pair{F}{\mu, \eta, \delta, \epsilon}$ is a Frobenius Algebra in the sense of \ref{def:FrobAlgCat}. This also ensures uniqueness of the relevant arrows. Again we could have just cited \ref{lemma:FrobToFrob} however the constructed mapping is important for the other direction of the proof and does a better job at capturing the functorial nature of the described equivalence.

    \item[$\Longleftarrow$]
      Now take a Frobenius Algebra in $\cat{C}$ denoted by $\pair{F}{\mu, \eta, \delta, \epsilon}$. We can easily construct a monoidal functor by taking the previously emerged prescription as its definition. By rule of the relations the arrows in $\chi$ satisfy, the thus described mapping is indeed well-defined in that it respects the relations of a Frobenius Algebra.
  \end{itemize}
  Further note that starting with an arbitrary monoidal functor $\cat{M}$, constructing the respective Frobenius Algebra $F$ using the first direction and then using $F$ to form a monoidal functor, we recover $\cat{M}$ by pure definition. Next we need to investigate the arrows of $\mon(\chi, \cat{C})$, namely monoidal natural transformations.\\

  To this end, let $\cat{M}, \cat{N}$ be two monoidal functors in $\mon(\chi, \cat{C})$ and $\tau \colon \cat{M} \nat \cat{N}$ a monoidal natural transformation between them. Let us denote the images of $1 \in \chi$ under $\cat{M}, \cat{N}$ by $F$ and $G$ respectively.
  Since for any element $n \in \chi$ the natural transformation needs to map $\tau_n \colon \bigotimes^n F \lra \bigotimes^n G$, namely it is the n-fold tensor product of the arrow $F \lra G$, we can restrict its investigation to this simple case.
  \begin{itemize}
    \item[$\Lra$]
      Looking at \ref{def:natural_transformation} we see that every induced map $\tau_n$ is compatible with arrows in $\chi$ which are in turn generated by a finite set of arrows. Thus we can translate the definition reduced to the generating maps to the following four diagrams:
      \begin{center}
      \begin{tikzcd}[sep=huge]
        F \tensor F \arrow[d, "\mu_F"'] \arrow[r, "\tau_2"] & G \tensor G \arrow[d, "\mu_G"] & F \arrow[r, "\tau"] & G \\
        F \arrow[r, "\tau_2"] & G & e \arrow[u, "\eta_F"] \arrow[r, leftrightarrow, "="] & e \arrow[u, "\eta_G"']
      \end{tikzcd}
      ~\\
      \begin{tikzcd}[sep=huge]
        F \arrow[d, "\delta_F"'] \arrow[r, "\tau_2"] & G \arrow[d, "\delta_G"] & e \arrow[r, leftrightarrow, "="] & e \\
        F \tensor F \arrow[r, "\tau_2"] & G \tensor G & F \arrow[u, "\epsilon_F"] \arrow[r, "\tau"] & G \arrow[u, "\epsilon_G"']
      \end{tikzcd}
      \end{center}
      Now the first two diagrams can be summarised by saying that $\tau$ is a monoid morphism in the sense of \ref{def:monoidal_monoid_arrows}. The second pair is just the respective statement for comonoids. Thus by rule of definition \ref{def:FrobAlgMorph} $\tau := \tau_1$ classifies as a Frobenius Algebra morphism.

    \item[$\Longleftarrow$]
      Again, the main work has been done in the converse direction. When given a Frobenius Algebra morphism, we simple use its action on the generating arrows to \emph{define} an obviously monoidal natural transformation using the four diagrams. Again this construction is the inverse of the previous by pure definition. Thus when starting with a monoidal natural transformation $\tau$ and using $(\Lra)$ to form a Frobenius Algebra morphism, which we then use to again construct a monoidal natural transformation, we again arrive at $\tau$.
  \end{itemize}

  Altogether we have proved that every object in $\mon(\chi, \cat{C})$ is in direct correspondence to an object in $\cat{Frob}_\cat{C}$ and further that every arrow of $\mon(\chi, \cat{C})$, namely every monoidal natural transformation, is in direct correspondence to an arrow in $\cat{Frob}_\cat{C}$, namely a Frobenius Algebra morphism. This correspondence defines an equivalence of categories in the sense of \ref{def:equivalence_categories} which proves the theorem.
\end{proof}
\end{theorem}

A direct and deep result following from the above theorem is that $\tqft$ is equivalent to a \emph{subcategory} or all of $\cat{Frob}_\kk$. Note that this statement is not really rigid as of now since we have completely omitted braidings i.e. symmetric structures from \ref{theorem:EquivalenceFrob}. This in turn makes it impossible to account for the symmetric structure of TQFTs and precisely pinpoint their counterparts in Frobenius Algebras. Since we introduced free symmetric monoidal catgeories over Frobenius Algebras, let us formulate a specialised variant of the previous theorem:

\begin{theorem}[Equivalence Theorem – Symmetric Monoidal Functor Categories]
\label{theorem:EquivalenceComFrob}\index{Functor!symmetric monoidal functor category!Equivalence Theorem}
  Let $\chi$ be a free symmetric monoidal category over the commutative Frobenius Algebra $1$ with braiding $\gamma$ and $\triple{\cat{C}}{\tensor}{e, \sigma}$ any symmetric monoidal category. Then there exists a natural isomorphism
  $$ \smon(\chi, \cat{C}) \cong \comfrob_\cat{C} $$
\begin{proof}
  All we need to show is that we can extend the constructions given in the proof of \ref{theorem:EquivalenceFrob} to incorporate the symmetric structures of $\chi$ and $\cat{C}$.\\

  First of all, we are now working with \emph{symmetric} monoidal functors. Thus the image of $1$ is a Frobenius Algebra in $\cat{C}$ and further the twist map $\gamma$ is mapped to the twist of the tensor product in $\cat{C}$, namely $\sigma$.
  Using the additional properties involving the twist (see \ref{def:FreeSymFrobAlgCat}) the image of $1$ classifies as a commutative Frobenius Algebra too. Note that this is just a broader explanation of the general fact that symmetric monoidal functors preserve the commutativity of Frobenius Algebras.\\

  On the other hand if we have any commutative Frobenius Algebra in $\cat{C}$ we already know that we can construct a monoidal functor which is well-defined due to the relations the arrows in $\chi$, now a free \emph{symmetric} monoidal category, satisfy. Using the additional relations involving $\gamma$, we further recover the symmetric structure defining a symmetric monoidal functor. Again the two constructions are inverse to each other.\\

  For arrows the additional structure is even easier to handle. Since any symmetric monoidal functor maps twist map to twist map, the diagram we could draw for $\sigma_F$ and $\sigma_G$ trivially commutes. This also completely determines the definition of a monoidal natural transformation given a Frobenius Algebra morphism in $\comfrob_\cat{C}$. Altogether we can incorporate the symmetric structure of the underlying categories into the equivalence theorem \ref{theorem:EquivalenceFrob} by making a transition to commutative Frobenius Algebras. Since we again form an equivalence of categories, this proves the theorem.
\end{proof}
\end{theorem}

This result finally fits the special case of TQFTs! We can simply apply the very general statement to our previous investigations regarding $\tqft$ and state the following theorem serving as the motivation and main result of this project:

\begin{theorem}[Equivalence Theorem – TQFTs]
\label{theorem:Equivalence2TQFT}\index{TQFT!Equivalence Theorem}
  $\twoft \cong \comk$
\end{theorem}

An almost bland looking result after all the groundwork we have done. However it rigorously states the at best "plausible" correspondence introduced in example \ref{construction:comFrobToCob2}. And while the proof indeed used the same "dictionary" as in said construction, the required abstraction from $\cobtwo$ and even TQFTs brought quite a few new challenges. The specialisation to $\twoft$ really begins to shine once we look at the physical implications of the above theorem:\\

Every $2$-dimensional TQFT is uniquely determined by a commutative Frobenius Algebra capturing the topological structure of a Quantum Field Theory. On the other hand, every commutative Frobenius Algebra uniquely determines a $2$-dimensional TQFT. Thus if a theoretical physicist, in his ongoing considerations of certain algebraic relations between say fields, antifields, their pairwise or individual destruction and creation, happens to recover the structure of a commutative Frobenius Algebra over a field $\kk$, he can directly deduce a $2$-dimensional TQFT and thus the topological skeleton of his theory. As discussed before this allows for the calculation of topological invariants thus analytic values which can be invaluable goods in increasingly complex and often perturbative field theories.\\

This result together with its more general parent \eqref{theorem:EquivalenceComFrob} also provides some fairly universal statements about functorial theories with different domain and/or codomain category. Namely every functorial theory comprised of symmetric monoidal functors that has a free symmetric monoidal category over a Frobenius Algebra as its source category is equivalent to the commutative Frobenius Algebras in its domain category. Furthermore, if the domain category should be $\vectk$, as discussed a very reasonable assumption for a QFT, the theory is necessarily equivalent to $\twoft$. This fact is completely invariant under the precise form of the source category as long as it stays a free symmetric monoidal category.\\

Altogether this draws the main conclusion of this chapter and thus the project. The next and last chapter will investigate applications of TQFTs and especially $2$-dimensional ones in physics making use of the deep results proven in this chapter.

\newpage
\subsection{TQFTs in Physics}

In this chapter a small exposition of physical applications of TQFT will be presented. Note that physical quantum field theories, even if they turn out to be TQFTs in the sense of \ref{def:FunctorTQFT}, are usually formulated and treated using action functionals, Lagrangians and observables. Namely the so-called \emph{path integral formalism} is the predominant language for physical quantum field theories. Since this project is mainly focused on TQFTs as functorial theories, the exposition will be of rather qualitative form not treating the well-definedness of the discussed notion. A short and insightful introduction to the motivation of TQFTs from a path integral perspective is given in \cite[chapter 2.1]{Intro_TQFT} which is subsequently compared to the axiomatic definition in \cite[chapter 2.3]{Intro_TQFT}.\\

After a quick note on the relation between the presented TQFT formalism and the usual path integral formulation following the lines of \cite[2.3]{Intro_TQFT}, a brief introduction to \emph{Chern-Simons} and \emph{Dijkgraaf-Witten theories} will be given. It follows the lines of \cite[chapter 3.2]{Valentino} and will somewhat illuminate the role of TQFTs in a physical context.

\subsubsection{TQFT vs. Path Integral Formalism}

Since there are quite a few ambigous terms in the usual path integral formulation, note that the choices in this section are based purely on convenience regarding the comparison to the functorial TQFTs. We will begin with a short overview of the ingredients required for a physical quantum field theory to then confront them with their counterparts or rather emergence in TQFT. This juxtaposition will make extensive use of and in parts repeat the discussion of the position TQFTs and their interpretation have in physics.\\

The first object we need to introduce are the fields and unsurprisingly so since they lend their name to the theories we aim to investigate. Note that we will restrict to one field to simplify notation. Mathematically a field is just a map $\phi \colon M \lra X$ between two (riemannian) manifolds $(M,g)$ and $(X,h)$. Here $M$ takes the role of the spacetime and $X$ that of the target. If we choose $X = \RR$ we call $\phi$ a scalar field. Now classically we would simply consider the action of the field $\phi$, thus a functional of the form

$$ S[\phi] = \int_M \LL (\phi, \partial_\mu \phi) \sqrt{\det(g)} d^nx $$

Here $\LL$ is the so-called \textbf{Lagrangian} or \textbf{Lagrange function} which depends on $x$ via $\phi$ and $g$ is the metric on $M$. Note that the Lagrangian only depends on the field and its first derivative. Now for a quantum field theory, we obviously need to go one step further. Capturing the quantum nature of the theory, we now investigate a formalized weighted sum over \emph{all possible fields on $M$} which takes the form of an integral over a highly non-trivial \textbf{space of fields} $\FF(M)$. This integral will be called a \textbf{partition function}. The spaces of fields are usually infinite-dimensional and thus require silk gloves when it comes to recovering usual notions from geometry and analysis. Since such considerations are far beyond the scope of this project, here is a partition function:

\begin{equation}\label{eq:Partition}\tag{$\diamond$}
  \ZC_M(\Phi) := \int\limits_{\FF(M)\ s.t. \ \phi|_{\partial M} = \Phi} e^{-S[\phi]} D\phi
\end{equation}

Note that as mentioned, we now integrate over the space of fields and thus formalize the idea of a weighted sum with respect to the value of the respective classical action functional. Note further that the measure and the integral itself require extensive work to be made sense of, however we require only the qualitative insights. In the above equation we technically only treat the special partition function for which the fields restrict to a certain field $\Phi$ on the boundaries. This will be used later on. The last object we need to introduce is the Hamiltonian $\HH$. In a quantum theory it governs the time evolution of the physical system and is thus closely related to any system dependence that is not of purely topological nature. Now we are ready to confront the above notions with their role in TQFTs.\\

First note that the letter $\ZC$ for the partition function was not randomly chosen. It stems from the german word "Zustandssumme" which literally means "sum of states". As such it not only embodies the nature of the integral at hand, it also resembles the sign we chose for TQFTs $\ZC \in \tqft$. And indeed the two are closely related. The axioms for a TQFT as written down in \ref{def:TQFT} translate into physically "reasonable" properties for a partition function:

\begin{itemize}\index{TQFT!comparison to path integrals}
  \item The partition function generates all possible correlation functionals. As such it does not only generalize the usual partition function from statistical mechanics, it also encodes the entire information about the theory, just as a TQFT does by pure definition.

  \item Partition functions on not simply connected manifolds can be calculated separately for each component. This is akin to axiom $d)$ for TQFTs. Note that this can also be adapted to partition functions respecting certain boundary conditions where the boundary is not simply connected.

  \item Given a manifold $\Sigma$ axiom $b)$ tells us that the cylinder cobordism $\Sigma \times I$ is mapped to the identity on $\ZC(\Sigma)$. This can be summarised as the "time evolution" being independent of non-topological factors which directly translates to the path integral setting as $\HH = 0$. This pinpoints the topological nature of our theory in the path integral formalism.

  \item By axiom $e)$ any closed oriented $n$-dimensional manifold $M$, interpreted as a cobordism between empty manifolds, is sent to a map $\kk \lra \kk$. Thus we obtain for any such manifold an element of $\kk$. This corresponds to the partition function evaluated on the manifold $M$ producing a (real) number.
\end{itemize}

Altogether we see immediate links between the two formalisms attempting to treat the same objects. And while the strength of the path integral formalism lies in its direct link to geometry and physical applications, the functorial formulation provides powerful tools to classify low-dimensional theories, in some cases even completely as seen in \ref{theorem:Equivalence2TQFT}. In the following chapter we will introduce a popular example of a topological quantum field theory arising from physical motivation and discuss its relation to the above comparison.

\subsubsection{Examples}

In this chapter a brief introduction to Chern-Simons theory and Dijkgraaf-Witten theory will be given to provide some applied context to the extensive background provided in this and the previous sections of the script. Both are famous examples of topological quantum field theories and the latter can even be used to formulate a $2$-dimensional TQFT. Chern-Simons theory could even be seen as a starting point for a motivation of functorial TQFTs, more on this later. The exposition of both theories mainly draws from \cite{Witten89a}, \cite{DijkgraafWitten}, \cite[chapter 3.2]{Valentino} and some insightful comments of \cite[chapter 3]{Atiyah_1}. Note that the exposition is in no way meant to be extensive, it rather serves to qualitatively embedd the acquired results and notions into actual applications.

\subsubsection*{Chern-Simons Theory}

As mentioned above, Chern-Simons Theory can be seen as a motivation for an axiomatic and functorial approach to TQFTs. We will first give a short introduction to the theory and then discuss this statement by arguing that it indeed is a TQFT. The physically most intriguing property of Chern-Simons theory is the following:\\

When considering a "topologic" theory one usually, especially in physics, just thinks of a theory that is constructed and working entirely without fixing a metric. Such a theory describes only those observables, namely physical properties, that are topologic invariants of the underlying manifold (or topological space). In physics this is usually expressed as \textbf{general covariance} which means that all derived physical laws are invariant under \emph{any} possible coordinate transformation, a rather strict demand. General Relativity achieves this feat by choosing a fixed metric and then integrating over all possible metrics to obtain a "mutable" metric that is an evolving variable rather than a fixed property of the theory. This interpretation has obviously had a tremendous impact on the perspective of modern physics; a generally covariant, i.e. topologic theory is almost synonymous with a theory in which the metric is dynamic.\\

This is where Chern-Simons theory comes in. Unlike General Relativity one does \emph{not} fix a metric and then integrate over all possible metrics but rather starts off without any metric-dependent terms to begin with. So apart from being a TQFT, Chern-Simons Theory also "breaks the mold" when it comes to achieving its purely topological nature.\\

First of all, \textbf{Chern-Simons Theory}\index{TQFT!Chern-Simons} is a $3$-dimensional theory in that it works with a closed $3$-dimensional manifold $M$. We further need a simply connected Lie group $G$ with corresponding Lie algebra $\gf$ and denote by $\tr$ an invariant trace map on $\gf$. Now we denote the space of connections on the trivial $G$-principal bundle over $M$ by $\Omega^1(M; \gf)$. This lets us define the following action functional called the \textbf{Chern-Simons functional}:

\begin{align*}
  S^k_M \colon \Omega^1(M; \gf) &\lra \RR \\
  A &\lmap \frac{k}{4 \pi} \int_M \tr\left( A \wedge dA + \frac{2}{3} A \wedge A \wedge A \right)
\end{align*}

Here $k \in \ZZ$ is called the "level" of the action. Looking at the discussion of the path integral formalism, we see that $A$ takes the role of a field and $S_M$ that of a classical action functional. The next step would thus be to "integrate over all connections" resulting in a partition function of the following form:

$$ Z_k(M) \approxeq \int_{\Omega^1(M;\gf) / \sim} e^{ik S^k_M(A)} DA $$

In the above sketch, $\sim$ represents gauge-equivalence and $DA$ is supposed to be a measure on the arising modulo space. Since such a measure can not necessarily be defined, the above equation should not be read as a definition but rather as an idea, hence the use of "$\approxeq$". However an integral of this form should produce an element of $\RR$ which only depends on the topological properties of $M$; a first nod to TQFT. Further note that since we chose a closed manifold, the above statement related to property $e)$ in the definition of a TQFT. However since the above sketch shows drastic flaws regarding well-definedness, one turns to the properties an object of the form of $Z_k(M)$ would need to display in certain situations.\\

First let us consider $M$ compact with boundary $\partial M$ and define for $\alpha \in \Omega^1(M; \gf)$ the subspace of connections that are gauge-equivalent to $\alpha$ on the boundary of $M$, denoted by $\Omega_\alpha^1(M;\gf)$. We can then, very much in the spirit of \eqref{eq:Partition} from the definition of a partition function, restrict to such equivalences and calculate

$$ Z_k(M)(\alpha) \approxeq \int_{\Omega_\alpha^1(M;\gf) / \sim} e^{ik S^k_M(A)} DA $$

This particular equation marks $Z_k(M)$ as a function $\Omega^1(\partial M;\gf) \lra \RR$ which can be read as a vector space element. Thus the partition function has values in objects of $\vectr$. Going even further we recover a form of property $d)$ in the axiomatic definition of a TQFT \ref{def:TQFT}:\\
For a not simply connected manifold $M = M_1 \coprod M_2$ we obtain a tuple of two vectors, one associated to $M_1$ and one to $M_2$. We can further recover property $c)$ by noting that we can "cut" a manifold by restricting the integration and separately integrate on the two parts to then add them together and receive the same result. Property $a)$ is included since we consider only the topological properties of the underlying manifold, thus the theory is indifferent to the representant of the diffeomorphism class of a manifold.\\

Altogether we have a theory displaying close similarities to the axiomatic properties of TQFTs in addition to being purely topological. Nonetheless the strict formulation of Chern-Simons Theory as a functorial TQFT is still not entirely clear. A short remark on contemporary approaches to solving this issue can be found in \cite{Valentino} who also provides further material on the well-definedness of the theory.

\subsubsection*{Dijkgraaf-Witten Theory}

Unlike Chern-Simons theory, this theory will be formulated entirely within the functorial framwork provided in earlier chapters of this section. While some notions and even categories will not be defined in this project, the user is again refered to \cite{DijkgraafWitten} and \cite{Valentino} for further sources. So let us define \textbf{Dijkgraaf-Witten Theory}\index{TQFT!Dijkgraaf-Witten}:\\

Let $G$ be a group with cardinality $c$, $\kk$ a field whose characteristic is coprime to $c$ and $n > 0$ a positive integer. Now given any closed $(n-1)$-dimensional manifold $\Sigma$ we denote its grupoid of $G$-bundles by $BG(\Sigma)$ and define the following prescription:
$$ \Sigma \lra \kk [\pi_0(BG(\Sigma))] =: V_\Sigma $$
Here $\pi_0$ stands for the $0$-th homotopy group and $\kk[a]$ for the smallest ring that contains $\kk$ and $a$. This can be naturally interpreted as a vector space with coefficients in $\kk$, thus an object of $\vectk$. Moreover this vector space is finite due to the compactnes of $\Sigma$ and the fact that $G$ is necessarily finite. If we take two such manifold $\Sigma_1$ and $\Sigma_2$ and consider their disjoint union we further obtain $V_{\Sigma_1 \coprod \Sigma_2} \cong V_{\Sigma_1} \tensor V_{\Sigma_2}$ which immediately relates to axiom $d)$ in the axiomatic definition of TQFTs.\\

Since we are already close to a TQFT, we next consider an $n$-dimensional cobordism $M \colon \Sigma_1 \Lra \Sigma_2$. Using the grupoids of the three manifolds and the restricting inclusion functors $\imath_1^*, \imath_2^*$ we can assign to $M$ a linear map by setting

\begin{align*}
  \phi_M \colon V_{\Sigma_1} &\lra V_{\Sigma_2} \\
  [y] & \lmap \sum\limits_{[x] \ s.t. \ [x] = [\imath_1^* y]} \frac{\imath_2^*[x]}{|\Aut(x)|}
\end{align*}

Note that $x \in BG(M)$ and $\Aut(x)$ denotes its automorphism group. Using the fact that glueing cobordisms together allows for a splitting of their respective grupoids of $G$-bundles, we can state property $c)$ as
$$ \phi_{MN} = \phi_N \circ \phi_M $$
It is evident that this construction does not depend on the representant of diffeomorphism classes, thus satisfying property $a)$, and further that $\Sigma \times I \colon \Sigma \Lra \Sigma$ is mapped to the identity on $V_\Sigma$. If $M$ is closed, thus a cobordism between empty manifolds, we naturally recover $V_\emptyset = \kk$, axiom $e)$ for a TQFT. Altogether we obtain a symmetric monoidal functor
$$ \ZC_G \colon \cobn{n} \lra \vectk $$
using the above assignments as its definition. Since we can construct a TQFT for any dimension using this procedure, note that $2$-dimensional Dijkgraaf-Witten Theory is equivalent to a commutative Frobenius Algebra as proven in \ref{theorem:Equivalence2TQFT}. As a brief nod towards the path integral formulation the theory, note that for a closed $n$-dimensional manifold $M$ we obtain an element of $\kk$ via the partition function
$$ \ZC_G(M) = \sum\limits_{[x] \in BG(M)} \frac{1}{|\Aut(x)|} $$
Note that $BG(M)$ takes the role of the space of fields while the sum can be seen as the integral of $1$ using the gauge-invariant measure $\frac{1}{|\Aut(x)|}$. For $M = S^n$ this immediately yields $\ZC_G(S^n) = \frac{1}{|G|}$.\\

 The main point of this chapter was to draw some connecting lines between physical applications and the mathematical formulations of TQFTs. We have seen that the translation of the path integral formulation to the functorial setting can be quite challenging while certainly rewarding if we look at the clean results of Dijkgraaf-Witten Theory. At the core of the chapter are the discussions in the comparison of functorial and path integral formalism as well as the comments on topological theories when discussing Chern-Simons Theories.

\newpage

\newpage
\section{Outlook}
\label{sec:Outlook}

At this point there won't be yet another recap of the topics presented in the project. After all a good amount of the script was dedicated to just that thus enabling the tight connections between the different sections. Rather we will single out certain points of the project and provide some insights into possible future or contemporary developments in different directions:\\

Propably the first interception point is the more abstract treatment of Frobenius Algebras (often called Frobenius Objects) in \ref{subsec:Frob}. While certainly not the standard approach to deal with Frobenius Algebras in $\vectk$, we provided quite a handful of results that directly mirror common proofs in a more general formulation. This peaked in the simple descend to the special case of $\vectk$ in \ref{subsec:CatComfrob}. The next natural step would be to investigate further properties of Frobenius Algebras in $\vectk$ using the tools of category theory. For instance one could incorporate rigidness into further proofs or consider the descend of left and right actions of monoids in category theory. Inspired from the usual formulation of these actions, one could investigate if in a compact closed category there exists an equivalence of categories between left and right actions.\\

Regarding TQFTs, the next natural step would be to consider even more general formulations. We provided another abstraction step by incorporating tensor categories and proving the first two TQFT equivalence theorems \ref{theorem:EquivalenceFrob} using free (symmetric) monoidal categories over a Frobenius Algebra. One could start here and further investigate this abstraction deviating from the streamlined exposition provided here. While this is done in parts in a bonus part in \cite{FrobAlgebraTQFT} a more extensive discussion of the benefits and properties of these generalisations might be benefitial for the understanding of TQFTs.\\

The last point I want to mention here is the physical application and interpretation of TQFTs. While certainly not without flaws, I tried to focus on a consistent and insightful discussion of the physical interpretation of the axiomatic definition of TQFTs, its general properties and popular examples like Chern-Simons Theory. Another related topic is the interpretation of TQFTs formulated using path integrals to a functorial setting. As mentioned before, many theories, e.g. Chern-Simons Theory, are still not fully understood in a functorial setting. To use the advantages and tools of both formulations to the fullest extent, an ongoing focus on such translations will be mandatory.

\newpage

\newpage
\section{Bibliography}
\label{sec:Biblio}

\addcontentsline{toc}{section}{Index}
\printindex

\end{document}